\documentclass{article}
\usepackage{graphicx} 
\usepackage{tikz}
\usepackage{tikz-cd}
\usepackage{quiver}
\usepackage{amssymb}
\usepackage{amsmath}
\usepackage{amsthm}
\usepackage{mathtools}
\usepackage{bm}
\usepackage{framed,enumitem} 
\usepackage{hyperref}
\usepackage{subfiles}
\usepackage{bbm}
\usepackage{varwidth}
\usepackage{pinlabel,tikz,listings}
\usepackage[all]{xy}
\usepackage{extarrows}
\usepackage{caption}
\usepackage{subcaption}
\usepackage{geometry}
\usepackage{indentfirst}

\hypersetup{colorlinks=true,citecolor=blue,
linkcolor=blue,urlcolor=black,filecolor=black}

\newtheorem{theorem}{Theorem}
\numberwithin{theorem}{section}
\newtheorem*{theoremA*}{\textbf{Main Theorem A}}

\newtheorem*{theoremB*}{\textbf{Main Theorem B}}

\newtheorem*{conjectureC*}{\textbf{Conjecture C}}

\theoremstyle{definition}
\newtheorem{definition}[theorem]{Definition}

\newtheorem{proposition}[theorem]{Proposition}

\newtheorem{lemma}[theorem]{Lemma}

\newtheorem{remark}[theorem]{Remark}

\newtheorem{example}[theorem]{Example}
\newtheorem{algorithm}[theorem]{Algorithm}

\newtheorem{conjecture}[theorem]{Conjecture}

\newcommand{\BP}{\operatorname{BP}}

\newtheorem{corollary}[theorem]{Corollary}

\newtheorem*{notation}{Notation}

\newcommand{\N}{\mathbb{N}}
\newcommand{\Z}{\mathbb{Z}}

\newcommand{\C}{\mathbb{C}}

\renewcommand{\1}{\mathbf{1}}
\newcommand{\Id}{\mathrm{Id}}

\newcommand{\matfont}[1]{\pmb{\mathsf{#1}}}
\newcommand{\mb}[1]{\mathbf{#1}}

\newcommand{\Addresses}{{
  \bigskip
  \footnotesize

 \textsc{Institute of Mathematics, University of Zurich, Winterthurerstrasse 190, CH-8057 Zurich, Switzerland}\par\nopagebreak
  \textit{E-mail address}: \texttt{quentin.faes@math.uzh.ch}

  \medskip

  \textsc{Institute of Mathematics, University of Zurich, Winterthurerstrasse 190, CH-8057 Zurich, Switzerland}\par\nopagebreak
  \textit{E-mail address}: \texttt{maksymilian.manko@math.uzh.ch}

}}

\title{\normalsize \textbf{NON-FACTORIZABLE RIBBON HOPF ALGEBRAS}}

\author{Quentin Faes and Maksymilian Manko}
\date{\normalsize{May 2025}}

\begin{document}
\maketitle

\begin{abstract}
Building on the work of Nenciu we provide examples of non-factorizable ribbon Hopf algebras, and introduce a stronger notion of non-factorizability. These algebras are designed to provide invariants of $4$-dimensional $2$-handlebodies up to $2$-deformations. We prove that some of the invariants derived from these examples are invariants dependent only on the boundary or on the presentation of the fundamental group of the $2$-handlebody.
\end{abstract}

\tableofcontents
\section{Introduction}
Unimodular, ribbon categories introduced by Turaev in \cite{turaev_94} constitute an important topic in quantum algebra, chiefly for their many applications in quantum topology. Since the 1990s, multiple constructions of 3-manifold and 4-manifold invariants have been constructed using these essential ingredients \cite{turaev_94}, \cite{kerler2001non}, \cite{crane1993categorical}. Ribbon Hopf algebras are a natural source of such categories through their representation theory. A ribbon Hopf algebra $H$ is a quasitriangular Hopf algebra (i.e. it comes with the extra data of an invertible element $R$ of $H\otimes H$ called the R-matrix) together with a central element of $H$ called the ribbon element. Both $R$ and $v$ are required to satisfy certain equations, see Section \ref{secprelim}.

In particular, one may define a $3$-manifold invariant from any so-called \emph{factorizable} unimodular ribbon Hopf algebra. Factorizability is a non-degeneracy condition related to the \emph{monodromy matrix} $M:= R_{21}R$ of the algebra, where $R_{21}$ is the R-matrix with permuted entries (see Proposition \ref{factorizability}). Recall also that a Hopf algebra $H$ is unimodular if it admits an integral $\lambda \in H^*$ and a cointegral $\Lambda \in H$ satisfying certain properties, as detailed in Definition \ref{defunimod}.

There is, similarly, a strong interest in  \textit{triangular} Hopf algebras. By the work of Bobtcheva \cite{bobtcheva_2023} these provide invariants of 2-dimensional CW-complexes, and are closely related to a well-known problem in combinatorial group theory, the \textit{Andrews-Curtis conjecture} \cite{andrews1965free}. Finite-dimensional triangular Hopf algebras over $\C$ were fully classified in \cite{andruskiewitsch_etingof_gelaki_2002}.

More recently there has been a considerable interest in the intermediate case, the unimodular, \textit{non-factorizable} ribbon Hopf algebras, as ingredients of non-semisimple 4-manifold invariants. In \cite{beliakova_derenzi_2023} Beliakova and De Renzi introduced Kerler-Lyubashenko-type invariants of 4-dimensional 2-handlebodies up to 2-deformations \cite{beliakova_derenzi_2023} based on the work of Bobtcheva and Piergallini \cite{bobtcheva_piergallini_2006}, and presented in a more streamlined way in \cite{beliakova_bobtcheva_derenzi_piergallini_2023}. These invariants are the main motivation for this paper. They are meant to provide an angle on a conjecture of Gompf, see Conjecture \ref{conj}, that asserts any diffeomorphic 4-dimensional 2-handlebodies are related by a sequence of Kirby moves involving handles of index at most 2, and is expected to be false. A TQFT construction producing invariants of 4-manifolds up to diffeomorphism introduced in \cite{cghp_2023} similarly relies on representation categories of non-semisimple, unimodular, non-factorizable ribbon Hopf algebras.

At the topological level, non-factorizability is directly related to the ability of the invariants from \cite{beliakova_derenzi_2023} to differentiate between $1$-handles and $2$-handles, which in turn is related to the ability of the invariant to see beyond the boundary of the $4$-dimensional $2$-handlebody. Non-(co)semisimplicity is related to the sensitivity of the invariant to a stabilization operation on handlebodies. This is detailed in Section \ref{sectopology}.

Several examples of non-factorizable unimodular ribbon Hopf algebras are known and include the small quantum group $u_q \mathfrak{sl}_2$ at even root of unity $q$. It was shown in \cite{beliakova2022refined} that the scalar Kerler-Lyubashenko type invariant of \cite{beliakova_derenzi_2023} applied to a connected $4$-dimensional $2$-handlebody only sees its $3$-dimensional boundary and homological information. Another example is the Hopf algebra of symplectic fermions $\operatorname{SF}_{2n}$, which admits a continuous family of non-factorizable ribbon structures, including a triangular one. The generic case, for $n=1$, was shown by Kerler in \cite{kerler2003homology} to give essentially invariants of the $3$-dimensional boundary with this approach. In both cases, the Hopf algebras are, in some sense, ``close to factorizable". Indeed, in the unimodular case, factorizability is equivalent to requiring that applying the integral to one side of the monodromy matrix gives the cointegral. In these non-factorizable cases we mentioned, this equation will still be verified up to an invertible element of the algebra.

\subsection{Main results}
Hence, in an attempt to avoid such ``degeneration to the boundary", in this paper, we introduce the stricter notion of \textit{strong non-factorizability}. We demand that the result of applying the integral to either side of the monodromy matrix, with the antipode previously applied to its first tensor component, is zero, and hence cannot be related to the cointegral by any invertible element of the algebra:
$$\big((\lambda \circ S)\otimes \Id\big) M = \big( S \otimes \lambda \big)M= 0.$$ 

Non-semisimple triangular Hopf algebras provide examples of this property, but yield handlebody invariants only related to the presentation of the fundamental group. Here, we present two constructions capable of producing ribbon Hopf algebras which admit only strongly non-factorizable quasitriangular structures, and are generically not triangular. 

First we discuss a family of Hopf algebras that was introduced in \cite{beattie_dascalescu_grünenfelder_2000} and analysed from the angle of quasitriangularity and ribbonness in \cite{nenciu_2004}. Thus, we will call these Hopf algebras \textit{Nenciu algebras}. Let $\mb{m}\in \Z^s$ be a row tuple of integers of length $s\in \Z$, $t\in \Z$ and $\matfont{d}, \matfont{u}$ be $t\times s$ matrices with integer entries. A Nenciu algebra is generated by grouplike generators $K_a$, $a=1, \dots, s$ such that $K^{m_a}_a=\1$, commuting with each other, as well as skew-primitive nilpotent generators $X_k$, $k=1, \dots, t$ such that $X^2_k=0$. The relations between the generators are prescribed by the matrices $\matfont{d}, \matfont{u}$ and are \textit{diagonal}, that is commuting any two elements results only in a constant in front of the swapped product. We denote the resulting algebra by $H(\mb{m}, t, \matfont{d}, \matfont{u})$, and a full definition is provided in Definition \ref{genDef}. In \cite{nenciu_2004}, an array of previously existing quasitriangular Hopf algebra constructions, including the family $\operatorname{SF}_{2n}$, were unified, and sufficient conditions for such an algebra to be ribbon where determined, but no new examples were given. In particular, the questions of unimodularity and factorizability were not adressed.

In this paper, we exhibit new examples. Moreover we determine sufficient conditions for unimodularity, and for a number of instances we establish when non-factorizable and strongly non-factorizable ribbon structures are admitted, relying in particular on a theorem of Radford \cite{radford_1994}, reproduced below as Theorem \ref{classicalTheorem}. We also study another property: the algebra is called \emph{anomaly-free} if $\lambda(v)$ is invertible and \textit{anomalous} otherwise (see Definition \ref{DefAF}). Topologically, as explained in \cite{beliakova_derenzi_2023}, when the algebra is both factorizable and anomaly-free, the $4$-dimensional $2$-handlebody invariant associated to the algebra is actually an invariant of the boundary of the handlebody. For completeness, we check that anomaly-freeness is not satisifed in the examples where the corresponding quasitriangular structure is strongly non-factorizable. Combining the results of \cite{nenciu_2004} with our own, we have the following.
\begin{theoremA*}
    Under an appropriate choice of parameters a Nenciu type Hopf algebra $H(\mb{m}, t, \matfont{d}, \matfont{u})$ can be simultaneously
    \begin{enumerate}
        \item quasitriangular (Theorem \ref{QTthm} and \cite[Theorem 3.4]{nenciu_2004}),
        \item ribbon (Theorem \ref{ribbonThm} and \cite[Proposition 4.9]{nenciu_2004}),
        \item unimodular (Proposition \ref{unimodTheorem}),
        \item strongly non-factorizable (Proposition \ref{ExAreSNF}),
        \item and anomalous (Corollary \ref{CorNenciuAF}).
    \end{enumerate}
\end{theoremA*}
However, we show in Theorem \ref{2dDegenProp} that the related Kerler-Lyubashenko invariants of 4-dimensional 2-handlebodies degenerate to invariants of the 2-dimensional \textit{spine}, the 2-dimensional CW-complex resulting from collapsing the cocores of the handles attached or equivalently to invariants of the presentation of the fundamental group of the handlebody (see the discussion above the aforementioned proposition).

In an attempt to remedy this, and inspired by the work of Majid \cite{majid_2000}, we present a second family of non-factorizable Hopf algebras where we form a sort of semi-direct product of $u_q \mathfrak{sl}_2$ where $q$ is a root of unity of order divisible by 4, and a Nenciu-type Hopf algebra $H(\mb{m}, t, \matfont{d}, \mb{u})$, where the former (co)acts on the latter on the right. While on the Hopf algebraic level, the construction utilises standard techniques, the fine-tuning required to retrieve the desired properties does not seem to have appeared in the literature before. The action and coaction involve only the $K$ generator of $u_q \mathfrak{sl}_2$ and are designed so that the result is unimodular and carries a ribbon structure being a slightly modified product of the R-matrices and ribbon elements of the two respective factors. It turns out this Hopf algebra has different properties than a trivial tensor product Hopf algebra $u_q \mathfrak{sl}_2 \otimes H(\mb{m}, t, \matfont{d}, \mb{u})$, or standard smash and cosmash products affecting only the algebra or coalgebra (and the antipode), respectively, presented in \cite{molnar_1977}. Thus, we dub it the \textit{semidirect biproduct} as both algebra and coalgebra structures interact, and denote $u_q \mathfrak{sl}_2 \ltimes H(\mb{m}, t, \matfont{d}, \mb{u})$ to emphasize that it is the $u_q \mathfrak{sl}_2$ factor that (co)acts, see Definition \ref{smashBiproduct}. Its key property is that if the Nenciu factor carries a strongly non-factorizable, anomalous ribbon structure, so does the semidirect biproduct with $u_q \mathfrak{sl}_2$. We summarise our results as follows.
\begin{theoremB*}
    Under an appropriate choice of parameters a semidirect biproduct $u_q \mathfrak{sl}_2 \ltimes H(\mb{m}, t, \matfont{d}, \mb{u})$ can be simultaneously
    \begin{enumerate}
        \item quasitriangular (Theorem \ref{extMainThm})
        \item ribbon (Theorem \ref{extMainThm}),
        \item unimodular (Theorem \ref{extMainThm}),
        \item strongly non-factorizable (Proposition \ref{SNFextension}),
        \item and anomalous (Corollary \ref{CorBiprodAF}).
    \end{enumerate}
\end{theoremB*}
However, we again observed that no true 4-dimensional information seems to be carried by the $4$-dimensional $2$-handlebody invariants given by this family of Hopf algebras. More precisely, computations for multiple examples suggest that the invariant might be given by the product of the 2-dimensional invariant of the spine of the handlebody obtained from $H(\mb{m}, t, \matfont{d}, \mb{u})$ and the refined 3-dimensional invariant of its boundary associated to $u_q \mathfrak{sl}_2$. 
Thus, we leave it as Conjecture \ref{failedConjecture}. 
\subsection{Structure of the paper}
The paper is organised as follows: in Section 2 we introduce all the relevant definitions pertaining to unimodular ribbon Hopf algebras and the notation used to describe the Nenciu algebras. In Section 3 the Nenciu algebras are introduced together with the conditions for a given one to admit unimodular, ribbon, and strongly non-factorizable structures. In Section 4, we introduce the second family of examples and once again give conditions for strong non-factorizability. In Section 5 we recall some topological background, including the Kerler-Lyubashenkno invariants of 4-dimensional 2-handlebodies, and discuss the efficacy of the invariants provided by our examples. Finally, in Appendix A we list all the examples of Hopf algebras found throughout the paper, and in Appendix B we collect some proofs that are too long to be presented in the main body.\\

\noindent\textbf{Acknowledgements:} 
The authors would like to thank Ivelina Bobtcheva, Marco de Renzi,  Azat Gainutdinov, Cris Negron and Matthieu Faitg for early comments on the paper and fruitful discussions, as well as Riccardo Piergallini for making his Mathematica computer program available, and Anna Beliakova for her help and introducing them to the topic. Both authors were supported by Simons Collaboration on New Structures in Low-Dimensional Topology
and Grant 200020\_207374 of the Swiss National Science Foundation.

\section{Preliminaries}
\label{secprelim}
Let us recall some definitions and results regarding Hopf algebras (in the category of complex vector spaces), and fix the notation used in the rest of the paper.
\subsection{Notations for Hopf algebras}
In this paper, we consider Hopf algebras over the field of complex numbers $\C$.
\begin{definition}
    A \textit{Hopf algebra} $(H, \mu, \eta, \Delta,\epsilon, S)$ is the data of a vector space (over $\C$) $H$ together with 
    \begin{itemize}
        \item a \textit{unit} $\1: \C \rightarrow H$ and a \textit{product} $\mu: H\otimes H \rightarrow H,$
        \item a \textit{counit} $\epsilon: H \rightarrow \C$, and a \textit{coproduct} $\Delta: H \rightarrow H\otimes H$ 
        \item an invertible \textit{antipode} $S : H\rightarrow H$,
    \end{itemize}
    that are $\C$-linear maps satisfying the axioms
$$
\begin{gathered}
\mu \circ(\mu \otimes \mathrm{id})=\mu \circ(\mathrm{id} \otimes \mu), \\
\mu \circ(\eta \otimes \mathrm{id})=\mathrm{id}=\mu \circ(\mathrm{id} \otimes \eta), \\
(\Delta \otimes \mathrm{id}) \circ \Delta=(\mathrm{id} \otimes \Delta) \circ \Delta, \\
(\varepsilon \otimes \mathrm{id}) \circ \Delta=\mathrm{id}=(\mathrm{id} \otimes \varepsilon) \circ \Delta, \\
(\mu \otimes \mu) \circ(\mathrm{id} \otimes \tau \otimes \mathrm{id}) \circ(\Delta \otimes \Delta)=\Delta \circ \mu, \\
\varepsilon \circ \mu=\varepsilon \otimes \varepsilon, \\
\Delta \circ \eta=\eta \otimes \eta, \\
\varepsilon \circ \eta=1, \\
\mu \circ(S \otimes \mathrm{id}) \circ \Delta=\eta \circ \varepsilon=\mu \circ(\mathrm{id} \otimes S) \circ \Delta, 
\end{gathered}$$
where $\tau:H\otimes H \rightarrow H \otimes H$ is the \textit{flip map} exchanging the factors in the tensor product.\\
We will often abusively denote the Hopf algebra by $H$, and write $gh$ instead of $\mu(g,h)$ for $g,h \in H$. We also note that $H\otimes H$ can be endowed with an algebra structure, using the product
\[
    (a\otimes b)(c\otimes d):= ac \otimes bd. 
\] where $a,b,c,d \in H$.
\end{definition}
\begin{remark}
    Hopf algebras may be defined in any braided monoidal category $\mathcal{C}$. The algebra structure of $H\otimes H$ is then defined using the braiding.
\end{remark}

\begin{notation}
     We use the Sweedler notation for the coproduct: if \[\Delta(h) = \sum_{i =1}^n h_{(1)i} \otimes h_{(2)i} \] then we write
    \[
        \Delta(h) = h_{(1)} \otimes h_{(2)},
    \]
    for some $n\in \Z_{>0}$.
By $\Delta^{cop}:=\tau \circ \Delta$ we denote the \textit{coopposite coproduct}, with $\tau:H\otimes H \rightarrow H \otimes H$ the flip map. Recall that if $H$ is a Hopf algebra, then $H^{cop}:= (H,\mu,\eta,\Delta^{cop}, \epsilon, S)$ is also a Hopf algebra.
\end{notation}Recall that if $H$ is a Hopf algebra, the dual vector space
    \[
        H^*:= \operatorname{Hom}_\C (H, \C).
    \]
    can also be provided with a Hopf algebra structure. 

\begin{definition}
     An element $g\in H$ is called \textit{grouplike} if $\Delta(g)= g\otimes g$. The subset $G(H)$ of all such elements carries a group structure, and generates a Hopf subalgebra which is isomorphic to the group algebra of $G(H)$. An element $X \in H$ is called \textit{skew-primitive with respect to} $L\in G(H)$ if $\Delta(X) = \1 \otimes X + X\otimes L$, and \textit{primitive} if $L=\1$. 
\end{definition}
\noindent
\begin{definition}
\label{adjointaction}
    The \emph{left adjoint action} of $H$ is the left $H$-module structure on $H$ given by
    \[
        \rhd: H\otimes H \rightarrow H, \;\; a\rhd b :=  a_{(1)}\,b \,S(a_{(2)})
    \]
    in the Sweedler notation.
\end{definition}
\subsection{Integrals and unimodularity}

Let $H$ be a finite dimensional Hopf algebra over an algebraically closed field. Then it admits a left (resp. right) integral $\lambda_L \in H^*$ (resp. $\lambda_R$) and cointegral $\Lambda_L\in H$ (resp. $\Lambda_R$), see \cite[Chapter 10.1]{radford_2012}. 
\begin{definition}
\label{defunimod}
    An element $\lambda_L \in H^*$ is a \textit{left} (resp. $\lambda_R \in H^*$ a \textit{right}) \textit{integral} of $H$ if for all $h\in H$,
\begin{align*}
    h_{(1)}\otimes \lambda_L(h_{(2)}) = \lambda_L(h) && \big(\text{resp.} && \lambda_R(h_{(1)}) \otimes h_{(1)}  = \lambda_R(h)\big).
\end{align*} 
\noindent An element $\Lambda_L \in H$ is a left (resp. $\Lambda_R \in H$ a right) cointegral if for all $h\in H$, 
\begin{align*}
     h\Lambda_L = \epsilon(h)\Lambda_L && \big(\text{resp.} && \Lambda_R h = \Lambda_R\epsilon(h)\big).
\end{align*}
An integral (resp. cointegral) is called two-sided if its both a left and right integral (resp. cointegral), or equivalently, if
\begin{align*}
    \lambda \circ S = \lambda && \big(\text{resp.} && S(\Lambda) = \Lambda\big).
\end{align*}
\end{definition}
\noindent It is also known that $\Lambda_L, \Lambda_R$ generate $1$-dimensional subalgebras of $H$ and $\lambda_L, \lambda_R$ generate $1$-dimensional subalgebras of $H^*$, thus different choices differ up to a constant in the ground field.
\begin{definition}[\cite{radford_2012} Definition 10.2.3]
    We will call $H$ \textit{unimodular} if $\Lambda$ is two-sided.
\end{definition}
\subsection{Quasitriangular and ribbon structures}
Let us recall the definition of the quasitriangular Hopf algebra. We need the following notation.

\begin{notation}
    For an element $R\in H\otimes H$ such as the \textit{R-matrix} of a quasitriangular Hopf algebra of Definition \ref{QTdef} below, we write
\[
    R = R' \otimes R'',
\]
where the summation is again implicit. Sometimes, we will need many copies of $R$, in which case they will be labeled with numbers. We also denote $R_{21}:= \tau(R)$.
\end{notation}

\begin{definition}
\label{QTdef}
    Let $H$ be a Hopf algebra and $R$ be an element of $H\otimes H$. The pair $(H,R)$ is quasitriangular if the following axioms are satisfied.
    \begin{center}
    \begin{varwidth}{\textwidth}
        \begin{itemize}
            \item[(QT1)] $\Delta(R') \otimes R'' = R'_1 \otimes R'_2 \otimes R''_1 R''_2$
            \item[(QT2)] $\epsilon(R')R'' = \1$
            \item[(QT3)] $R' \otimes \Delta^{cop} (R'') = R'_1 R'_2 \otimes R''_1 \otimes R''_2$
            \item[(QT4)] $\epsilon(R'')R' = \1$
            \item[(QT5)] $\Delta^{cop} (h) R = R\Delta(h), \forall h \in H$
        \end{itemize}
    \end{varwidth}
    \end{center}
    \noindent Then $R$ is called a  \textit{universal R-matrix} (where we usually drop the ``universal"), and it follows from (QT2) and (QT4) that it is invertible in the algebra $H\otimes H$, with the inverse denoted $R^{-1}$.
\end{definition} 

\begin{definition}
    The \textit{monodromy matrix} is defined as
    \[
        M:= R_{21} R = R''_1 R'_2 \otimes R'_1 R''_2 .
    \] If $M =  \1\otimes \1$, the quasitriangular structure is called \textit{triangular}.
\end{definition}
\noindent 
In order to exhibit quasitriangular structures on Hopf algebras, we will strongly rely on the following result.
\begin{theorem}[\cite{radford_1994}, Section 2.1] \label{QTthm} 
\label{classicalTheorem}
    Let $H$ be a Hopf algebra. Any element $R \in H \otimes H$ satisfying (QT1)-(QT4) induces a Hopf algebra map 

    \begin{align*}
         f_R: H^* &\rightarrow H^{cop} \\ h^* &\mapsto h^*(R')R'',
    \end{align*}
    where $H^{cop}$ is the \textit{coopposite Hopf algebra} of $H$, defined above.
\end{theorem}
\noindent By the contrapositive of Theorem \ref{classicalTheorem}, if there are no bialgebra maps between $H^*$ and $H^{cop}$, then there can be no quasitriangular structures on $H$.

\begin{proposition}
\label{factorizability}[\cite{shimizu_2016} Theorem 1.1, \cite{etingof_gelaki_nikshych_ostrik_2016}, Exercise 8.6.4]
    For a finite-dimensional quasitriangular Hopf algebra $(H, R)$ the following two conditions are equivalent
    \begin{enumerate}
        \item the monodromy matrix $M$ induces a non-degenerate pairing on $H^*\otimes H^* \rightarrow \C$ given by $h^*_1 \otimes h^*_2 \mapsto h^*_1 (M') h^*_2(M'')$, or
        \item the Drinfeld-Reshetikhin map 
        \begin{align*}
             H^* &\rightarrow H \\
             h^* &\mapsto h^*(M')M''
        \end{align*}
        is an isomorphism of vector spaces.

    \end{enumerate}
    If either item is fulfilled $(H, R)$ is called \textit{factorizable}, and otherwise \textit{non-factorizable}. 
\end{proposition}
\noindent For the purpose of this discussion we introduce a stronger notion.
\begin{definition}
\label{SNFdef}
    A unimodular, quasitriangular Hopf algebra $(H, R)$ with an integral $\lambda \in H^*$ is called \textit{strongly non-factorizable} if any of the two equivalent conditions is fulfilled
    \[
        \lambda(S(M')) M'' = 0,
    \]
    or
    \[
        S(M') \lambda(M'') = 0.
    \]
\end{definition}
\noindent It is easy to see from Definition \ref{factorizability} that any strongly non-factorizable Hopf algebra is non-factorizable.
\begin{remark}
    Note that if $(H, R)$ is non-cosemisimple, triangular, and has a two-sided integral $\lambda$, then it is strongly non-factorizable. Indeed $\lambda(\1)=0$ is equivalent to the non-cosemisimplicity condition for a unimodular Hopf algebra (see \cite{radford_2012}, Theorem 10.3.2).
\end{remark}

We now turn to the ribbon structure.
\begin{definition}
\label{ribbonAxioms} Let $(H, R)$ be a finite dimensional, quasitriangular Hopf algebra with monodromy matrix $M$. A \textit{ribbon element} $v\in Z(H)$ is a central element of $H$ such that
\begin{center}
\begin{varwidth}{\textwidth}
\begin{itemize}
    \item[(R1)] $S(v) = v$
    \item[(R2)] $\epsilon(v)=1$
    \item[(R3)] $M\Delta(v) = v\otimes v$.
\end{itemize}
\end{varwidth}
\end{center}
The triple $(H, R, v)$ is called a \textit{ribbon Hopf algebra}.
\end{definition}

\begin{definition}
\label{drinfeldelem}
    The \textit{Drinfeld element} of a quasitriangular Hopf algebra $(H, R)$ is
    \[
        u:= S(R'') R'.
    \]
    It is invertible with inverse given by 
    \[
        u^{-1} = S^{-1}\left((R^{-1})''\right) (R^{-1})'.
    \]
\end{definition}
\begin{proposition}[\cite{kassel_2012}, VIII.4.1]\label{DrinfeldProp1} Let $(H, R)$ be a quasitiangular Hopf algebra and $u$ be the corresponding Drinfeld element. For any $h \in H$
    \[
        S^2(h) = u h u^{-1}.
    \]
\end{proposition}
\begin{definition}[\cite{radford_1994}]
    Let $H$ be a Hopf algebra. A \textit{pivotal element} is a grouplike element $g\in G(H)$ such that for any $h\in H$.
    \[
        S^2(h) = g h g^{-1}.
    \]
\end{definition}
\noindent As explained in \cite[Prop. 2]{radford_1994} in a ribbon Hopf algebra, a ribbon element $v=g^{-1}u$ is the product of the Drinfeld element $u$ and the inverse of a pivotal element. However, not every pivotal element can be used to define a ribbon element. 
\begin{definition}
    Let $H$ be  Hopf algebra admitting a left integral $\lambda_L \in H^*$. The the \textit{distinguished grouplike element} $a\in G(H)$ is the unique element such that $\lambda_L \circ S(h) = \lambda_L(ah)$. 
\end{definition}
\begin{lemma}[\cite{radford_1994} Proposition 2.(b)] \label{2SidedCoIntLemm} Let $(H, R)$ be a finite-dimensional quasitriangular Hopf algebra with a two-sided cointegral and denote the Drinfeld element by $u$. Let $g \in G(H)$ be a pivotal element and $a$ be the distinguished grouplike element. Then $v:=g^{-1}u$ is a ribbon element if and only if 
\[
    g^2 = a^{-1}.
\]
\end{lemma}

\begin{lemma}
\label{2SidedIntLemm}
    Let $(H, R)$ be a finite-dimensional quasitriangular Hopf algebra with a two-sided integral $\lambda$. If $g\in G(H)$ is a pivotal element, then $g^2=\1$.
\end{lemma}
\begin{proof}
    If $\lambda$ is two-sided, then the distinguished grouplike element $a$ of Lemma \ref{2SidedCoIntLemm} satisfies $a=\1$. Thus, $g^2 = \1$.
\end{proof}
\noindent We introduce one more property of the ribbon structure that bears significance in many topological constructions.
\begin{definition}
\label{DefAF}
    Let $H$ be a Hopf algebra carrying a (left) integral $\lambda$ and a ribbon element $v$. We call $H$ \textit{anomaly-free} if 
    \[
        \lambda(v) \neq 0.
    \]
    Otherwise we will refer to it as \textit{anomalous}.
\end{definition}
\subsection{Algebra generators and tuple notation}
All the constructions considered in this paper produce Hopf algebras generated by grouplike  and skew-primitive generators. 
\begin{notation}
    We denote by $K_a$, $a=1, \dots, s$ the grouplike generators that commute with one another and are of respective orders $m_1, \dots, m_s$. By $X_k, X^\pm_k, Z_l^\pm$ for $k=1, \dots, t_1$ and $l=1, \dots, t_2$ we denote skew-primitive generators such that $(X_k)^2 = (X^\pm_k)^2 = (Z^\pm_l)^2=0$.
\end{notation}
\noindent We will repeatedly appeal to the following notion.
\begin{definition}
    Let $g, h\in H$. We say $g, h$ have a diagonal relation if there exists $\gamma \in \C$ such that 
    \[
        gh = \gamma hg.
    \]
\end{definition}
\begin{notation}
     We will use boldface for tuples of numbers and algebra generators, for instance, let $\mb{0}=(0, \dots, 0)$, and let
     \[
        \pmb{\xi} = (\xi_1, \dots, \xi_s):= (\exp^{2\pi/m_1}, \dots, \exp^{2\pi/m_s}) \in \C^s
     \]
     be the tuple of primitive roots of unity of orders $\mb{m} = (m_1, \dots, m_s) \in \Z^s_{>0}$. We similarly express the lists of algebra generators as $\mb{K}=(K_1, \dots, K_s)$, or $\mb{X} = (X_1, \dots, X_t)$.
\end{notation}
\noindent We will often use the following algebraic structures.
\begin{notation}
    Let $\Z_\mb{m} := \Z_{m_1}\times\dots\times \Z_{m_s}$, where $\Z_{m_a}:=\Z/m_a\Z$, and let $\C[\mb{K}] := \C[K_1, \dots, K_{m_s}] \subset H$ be the group subalgebra generated by $K_1, \dots, K_s$. It is isomorphic to $\C[\Z_\mb{m}]$ as a Hopf algebra.
\end{notation}
We also introduce some operations on tuples.
\begin{notation}
    For two tuples of the same shape, for instance $\mb{w}, \mb{v} \in \Z_\mb{m}$ we define the element-wise sum
    \[
        \mb{w} + \mb{v}:= (w_1+v_1, \dots, w_s+v_s),
    \]
    and the element-wise product
    \[
        \mb{w} \cdot \mb{v}:=(w_1 v_1, \dots, w_s v_s).
    \]
    By a power of a tuple by a tuple, for instance for $\mb{v} \in \Z_{\mb{m}}$, $\pmb{\xi}\in \C^s$ and $\mb{K}$, we mean
    \[
        \pmb{\xi}^\mb{v} = (\xi_1)^{v_1}(\xi_1)^{v_2}\dots(\xi_s)^{v_s},
    \]
    and
    \[
        \mb{K}^\mb{v} = (K_1)^{v_1}(K_2)^{v_2}\dots(K_s)^{v_s}.
    \]
\end{notation}
\begin{notation}
    We will also encounter \textit{matrices}, denoted with boldface sans-serif font, for instance 
    \begin{align*}
        \pmb{\alpha} = (\alpha_{kl})_{1 \leq k, l \leq t} \in \operatorname{Mat}_{t\times t}(\C) && \matfont{d}  = (d_{ka})_{\substack{1 \leq a \leq s \\ 1 \leq k \leq t}},\; \matfont{u} = (u_{ka})_{\substack{1 \leq a \leq s \\ 1 \leq k \leq t}} \in \operatorname{Mat}_{t\times s}(\Z_\mb{m}).
    \end{align*}
     The rows $\matfont{d}_k = (d_{ka})_{0\leq a \leq s}$, $\matfont{u}_k=(u_{ka})_{0\leq a \leq s}$ are considered as tuples (and conform to the boldface notation). Products with no operation symbols mean the usual matrix multiplication, provided the participating matrices have appropriate dimensions and entries in the same ring. For $\mb{w}, \mb{v} \in \Z_\mb{m}$,  
\[
  \mb{w} \mb{v}^T = w_1 v_1+\dots+w_s v_s,
\]
is the inner product of row tuples, valued in $\Z_{\mathrm{lcm}(m_1, \dots, m_s)}$, where $\mathrm{lcm}(m_1, \dots, m_s)$ is the least common multiple of $m_1, \dots, m_s$. For a matrix $\matfont{z}\in \operatorname{Mat}_{s\times s}(\Z_\mb{m})$, the expression $\mb{v} \matfont{z}$ is a row tuple again, as is $(\mb{w}\mb{z})\cdot \mb{v}$, but $\mb{w}\matfont{z}\mb{v}^T$ is a scalar in $\Z_{\mathrm{lcm}(m_1, \dots, m_s)}$. 
\end{notation}

\section{The Nenciu construction}
The Nenciu construction of \cite{nenciu_2004} produces finite-dimensional ribbon Hopf algebras generated by grouplike generators collected in a tuple $\mb{K} = (K_1, \dots, K_s)$ and nilpotent generators collected in a tuple $\mb{X}=(X_1, \dots, X_t)$, such that the group subalgebra $\C[\mb{K}]$ is abelian, the nilpotent generators are skew-primitive with respect to $\C[\mb{K}]$ and all commutation relations between generators are diagonal. The choice of the generators and relations is strongly constrained, and sufficient conditions for the result to be a ribbon Hopf algebra were determined in \cite{nenciu_2004}. We will see that strongly non-factorizable examples of this construction exist.
\subsection{Hopf algebra structure}
We start by recalling Nenciu's construction: the fact that Definition \ref{genDef} indeed defines a Hopf algebra is a Theorem of \cite{nenciu_2004}.
\begin{definition}
\label{genDef}
    Let $\mb{m} \in \Z_{>0}^s$ be a tuple of positive integers of length $s$, and $t\in \Z_{>0}$ and $\matfont{d}$, $\matfont{u}$ be dimension $t\times s$, $\Z$-valued matrices, such that
    \begin{equation}
\label{invConstr}
    \pmb{\xi}^{\matfont{d}_k \cdot \matfont{u}_l}\pmb{\xi}^{\matfont{d}_l \cdot \matfont{u}_k}=1, \text{ and } \pmb{\xi}^{\matfont{d}_k \cdot \matfont{u}_k}=-1,
\end{equation}
    where $\pmb{\xi}:=(\exp^{2\pi/m_1}, \dots, \exp^{2\pi/m_s})$ is the tuple of primitive roots of unity of orders given by the tuple $\mb{m}\in \Z^s$.
    Define $H(\mb{m}, t, \matfont{d}, \matfont{u})$ to be the Hopf algebra generated by grouplike generators $\mb{K}=(K_1, \dots, K_s)$, and skew-primitive generators $\mb{X}=(X_1, \dots, X_t)$, and with the relations 
    \begin{align*}
        K^{m_a}_a=\1 && K_a K_b = K_b K_a && K_a X_k := \xi_a^{d_{ka}} X_k K_a\\
        X^2_k = 0 &&  X_l X_k := \pmb{\xi}^{\matfont{d}_k \cdot \matfont{u}_l} X_k X_l,
    \end{align*}
    coalgebra structure determined by
    \begin{align*}
        \epsilon(K_a) := 1 && \Delta(K_a) := K_a \otimes K_a \\
        \epsilon(X_k) := 0 && \Delta(X_k) := \1 \otimes X_k + X_k \otimes \mb{K}^{\matfont{u}_k},
    \end{align*}
    and antipode determined by
    \begin{align*}
        S(K_a) = K^{-1}_a := K^{m_a-1}_a && S(X_k) := -X_k \mb{K}^{-\matfont{u}_k},
    \end{align*}
    for $a, b = 1, \dots, s$, $s:=|\mb{m}|$ and $k, l = 1, \dots, t$. Here $\mb{m}$ is the tuple of orders of the grouplike elements, $\matfont{d}$ encode the diagonal relations for the grouplike elements, and $\matfont{u}$ the relations for the skew-primitive elements.
\end{definition}
\begin{remark}
    Note that in Definition \ref{genDef} the relations $X^2_k = 0$, for $k=1, ..., t$, result from the requirement that
    \[
        \pmb{\xi}^{\matfont{d}_k \cdot \matfont{u}_k}=-1
    \]
    holds, rather than being imposed independently. We include them in the definition for clarity.
\end{remark}
\noindent Sometimes, we will choose several variables $X_k$ to have the same commutation properties. We refer to this as \textit{type}, and in examples we will indicate each type with a different letter $X, Y, Z$ etc. if necessary. 
\begin{definition}
    By the \textit{type} of $X_k$ we mean the prescription of $\matfont{d}_k$ and $\matfont{u}_k$, that is $X_l$ is of the same type as $X_k$ if $\matfont{d}_k=\matfont{d}_l$ and $\matfont{u}_k=\matfont{u}_l$. Since all relations involving $X_k$ and $X_l$ are diagonal, we will call the two generators of \textit{opposite type} if for $a=1, \dots, s$
\[
    d_{ka} \equiv -d_{la} \;\; \mod \; m_a
\]
and 
\[
    u_{ka} \equiv -u_{la} \;\; \mod \; m_a.
\]
This essentially means that generators of opposite type commute with all other generators over reciprocal constants. We will sometimes indicate variables of opposite type using superscripts, for instance $X^+_k$, $X^-_k$.
\end{definition}
\begin{proposition}
     Let $X^+_k, X^-_k$ be nilpotent generators of opposite type, then $\{X^+_k, X^-_k\}=0$.
\end{proposition}
\begin{proof}
We have directly
    \[
    X^+_k X^-_k= \pmb{\xi}^{\matfont{d}_k \cdot \matfont{u}_k}X^-_k X^+_k = (-1)^{-1} X^-_k X^+_k = -X^-_k X^+_k.
    \]
\end{proof}
\begin{remark}
    When making general statements about $H(\mb{m}, t, \matfont{d}, \matfont{u})$ that do not invoke the types we will use the notation $X_k$ for any nilpotent generator, with no superscripts. In such sense, two nilpotent generators with distinct indices, say $X_k$, $X_l$, can be of different, opposite or same type. 

\end{remark}
\noindent We can define a \textit{monomial basis} for $H(\mb{m}, t, \matfont{d}, \matfont{u})$.
\begin{proposition}
    \label{monomialBasis} The Hopf algebra $H(\mb{m}, t, \matfont{d}, \matfont{u})$  has a monomial basis
    \[
        \{\mb{K}^\mb{v} \mb{X}^\mb{r} |\mb{v}\in \Z_\mb{m},  \mb{r} \in \Z^t_2\}.
    \]
\end{proposition}
\begin{proof}
    The family  $\{\mb{K}^\mb{v} \mb{X}^\mb{r} |\mb{v}\in \Z_\mb{m},  \mb{r} \in \Z^t_2\}$ is free since different elements involve different sets of generators, and the relations are diagonal. Furthermore, any element in the Hopf algebra is, by definition, a linear combination of products of elements of $\mb{K}$ and $\mb{X}$. Since the relations are diagonal, up to a scalar, each term can be reordered so that the grouplike generators come first with a fixed order for all generators. Finally, the term is $0$ if a given skew-primitive generator appears twice. Indeed, the relations are diagonal and the nilpotent generators square to $0$.
    
\end{proof}

\begin{definition}
    We call $T := \mb{X}^\mb{r}$ for $\mb{r}=(1, \dots, 1)$ the \textit{top element}. 
\end{definition}
\noindent It will play an important role in the unimodular structure we will sometimes define on $H(\mb{m}, t, \matfont{d}, \matfont{u})$.
\noindent We will also need to use the adjoint action of generators on one another (see Definition \ref{adjointaction}). By direct computation, we get:
\begin{proposition}
\label{adjointAct}
     The adjoint action of the generators of $H$ is given by the following formulas.
     A grouplike $\mb{K}^\mb{w}$, $\mb{w} \in \Z_{\mb{m}}$ acts on the skew-primitive generator $X_k$ as
    \[
        \mb{K}^\mb{w} \rhd X_k = \mb{K}^\mb{w} X_k (\mb{K}^\mb{w})^{-1} = \pmb{\xi}^{\mb{w} \cdot \matfont{d}_k} X_k,
    \]
    a skew-primitive $X_k$ acts on a grouplike $\mb{K}^\mb{w}$ as
    \[
        X_k \rhd \mb{K}^\mb{w} = \mb{K}^\mb{w} S(X_k) + X_k \mb{K}^\mb{w} \mb{K}^{-\matfont{u}_k} = (1-\pmb{\xi}^{\mb{w} \cdot \matfont{d}_k})X_k \mb{K}^{-\matfont{u}_k}\mb{K}^{\mb{w}},
    \]
    and a skew-primitive $X_k$ acts on a skew-primitive $X_l$ as 
    \[
        X_k \rhd X_l = X_l S(X_k) + X_k X_l \mb{K}^{-\matfont{u}_k} = (\pmb{\xi}^{\matfont{u}_k \cdot \matfont{d}_l}-1)X_k \mb{K}^{-\matfont{u}_k}X_l.
    \]
\end{proposition}

\subsection{Examples}
     In this subsection we study explicit examples of the Nenciu construction. We start with a (not strongly) non-factorizable Hopf algebra, which can be retrieved from Nenciu's construction, and appears in \cite{nenciu_2004} as $E(n)$, first studied in \cite{panaite_oystaeyen_1999}. Here we will refer to it as \textit{symplectic fermions}, after \cite{gainutdinov_runkel_2017}, also appearing in \cite{faitg_gainutdinov_schweigert_2024}. The notation using the variables $L, Z^\pm_l$ is introduced to be used in further examples.
 \begin{definition}
 \label{SFdef}
     The \textit{Hopf algebra of symplectic fermions}, $\operatorname{SF}_{2n}$ is the Hopf algebra generated by $K_1 = L$ and, $X_{2l-1} = Z^+_l$, $X_{2l} = Z^-_l$, for $k, l=1, \dots, n$, subject to the relations
    \[
         L^2 = \1, \;\; L Z^\pm_l = - Z^\pm_l L, \;\; \{Z^\pm_l, Z^\pm_k\}=\{Z^\pm_l, Z^\mp_k\}=0.
    \]
    The Hopf structure is defined by
    \[
        \epsilon(L) = 1, \;\; \epsilon(Z^{\pm}_l)= 0,
    \]
    \[
        \Delta(L) = L \otimes L, \;\; \Delta(Z^{\pm}_l) = \1 \otimes Z^{\pm}_l + Z^{\pm}_l \otimes L,
    \]
    \[
        S(L)=L^{-1}=L, \;\; S(Z^{\pm}_l) = -Z^{\pm}_l L.
    \]
    The dimension of the algebra is $2^{2n+1}$. The corresponding Nenciu data is 
    \begin{itemize}
        \item $\mb{m}=(2)$, so $s=1$
        \item $t=2n$
        \item $\matfont{d} = \matfont{u} = (1, \dots, 1)^T$ are column vectors wherein each row is $d_k = 1$ and $u_k = 1$, $k=1, \dots, t$
         \begin{align*}
            \matfont{d}=\begin{pmatrix}
                1\\
                \multicolumn{1}{c}{$\vdots$} \\
                1\\
            \end{pmatrix},&&
            \matfont{u}=\begin{pmatrix}
                1\\
                \multicolumn{1}{c}{$\vdots$} \\
                1\\
            \end{pmatrix}.
        \end{align*}
    \end{itemize}
 \end{definition}
\noindent We now increase the number of grouplike and nilpotent generators, to produce previously unknown examples that will all turn out to admit only strongly non-factorizable quasitriangular structures.
\begin{example}
\label{ex1}
    Let $\operatorname{N}_1$ be the Hopf algebra generated by $K_a$, $X_k$, $Z^{\pm}_l$ for $a\in \{1, 2\}$, $j, k\in \{1, \dots, t_1\}$; $ t_1\in 4\N$ and $l, m\in \{1, \dots, t_2\}$; $t_2\in \N$, subject to the following relations
    \[
         K_a^4 = \1, \;\; K_a X_j = i X_j K_a, \;\; K_1 Z^{\pm}_l = \pm iZ^\pm_l K_1, \;\; K_2 Z^{\pm}_l =   Z^\pm_l K_2,
    \]
    \[
        Z^\pm_l X_k = \mp i X_k Z^\pm_l, \;\;\;  \{X_j, X_k\}= \{Z^{\pm}_l, Z^{\pm}_m\}=\{Z^{\pm}_l, Z^{\mp}_m\} =  0,
    \]
    where $i^2=-1$. The Hopf structure is defined by
    \[
        \epsilon(K_a) = 1, \;\; \epsilon(X_k) =\epsilon(Z^{\pm}_l)= 0,
    \]
    \[
        \Delta(K_a) = K_a \otimes K_a, \;\; \Delta(X_k) = \1 \otimes X_k + X_k \otimes K_1 K_2, \;\; \Delta(Z^{\pm}_l) = \1 \otimes Z^{\pm}_l + Z^{\pm}_l \otimes K_1^{\pm 2} K_2^{\pm 1} ,
    \]
    \[
        S(K_a)=K^{-1}_a, \;\; S(X_k) = -X_k (K_1 K_2)^{-1}, \;\; S(Z^{\pm}_l) = -Z^{\pm}_l K_1^{\mp 2} K_2^{\mp 1}.
    \]
    The Nenciu data for this Hopf algebra is 
    \begin{itemize}
        \item  $\mb{m} = (4, 4)$, with $s=2$
        \item  $t = t_1+2t_2$
        \item  $\matfont{d}$ and $\matfont{u}$ are $t\times 2$ matrices, such that for $k=1, \dots, t_1$; $\matfont{d}_k = (1, 1)$, and  $\matfont{u}_k = (1, 1)$ corresponding to $X_k$. For $l=t_1+1, \dots, t_2$ the rows are $\matfont{d}_{2l-1} = (1, 0)$, $\matfont{u}_{2l-1} = (2, 1)$ corresponding to $Z^+_l$ and $\matfont{d}_{2l} = (-1, 0)$, $\matfont{u}_{2l} = (-2, -1)$ corresponding to $Z^-_l$.
        \begin{align*}
            \matfont{d}=\begin{pmatrix}
                1&&1\\
                \multicolumn{3}{c}{$\vdots$} \\
                1&&1\\
                1&&0\\
                -1&&0\\
                \multicolumn{3}{c}{$\vdots$} \\
                 1&&0\\
                -1&&0\\
            \end{pmatrix},&&
            \matfont{u}=\begin{pmatrix}
                1&&1\\
                \multicolumn{3}{c}{$\vdots$} \\
                1&&1\\
                2&&1\\
                -2&&-1\\
                \multicolumn{3}{c}{$\vdots$} \\
                 2&&1\\
                -2&&-1\\
            \end{pmatrix}.
        \end{align*}
    \end{itemize}
\end{example}
\noindent The instance of Example \ref{ex1} with the smallest dimension occurs for $t_1 = 4$, $t_2 = 1$, and its dimension is $4\times 4 \times 2^6 = 2^{10} = 1024$. It will be shown later this Hopf algebra admits only triangular ribbon structures.
\\ 
The following examples were designed to admit strongly non-factorizable quasitriangular structures which are not triangular, as will be shown later.
\begin{example}
    \label{ex2}
    Let $\operatorname{N}_2$ be the Hopf algebra generated by $K_a$, $X_k^{\pm}$ and $Z^{\pm}_l$, for $a=1, 2, 3$, $j, k=1, \dots, t_1$; $ t_1\in \N$ and $l, m =1, \dots, t_2$, $t_2\in \N$, subject to the following relations
    \[
        K_a^4 = \1, \;\;  \;\; K_a X^{\pm}_k = \pm i X^{\pm}_k K_a
    \]
    \[
       K_1 Z^{\pm}_l = Z^{\pm}_l K_1, \;\; K_2 Z^{\pm}_l = - Z^{\pm}_l K_2, \;\;  K_3 Z^{\pm}_l = \pm i Z^{\pm}_l K_3
    \]
    \[
        \{X^{\pm}_j, X^{\pm}_k\} = \{X^{\pm}_j, X^{\mp}_k\} = \{Z^{\pm}_l, X^{\pm}_k\} = \{Z^{\pm}_l, X^{\mp}_k\} = \{Z^{\pm}_l, Z^{\pm}_m\}=\{Z^{\pm}_l, Z^{\mp}_m\} =  0.
    \]
    Let also $L:=K_3^2$ as a shorthand, note that $L^2 = \1$. The Hopf structure is defined by
    \[
        \epsilon(K_a) = 1, \;\; \epsilon(X^{\pm}_k)=\epsilon(Z^{\pm}_l) = 0,
    \]
    \[
        \Delta(K_a) = K_a \otimes K_a, \;\; \Delta(X^{\pm}_k) = \1 \otimes X^{\pm}_k + X^{\pm}_k \otimes (K_1 K_2)^{\pm 1}, \;\; \Delta(Z^{\pm}_l) = \1 \otimes Z^{\pm}_l + Z^{\pm}_l \otimes L
    \]
    \[
       S(K_a)=K^{-1}_a \;\; S(X^{\pm}_k) = -X^{\pm}_k (K_1 K_2)^{\mp1} \;\;  S(Z^{\pm}_l) = -Z^{\pm}_l L.
    \]
    The Nenciu data for this Hopf algebra is 
    \begin{itemize}
        \item  $\mb{m} = (4, 4, 4)$, so $s=3$
        \item  $t = 2t_1+2t_2$
        \item Then $\matfont{d}$ and $\matfont{u}$ are $t\times 3$ matrices, such that for $j=1, \dots, t_1$ $\matfont{d}_{2j-1} = (1, 1, 1)$, and  $\matfont{u}_{2j-1} = (1, 1, 0)$ corresponding to $X^+_j$ and $\matfont{d}_{2j} = (-1, -1, -1)$, and  $\matfont{u}_{2j} = (-1, -1, 0)$ corresponding to $X^-_j$. \\
        For $l=2t_1+1, \dots, t_2$ the rows are $\matfont{d}_{2l-1} = (0, 2, 1)$, $\matfont{u}_{2l-1} = (0,0,2)$ corresponding to $Z^+_l$ and $\matfont{d}_{2l} = (0, -2, -1)$, $\matfont{u}_{2l} = (0, 0, -2)$ corresponding to $Z^-_l$.
         \begin{align*}
            \matfont{d}=\begin{pmatrix}
                1&&1&&1\\
                -1&&-1&&-1\\
                \multicolumn{5}{c}{$\vdots$} \\
               1&&1&&1\\
                -1&&-1&&-1\\
                0&&2&&1\\
                0&&-2&&-1\\
                \multicolumn{5}{c}{$\vdots$} \\
                 0&&2&&1\\
                0&&-2&&-1\\
            \end{pmatrix},&&
            \matfont{u}=\begin{pmatrix}
               1&&1&&0\\
                -1&&-1&&0\\
                \multicolumn{5}{c}{$\vdots$} \\
               1&&1&&0\\
                -1&&-1&&0\\
                0&&0&&2\\
                0&&0&&-2\\
                \multicolumn{5}{c}{$\vdots$} \\
                 0&&0&&2\\
                0&&0&&-2\\
            \end{pmatrix}.
        \end{align*}
    \end{itemize}
    
\end{example}
\noindent The instance of Example \ref{ex2} with the smallest dimension occurs for $t_1 = 1$, $t_2 = 1$, and its dimension is $4\times 4\times 4 \times 2^4 = 2^{10} = 1024$. 

\noindent Note that in Example \ref{ex2} all the relations between skew-primitive generators are anti-commutations. By introducing new skew-primitive generators $Y^{\pm}_k$, we can produce commutations involving constants of higher order, while retaining the strongly non-factorizable, but not triangular, quasitriangular structure. 
\begin{example}
    \label{ex3}
    Let $\operatorname{N}_3$ be the Hopf algebra generated by $K_a$, $X_k^{\pm}$, $Y^{\pm}_k$ and $Z^{\pm}_l$, for $a=1, 2, 3$, $j, k=1, \dots, t_1$; $ t_1\in \N$ and $l, m =1, \dots, t_2$; $t_2\in \N$, subject to the following relations
    \[
        K_a^4 = \1, \;\;  \;\; K_a X^{\pm}_k = \pm i X^{\pm}_k K_a, 
    \]
    \[
       K_1 Y^{\pm}_k = \pm i Y^{\pm}_k K_1, \;\; K_2 Y^{\pm}_k = Y^{\pm}_k K_2, \;\;  K_3 Y^{\pm}_k = \pm i Y^{\pm}_l K_3
    \]
    \[
       K_1 Z^{\pm}_l = Z^{\pm}_l K_1, \;\; K_2 Z^{\pm}_l = -Z^{\pm}_l K_2, \;\;  K_3 Z^{\pm}_l = \pm i Z^{\pm}_l K_3
    \]
    \[
         \{X^{\pm}_j, X^{\pm}_k\} = \{X^{\pm}_j, X^{\mp}_k\} =  \{Y^{\pm}_j, Y^{\pm}_k\} = \{Y^{\pm}_j, Y^{\mp}_k\}, \;\; X_j Y_k = i Y_k X_j
    \]
    \[
       \{Z^{\pm}_l, X^{\pm}_k\} = \{Z^{\pm}_l, X^{\mp}_k\} = \{Z^{\pm}_l, Y^{\pm}_k\} = \{Z^{\pm}_l, Y^{\mp}_k\} = \{Z^{\pm}_l, Z^{\pm}_m\}=\{Z^{\pm}_l, Z^{\mp}_m\} =  0.
    \]
    Let also $L:=K_3^2$ as a shorthand, note that $L^2 = \1$. The Hopf structure is defined by
    \[
        \epsilon(K_a) = 1, \;\; \epsilon(X^{\pm}_k)= \epsilon(Y^{\pm}_k)=\epsilon(Z^{\pm}_l) = 0, 
    \]
    \[
        \Delta(K_a) = K_a \otimes K_a, \;\; \Delta(X^{\pm}_k) = \1 \otimes X^{\pm}_k + X^{\pm}_k \otimes (K_1 K_2)^{\pm 1},
    \]
    \[
         \Delta(Y^{\pm}_k) = \1 \otimes Y^{\pm}_k + Y^{\pm}_k \otimes (K_1 K^2_2)^{\pm1}, \;\; \Delta(Z^{\pm}_l) = \1 \otimes Z^{\pm}_l + Z^{\pm}_l \otimes L
    \]
    \[
       S(K_a)=K^{-1}_a \;\; S(X^{\pm}_k) = -X^{\pm}_k (K_1 K_2)^{\mp1}, \;\; S(Y^{\pm}_k) = -Y^{\pm}_k (K^2_1 K_2)^{\mp1} \;\;  S(Z^{\pm}_l) = -Z^{\pm}_l L.
    \]
    
     The Nenciu data for this Hopf algebra is 
    \begin{itemize}
        \item  $\mb{m} = (4, 4, 4)$, so $s=3$
        \item  $t = 4t_1+2t_2$
        \item Then $\matfont{d}$ and $\matfont{u}$ are $t\times 3$ matrices, such that for $j=1, \dots, t_1$ $\matfont{d}_{2j-1} = (1, 1, 1)$, and  $\matfont{u}_{2j-1} = (1, 1, 0)$ corresponding to $X^+_j$ and $\matfont{d}_{2j} = (-1, -1, -1)$, and  $\matfont{u}_{2j} = (-1, -1, 0)$ corresponding to $X^-_j$. \\
        For $k=t_1+1, \dots, 2t_1$, $\matfont{d}_{2k-1} = (1, 0, 1)$, and  $\matfont{u}_{2j-1} = (2, 1, 0)$ corresponding to $Y^+_k$ and $\matfont{d}_{2k} = (-1, -0, -1)$, and  $\matfont{u}_{2k} = (-2, -1, 0)$ corresponding to $Y^-_k$.\\
        For $l=2t_1+1, \dots, t_2$ the rows are $\matfont{d}_{2l-1} = (0, 2, 1)$, $\matfont{u}_{2l-1} = (0,0,2)$ corresponding to $Z^+_l$ and $\matfont{d}_{2l} = (0, -2, -1)$, $\matfont{u}_{2l} = (0, 0, -2)$ corresponding to $Z^-_l$.
       \begin{align*}
            \matfont{d}=\begin{pmatrix}
                1&&1&&1\\
                -1&&-1&&-1\\
                \multicolumn{5}{c}{$\vdots$} \\
               1&&1&&1\\
                -1&&-1&&-1\\
                1&&0&&1\\
                -1&&0&&-1\\
                \multicolumn{5}{c}{$\vdots$} \\
               1&&0&&1\\
                -1&&0&&-1\\
                0&&2&&1\\
                0&&-2&&-1\\
                \multicolumn{5}{c}{$\vdots$} \\
                 0&&2&&1\\
                0&&-2&&-1\\
            \end{pmatrix},&&
            \matfont{u}=\begin{pmatrix}
               1&&1&&0\\
                -1&&-1&&0\\
                \multicolumn{5}{c}{$\vdots$} \\
               1&&1&&0\\
                -1&&-1&&0\\
                2&&1&&0\\
                -2&&-1&&0\\
                \multicolumn{5}{c}{$\vdots$} \\
               2&&1&&0\\
                -2&&-1&&0\\
                0&&0&&2\\
                0&&0&&-2\\
                \multicolumn{5}{c}{$\vdots$} \\
                 0&&0&&2\\
                0&&0&&-2\\
            \end{pmatrix}.
        \end{align*}
    \end{itemize}
\end{example}
\noindent The instance of Example \ref{ex3} with the smallest dimension occurs for $t_1 = 1$, $t_2 = 1$, and its dimension is $4\times 4\times 4 \times 2^6 = 2^{12} = 4096$.

\subsection{Unimodular structure}
In this short section we construct integrals and cointegrals on $H(\mb{m}, t, \matfont{d}, \matfont{u})$ and determine the sufficient condition for unimodularity. 
\begin{definition}[\cite{nenciu_2004}, Proposition 4.7]
Let $ H =  H(\mb{m}, t, \matfont{d}, \matfont{u})$  be a Hopf algebra such as in Definition \ref{genDef}. Then $H$ admits a left cointegral
        \[
            \Lambda_L := \prod^s_{a=1} \left(\sum^{m_a-1}_{b=0} K^b_a \right) \prod^t_{k=1} X_k,
        \]
        and a right cointegral
        \[
            \Lambda_R := \prod^t_{k=1} X_k \prod^s_{a=1} \left(\sum^{m_a-1}_{b=0} K^b_a \right).
        \]
        Let us denote the top element by $T:=\prod^t_{k=1} X_k = \mb{X}^{(1, \dots, 1)}$ then a two-sided integral $\lambda$ is defined on the monomial basis by the formula
        \[
            \lambda(\mb{K}^\mb{v} \mb{X}^\mb{r} ):=
            \begin{cases}
                1 \text{\, if \,}\mb{K}^\mb{v} \mb{X}^\mb{r}  = T \\
                0 \text{\, otherwise}.
            \end{cases}
        \]
        
\end{definition}
\begin{remark}
    Note that $ \prod^s_{a=1}\left(\sum^{m_a-1}_{b=0} K^b_a \right)$ is simply the sum of all monomials in the Hopf subalgebra $\C[\mb{K}]$.
\end{remark}
\noindent 
\begin{proposition}
\label{unimodTheorem}
    We have $\Lambda_L=\Lambda_R$ if $T$ is central, which happens if and only if for all $a=1, \dots, s$
    \[
        \sum^t_{k=1} d_{ka} \equiv 0 \;\; \mod \;\; m_a.
    \] In particular $\Lambda_L$ is then a two-sided cointegral.
\end{proposition}
\begin{proof}
Firstly, using the monomial basis, and the fact that relations are diagonal, we see that $\Lambda_L = \Lambda_R$ if and only if $ \mb{K}^\mb{v}$ commutes with $T$ for every $\mb{v}$. This is tantamount to $T$ being central, since $T$ always commutes with any $X_j$. Indeed $TX_j = X_jT = 0$, as all skew-primitive generators square to $0$ and relations are diagonal. Secondly, this happens if and only if $T$ commutes with $K_a$ for every $a$. But we can compute directly that $$K_aT = \prod_{k=1}^t \xi_a^{d_{ka}} TK_a.$$ Thus, $T$ is central if and only if  $\sum^t_{k=1} d_{ka} \equiv 0$ in $\Z_{m_a}$.
\end{proof}

\subsection{Quasitriangular structure}
In \cite{nenciu_2004} a sufficient list of conditions is given for $H(\mb{m}, t, \matfont{d}, \matfont{u})$ to admit a ribbon structure. In this section, we recall these results, and we fix an imprecision concerning the existence of the pivotal element in the unimodular case, which we discuss in further sections. These rely on Theorem \ref{classicalTheorem}. We then revisit some examples of Nenciu and prove they admit unimodular and strongly non-factorizable structures. We end the section by computing the monodromy matrix associated to these examples.

As we have recalled for any Hopf algebra $H$ its dual $H^*$ is again a Hopf algebra. In our setting we have the following result:
\begin{proposition}[\cite{nenciu_2004}, Proposition 3.1] 
\label{dualProp}
Let $H(\mb{m}, t, \matfont{d}, \matfont{u})$ be as in Definition \ref{genDef}, then $H(\mb{m}, t, \matfont{d}, \matfont{u})^*$ is isomorphic to  $H(\mb{m}, t, \matfont{u}, \matfont{d})$ as a Hopf algebra, and the duals of the generators $K_a$ and $X_k$, denoted $\mathcal{K}_a$ and $\mathcal{X}_k$, correspond to the generators of $H(\mb{m}, t, \matfont{u}, \matfont{d})$ through this isomorphism.
\end{proposition}
\begin{notation}
    In the discussion that follows we will distinguish the algebra $H(\mb{m}, t, \matfont{d}, \matfont{u})$ and its dual, so we introduce a notation where $H(\mb{m}, t, \matfont{d}, \matfont{u})^* \cong H(\mb{m}, t, \matfont{u}, \matfont{d})$ is generated by grouplike generators $\mathcal{K}_a$ dual to $K_a$, for $a=1, \dots, s$ and skew-primitive generators $\mathcal{X}_k$ dual to $X_k$ for $k=1, \dots, t$ in the caligraphic font, satisfying the relations
    \[
         \mathcal{K}^{m_a}_a=\1, \;\;\; \mathcal{K}_a \mathcal{K}_b =  \mathcal{K}_b \mathcal{K}_a, \;\;\; \mathcal{K}_a \mathcal{X}_k := \pmb{\xi}^{u_{ka}} \mathcal{X}_k \mathcal{K}_a, \;\;\;  \mathcal{X}_l \mathcal{X}_k=\pmb{\xi}^{\mb{u_k}\cdot \mb{d_l}} \mathcal{X}_k \mathcal{X}_l,
    \]
    \[
         \epsilon(\mathcal{K}_a) := 1, \;\;\; \Delta(\mathcal{K}_a) := \mathcal{K}_a \otimes \mathcal{K}_a, \;\;  \epsilon(\mathcal{X}_k) := 0, \;\;\; \Delta(\mathcal{X}_k) := \1 \otimes \mathcal{X}_k + \mathcal{X}_k \otimes \pmb{\mathcal{K}}^{\matfont{d}_k}
    \]
    and antipode
    \begin{align}
        S(\mathcal{K}_a) := \mathcal{K}^{-1}_a =\mathcal{K}^{m_a-1}_a && S(\mathcal{X}_k) := -\mathcal{X}_k \pmb{\mathcal{K}}^{-\matfont{d}_k},
    \end{align}
    for $a, b=1, \dots, s$, and $k, l=1, \dots, t$. We also have the the monomial basis of the dual $H(\mb{m}, t, \matfont{d}, \matfont{u})$ denoted by 
    \[
        \{\pmb{\mathcal{K}}^\mb{v} \pmb{\mathcal{X}}^\mb{r} |\mb{v}\in \Z_\mb{m}, \mb{r} \in \Z^t_2\},
    \]
    with the duality pairing for any $\mb{w},\mb{v} \in \Z_\mb{m}$ and $\mb{p},\mb{r} \in \Z^t_2$ we have 
    \[
        \pmb{\mathcal{K}}^\mb{w}\pmb{\mathcal{X}}^\mb{p}(\mb{K}^\mb{v}\mb{X}^\mb{r}) = \prod^s_{a=1} \delta_{w_a, v_a}\prod^t_{k=1} \delta_{p_k, r_k}.
    \]
\end{notation}

\noindent The following is a Corollary of Proposition \ref{dualProp}.

\begin{corollary}[\cite{nenciu_2004}, Section 3]
\label{paramCor}
Any Hopf algebra map $f: H(\mb{m}, t, \matfont{d}, \matfont{u})^* \rightarrow H(\mb{m}, t, \matfont{d}, \matfont{u})^{cop}$ is parametrized by a $\Z$-valued $s\times s$ matrix $\matfont{z}$ and a $\C$-valued $t \times t$ matrix $\pmb{\alpha}$ as follows:
\[
    \mathcal{K}_a \mapsto \mb{K}^{\matfont{z}_a},
\]
where $\matfont{z}_a\in \Z^{s}$ is the $a$-th row of $\matfont{z}$ for $a=1, \dots, s$, and
\[
    \mathcal{X}_k \mapsto \sum^t_{l=1} \alpha_{kl} \mb{K}^{-\matfont{u}_l} X_l
\]
for $k, l=1, \dots, t$. With these definitions the image of an arbitrary element $\pmb{\mathcal{K}}^{\mb{v}} \pmb{\mathcal{X}}^\mb{p}$ for $\mb{v} \in \Z_{\mb{m}}$ and $\mb{p}\in \Z_2^t$, is 
\[
    f(\pmb{\mathcal{K}}^{\mb{v}} \pmb{\mathcal{X}}^\mb{p}) = \sum_{\mb{q} \in \Z_2^t} \left( \prod^t_{j, k=1} \alpha_{jk} p_j q_k \right)\mb{K}^{\mb{v}\matfont{z}} (\mb{K}^{-\matfont{u}_1} X_1, \dots, \mb{K}^{-\matfont{u}_t} X_t)^{\mb{q}}.  
\]
\end{corollary}
Thus, we have a general parametrization for any Hopf map $f: H(\mb{m}, t, \matfont{d}, \matfont{u})^* \rightarrow H(\mb{m}, t, \matfont{d}, \matfont{u})^{cop}$, where $H(\mb{m}, t, \matfont{d}, \matfont{u})^{cop}$ is the Hopf algebra $H(\mb{m}, t, \matfont{d}, \matfont{u})$ with the coproduct modified to $\Delta^{cop}$. We denote the elements of $H(\mb{m}, t, \matfont{d}, \matfont{u})^{cop}$ with the same symbols as those of $H(\mb{m}, t, \matfont{d}, \matfont{u})$, as the algebra structure remains unchanged.

The $\mb{p}$ and $\mb{q}$ are row tuples with values in $\{0, 1\}$ and they select a subset of skew-primitive generators present. In particular, 
    $(\mb{K}^{-\matfont{u}_1} X_1, \dots, \mb{K}^{-\matfont{u}_t} X_t)^{\mb{q}}$ is a product of terms of the form $\mb{K}^{-\matfont{u}_k} X_k$ according to the selection by $\mb{q}$
    \[
        (\mb{K}^{-\matfont{u}_1} X_1, \dots, \mb{K}^{-\matfont{u}_t} X_t)^{\mb{q}} = \prod^t_{k=1} q_k \mb{K}^{-\matfont{u}_k} X_k .
    \]

\begin{remark}
    Not any choice of matrices $\matfont{z}$ and $\pmb{\alpha}$ gives a valid bialgebra map, and we will list the necessary criteria below, in Theorem \ref{R-matrixConstr}. 
\end{remark}

\noindent This construction provides a weak converse to Theorem \ref{QTthm}, by reconstructing the R-matrix from the bialgebra map in the chosen parametrization. In essence, every term in the R-matrix expansion is a tensor product of a monomial in $H(\mb{m}, t, \matfont{d}, \matfont{u})$ and the image of its dual under $f$.

\begin{theorem}[\cite{nenciu_2004}, Theorem 4.1]
    \label{R-matrixConstr} 
     Let $R \in H(\mb{m}, t, \matfont{d}, \matfont{u}) \otimes H(\mb{m}, t, \matfont{d}, \matfont{u})$ be  the element defined by
    \[
        R:= R_{\matfont{z}} R_{\pmb{\alpha}},
    \]
    for
    \[
        R_{\matfont{z}} = \frac{1}{\prod^s_{a=1} m_a}\sum_{\mb{v}, \mb{w} \in \Z_\mb{m}} \pmb{\xi}^{-\mb{v}\cdot \mb{w}} \mb{K}^{\mb{w}} \otimes \mb{K}^{\mb{v}\matfont{z}}, 
    \]
    and 
    \[
        R_{\pmb{\alpha}} = \exp \left( \sum^{t}_{k, l=1} \alpha_{kl} X_k \otimes \mb{K}^{-\matfont{u}_l}X_l \right), 
    \]
    where $\matfont{z}\in \operatorname{Mat}_{s\times}(\Z_\mb{m})$ and $\pmb{\alpha} \in \operatorname{Mat}_{t\times t}(\C)$, are matrices as in Corollary \ref{paramCor}. Then $R$ is an R-matrix, if and only if the following constraints are satisfied.\\ 
    \begin{itemize}
        \item[(A1)] for any $a=1, \dots, s$, $b=1, \dots, s$
        \[
            m_a z_{ab} \equiv 0 \mod m_b
        \]
        \item[(A2)] for any $a=1, \dots, s$ and $k, l=1, \dots, t$
            \[
                \alpha_{kl} \pmb{\xi}^{\matfont{d}_l \cdot \matfont{z}_a}=\alpha_{kl} \xi_a^{u_{ka}}
            \]
        \item[(A3)] for any $i, j, k, l=1, \dots, t$
            \[
                (\alpha_{jl}\alpha_{ik}\pmb{\xi}^{\matfont{d}_k \cdot \matfont{u}_l} + \alpha_{jk}\alpha_{il}) \pmb{\xi}^{\matfont{d}_k\cdot \matfont{u}_l}= (\alpha_{jl}\alpha_{ik}\pmb{\xi}^{\matfont{d}_k \cdot \matfont{u}_l} + \alpha_{jk}\alpha_{il}) \pmb{\xi}^{\matfont{d}_j\cdot \matfont{u}_i}
            \]
        \end{itemize}
        
        \begin{itemize}
            \item[(B)] for any $k = 1, \dots, t$
                \[ 
                    \sum^t_{l=1}\alpha_{kl}\left(\1 \otimes \mb{K}^{-\matfont{u}_l} X_l + \mb{K}^{-\matfont{u}_l} X_l \otimes \mb{K}^{ \matfont{d}_k \matfont{z}} \right)=  \sum^t_{l=1} \alpha_{kl} \left(  \1 \otimes \mb{K}^{-\matfont{u}_l} X_l + \mb{K}^{-\matfont{u}_l} X_l \otimes \mb{K}^{-\matfont{u}_l} \right),  
                \]
        \end{itemize}

        \begin{itemize}
            \item[(C1)] for all $k = 1, \dots, t$
                \[
                    \matfont{d}_k \matfont{z}  = - \matfont{u}_k
                \]
            \item[(C2)] for all $k=1, \dots, t$, and $\mb{v} \in \Z_{\mb{m}}$
                \[
                    \pmb{\xi}^{\mb{v} \cdot \matfont{u}_k} = \pmb{\xi}^{\matfont{d}_k \cdot  \mb{v} \matfont{z}}
                \]
            \item[(C3)] whenever $(\mb{p}\pmb{\alpha})^{\mb{q}} \neq 0$, for some $\mb{p}, \mb{q} \in \Z^t_2$ and for all $a=1, \dots, s$
                \[
                    \sum^t_{l=1}(p_l + q_l) d_{la} \equiv 0 \mod m_a
                \]
            \item[(C4)]  whenever $(\mb{p}\pmb{\alpha})^{\mb{q}} \neq 0$, for some $\mb{p}, \mb{q} \in \Z^t_2$, then for $k=1, \dots, t$ 
                \[
                    \pmb{\xi}^{\matfont{d}_k \cdot (\mb{p} \matfont{u})} \pmb{\xi}^{\matfont{d}_k \cdot (\mb{q} \matfont{u})} =1.
                \]
    \end{itemize}
\end{theorem}
\noindent In (B) we have in particular that for any $k, l$, either $\alpha_{kl}=0$ or $\matfont{d}_k \matfont{z} = -\matfont{u}_l$. Note also that the latter two items (C3), (C4) are required if $R_{\pmb{\alpha}} \neq \1\otimes \1$, otherwise they are void. They are involved in realizing (QT5) for the grouplike $K_a$ and nilpotent $X_k$ generators for $a=1, ..., s$ and $k=1, ..., t$ respectively.
\begin{proof}
Let $f: H(\mb{m}, t, \matfont{d}, \matfont{u})^* \rightarrow H(\mb{m}, t, \matfont{d}, \matfont{u})^{cop}$ be the map induced by $R$ as in Theorem \ref{classicalTheorem}. In particular, (A1)-(A3) together with (B) are equivalent to (QT1)-(QT4) and (C1)-(C4) are equivalent to (QT5). 
We first of all note that (A1), (C1)-(C4) are a reformulation of \cite{nenciu_2004} Theorem 4.1 in our notation. We discuss the remaining conditions that, while implied in \cite{nenciu_2004}, were not collected in Theorem 4.1 there.
\begin{itemize}
\item[] Conditions (A1)-(A3) imply that $f$ is an algebra map. Towards (A2) let $\mathcal{K}_a$, $\mathcal{X}_k$ be the images of $K_a$, $X_k$, respectively, under the duality using the notation we introduced. Item (A2) is the expression of the fact that the commutation of $\mathcal{K}_a$, $\mathcal{X}_k$ has to be respected by the map $f$. Explicitly, if $\mathcal{K}_a \mathcal{X}_k = \pmb{\xi}^{\matfont{d}_a\cdot \matfont{u}_k} \mathcal{X}_k \mathcal{K}_a$, then
\[
    f(\mathcal{K}_a \mathcal{X}_k) = \mb{K}^{\matfont{z}_a}  \sum^t_{l=1} \alpha_{kl} \mb{K}^{-\matfont{u}_l} X_l = \sum^t_{l=1} \alpha_{kl} \pmb{\xi}^{\matfont{z}_a \cdot \matfont{d}_l} \mb{K}^{-\matfont{u}_l} X_l \mb{K}^{\matfont{z}_a}. 
\]
On the other hand
\[
    f(\mathcal{X}_k \mathcal{K}_a) = \xi_a^{u_{ka}}f(\mathcal{K}_a\mathcal{X}_k) =  \sum^t_{l=1} \alpha_{kl} \xi_a^{u_{ka}}\mb{K}^{-\matfont{u}_l} X_l.
\]
Comparing coefficients elementwise yields the required condition. This is enough as we defined the map on the monomial basis and all the relations are diagonal.\\
Similarly, (A3) provides this for the commutation of $\mathcal{X}_l$ with $\mathcal{X}_k$. The calculation is analogous to the case (A2).
\item[] Condition (B) expresses the fact that $f$ is a coalgebra map for $\mathcal{X}_k$, for $\mathcal{K}_a$ this already follows from (A1). We have that
\begin{align*}
    \Delta^{cop}(f(\mathcal{X}_k)) &= \Delta^{cop}\left(\sum^t_{l=1} \alpha_{kl} \mb{K}^{-\matfont{u}_l} X_l\right) \\&= \sum^t_{l=1} \alpha_{kl} \left( \mb{K}^{-\matfont{u}_l} X_l \otimes \mb{K}^{-\matfont{u}_l}  + \1 \otimes \mb{K}^{-\matfont{u}_l} X_l \right),
\end{align*}
but
\begin{align*}
    f(\Delta(\mathcal{X}_k)) &= f(\1)\otimes f(\mathcal{X}_k) + f(\mathcal{X}_k) \otimes f(\mb{K}^{\matfont{d}_k}) \\&= \1\otimes \sum^t_{l=1} \alpha_{kl} \mb{K}^{-\matfont{u}_l} X_l + \sum^t_{l=1} \alpha_{kl} \mb{K}^{-\matfont{u}_l} X_l \otimes \mb{K}^{-\matfont{u}_k}
\end{align*}
since $\matfont{d}_k\matfont{z} =-\matfont{u}_k$. In the latter expression combining the two sums yields
\[
    \sum^t_{l=1}\alpha_{kl}\left(\1 \otimes \mb{K}^{-\matfont{u}_l} X_l + \mb{K}^{-\matfont{u}_l} X_l \otimes \mb{K}^{ \matfont{d}_k \matfont{z}}  \right).
\]
Elementwise comparison of the two expressions retrieves (B), as above. 
\item[] The exponential in the piece $R_{\pmb{\alpha}}$ expands to the summation of all partial products between the elements of the sum in the exponent. In more detail, it can be expanded as
\[
     R_{\pmb{\alpha}} = \exp \left( \sum^{t}_{k, l=1} \alpha_{kl} X_k \otimes \mb{K}^{-\matfont{u}_l}X_l \right) = \sum_{j\in \operatorname{Mat}_{t\times t}(\Z_2)} \prod^{t}_{k, l=1} \alpha_{kl} (X_k)^{j_{kl}} \otimes (\mb{K}^{-\matfont{u}_l}X_l)^{j_{kl}}. 
\]
 Each piece of the rightmost sum corresponds in the former tensor factor to a monomial in $X_k$ generators, $\mb{X}^\mb{r}$, and in the latter tensor factor - its image under the composition of the duality (which, recall, is only a linear map), and the bialgebra map so the element $f((\mb{X}^\mb{r})^*) \in H^{cop}$.
The parameterization coefficients $\alpha_{kl}$ control whether a term appears in the sum, or equivalently in the image of $f$, or not, and with which coefficient.\\ 
To complete the argument, we can verify that $f_R: H^* \rightarrow H$ of Theorem \ref{classicalTheorem} retrieved from $R_{\matfont{z}}R_{\pmb{\alpha}}$ in fact coincides with $f$. This is because for any $\mb{u},\mb{o} \in \Z_\mb{m}$ and $\mb{p},\mb{r} \in \Z^t_2$ we have 
\[
    \pmb{\mathcal{K}}^\mb{u}\pmb{\mathcal{X}}^\mb{p}(\mb{K}^\mb{o}\mb{X}^\mb{r}) = \prod^s_{a=1} \delta_{u_a, o_a}\prod^t_{k=1} \delta_{p_k, r_k}.
\]
Therefore

\begin{align*}
    &f_R(\pmb{\mathcal{K}}^\mb{u}\pmb{\mathcal{X}}^\mb{p}) \\
    =&\sum_{\mb{w}, \mb{v}\in \Z_\mb{m}}\sum_{j\in \operatorname{Mat}_{t\times t}(\Z_2)} \prod^{t}_{k, l=1} \alpha_{kl} \pmb{\xi}^{-\mb{w}\cdot \mb{v}} \left(\pmb{\mathcal{K}}^\mb{u}(\pmb{\mathcal{X}}_k)^{p_k} \right)(\mb{K}^\mb{w}(X_k)^{j_{kl}}) \otimes \mb{K}^{\mb{v}\matfont{z}} (\mb{K}^{-\matfont{u}_l}X_l)^{j_{kl}}\\
    =&\sum_{\mb{w}, \mb{v}\in \Z_\mb{m}}\sum_{j\in \operatorname{Mat}_{t\times t}(\Z_2)} \prod^{t}_{k, l=1} \alpha_{kl} \pmb{\xi}^{-\mb{w}\cdot \mb{v}}  \prod^s_{a=1} \delta_{u_a, w_{a}} \delta_{p_k, j_{kl}} \mb{K}^{\mb{v}\matfont{z}} (\mb{K}^{-\matfont{u}_l}X_l)^{j_{kl}}
\end{align*}

which is exactly $f(\pmb{\mathcal{K}}^\mb{u}\pmb{\mathcal{X}}^\mb{p})$, as expected. 
\end{itemize}
\end{proof}


\subsection{Ribbon structure}
We established in Lemma \ref{2SidedIntLemm} that a unimodular quasitriangular Hopf algebra with two-sided integral and the Drinfeld element $u$ can be endowed with a ribbon structure if we can find a pivotal element $g\in G(H)$ such that $g^2=\1$.   
\begin{lemma}
    \label{DrinfledLemma2}
    Let $K_a$ and $X_k$ be generators of $H(\mb{m}, t, \matfont{d}, \matfont{u})$ and $R$ be an R-matrix. Denote $u$ the Drinfeld element as in Definition \ref{drinfeldelem}. Then
    \[
        [u, K_a] = \{u, X_k\} = 0.
    \]
\end{lemma}
\begin{proof}
By Proposition \ref{DrinfeldProp1}
    \[
        S^2(K_a) = K_a = u K_a u^{-1},  
    \]
and 
    \[
        S^2(X_k) = -X_k =  u X_k u^{-1},
    \]
the conclusion is immediate.    
\end{proof}
\noindent Hence we have the following result.
\begin{theorem}
\label{ribbonThm}
    Assume $H(\mb{m}, t, \matfont{d}, \matfont{u})$ admits a two-sided integral, a two-sided cointegral, an R-matrix $R$, and a pivotal element $g$. Then $g^2=1$ and $v = g^{-1}u = gu$ is a ribbon element.
\end{theorem}
\begin{proof}
    This follows from Lemmas \ref{2SidedCoIntLemm}, \ref{2SidedIntLemm} and \ref{DrinfledLemma2}.
\end{proof}

\subsection{Examples revisited}
 We revisit examples of the Nenciu construction, and show they all admit unimodular, strongly non-factorizable ribbon structures. 
 \begin{notation}
     In this section the expression $R_{\pmb{\alpha}}$ will recur in the same form given by
     \[
        R_{\pmb{\alpha}} := \exp\left(\sum^{t_2}_{l = 1} \alpha_l( Z^+_l \otimes L Z^-_l - Z^-_l \otimes L Z^+_l)\right).
     \]
     Here we stick to the notation for the nilpotent generators from the previous examples emphasizing the types. In particular, let $Z^+_l$ and $Z^-_l$ be nilpotent generators in $H(\mb{m}, t, \matfont{d}, \matfont{u})$ of opposite type, with neighbouring indices, that is $Z^+_l = X_k$ and $Z^-_l = X_{k+1}$ for some $k$. Let moreover $\{Z_l^\pm, X_k\} = 0$ for all 
     $X_k$, including $Z_l^\pm$ themselves.
     We also made a change in the notation in the expression above and we denoted $\alpha_{k, k+1} = -\alpha_{k+1, k}:=\alpha_l$. We collect $\alpha_l$ in the tuple $\pmb{\alpha}$. 
 \end{notation}
 \noindent We again start the list of examples with $\operatorname{SF}_{2n}$.
\begin{proposition}
\label{SFprop}
    The algebra $\operatorname{SF}_{2n}$ of Definition \ref{SFdef} is
    \begin{enumerate}
        \item unimodular, with a two-sided cointegral
        \[
            \Lambda := (\1+L)\prod^{n}_{l=1}Z^+_l
            \prod^{n}_{l=1}Z^-_l,
        \]
        and a two-sided integral defined on the monomial basis by
        \[
            \lambda(L^v \mb{X}^\mb{r} ):=
            \begin{cases}
                1 \text{\, if \,} v = 0, \mb{r} = (1, 1,\dots,1) \\
                0 \text{\, otherwise},
            \end{cases}
        \]
        \item quasitriangular, with the R-matrix defined for the grouplike generators as 
        \[
            R_{\matfont{z}} := \frac{1}{2}\sum_{v, w \in \Z_2} (-1)^{-vw} L^{w} \otimes L^{v z} = \frac{1}{2}(\1\otimes \1+ \1\otimes L + L\otimes \1 - L\otimes L), 
        \]
        where $z=1$.
        for some $\pmb{\alpha}=(\alpha_1, \dots, \alpha_{t_2}) \in \C^{2t_2}$, 
        so that 
        \[
            R = \frac{1}{2}(\1\otimes \1+ \1\otimes L + L\otimes \1 - L\otimes L) \exp\left(\sum^{n}_{l = 1} \alpha_l( Z^+_l \otimes L Z^-_l - Z^-_l \otimes L Z^+_l)\right),
        \]
        \item ribbon, with the ribbon element
        \[
            v:= \exp \left( - 2\sum^{n}_{k = l} \alpha_l Z^+_l Z^-_l  \right),
        \]
        corresponding to the pivotal element $g = L$, which in this example is the only such element available. 
    \end{enumerate}
\end{proposition}
\begin{proof}
    See \cite{panaite_oystaeyen_1999}, Propositions 8, 2 and 11, respectively.
\end{proof}
\begin{remark}
    Note that Theorem \ref{QTthm} allows multiple quasitriangular structures for $H(\mb{m}, t, \matfont{d}, \matfont{u})$. Indeed, it is known from \cite{gainutdinov_runkel_2017} that, for instance, $\operatorname{SF}_2$ admits a quasitriangular structure where $R:=R_z \tilde{R}_{\pmb{\alpha}}$, with $R_z$ as in Proposition \ref{SFprop}, but  
    \[
        \Tilde{R}_{\pmb{\alpha}} = \exp \left( \alpha^+ Z^+ \otimes L Z^+ + \alpha^- Z^- \otimes L Z^- \right)
    \]
    where $\pmb{\alpha} = (\alpha^+, \alpha^-) \in \C^2$. It admits only a trivial ribbon element $v=\1$. 
\end{remark}
\noindent We now move on to more complicated examples, that did not appear in \cite{nenciu_2004}, starting with one that is in fact \textit{triangular}.

\begin{proposition}
\label{ex1prop}
    The algebra $\operatorname{N}_1$ of Example \ref{ex1} is
    \begin{enumerate}
        \item unimodular, with a two-sided cointegral
        \[
            \Lambda := \left(\sum^3_{a,b=0} K^a_1 K^b_2 \right) \prod^{t_1}_{k=1}X_k \prod^{t_2}_{l=1}Z^+_l
            \prod^{t_2}_{l=1}Z^-_l,
        \]
        and a two-sided integral defined on the monomial basis by
        \[
            \lambda(\mb{K}^\mb{v}\mb{X}^\mb{r} ):=
            \begin{cases}
                1 \text{\, if \,} \mb{v} = (0, 0), \mb{r} = (1, 1,\dots,1) \\
                0 \text{\, otherwise}.
            \end{cases}
        \]
        \item triangular, with R-matrix defined for the grouplike elements as 
        \[
            R_{\matfont{z}} := \frac{1}{16}\sum_{\mb{v}, \mb{w} \in \Z^2_4} i^{-\mb{v}\mb{w}^T} \mb{K}^{\mb{w}} \otimes \mb{K}^{\mb{v} \matfont{z}}, 
        \]
        where $\mb{K} = (K_1, K_2)$, $\pmb{\xi} = (i, i)$, and $\matfont{z} =\begin{pmatrix} 2 && 3\\ 1 &&0 \end{pmatrix}$. 
        \item and ribbon, with ribbon element $v =\1$, corresponding to the unique pivotal element $g = K^2_1$, 
    \end{enumerate}
\end{proposition}
\begin{remark}
    The Hopf algebra of Example \ref{ex1} admits no non-triangular structures. The behaviour of $Z^{\pm}$ generators suggests that a non-triangular R-matrix could be built of terms of the shape $Z_l^{\pm} \otimes K_1^{\pm2} K_2^{\pm1} Z_l^{\mp}$, as the two generators are of opposite type. However, any Hopf algebra map $f_R$, sending $\mathcal{Z}_l^\pm\mapsto K_1^{-2} K_2^{-1} (Z_l^\mp)$, satisfies all constraints apart from (B1), which is equivalent to $(K_1^{2} K_2)^2=\1$, which is false in this case. The element $X_k$ similarly fails (B1) for any such Hopf algebra map. 
\end{remark}
With this observation, we can build a non-triangular example, at the cost of introducing a $K_3$, another grouplike element of order $4$. We can choose $L = K_3^2$, which implies that the corresponding $Z_l^{\pm}$ will at most anti-commute with all other nilpotent generators, according to the notation we established.

\begin{proposition}
\label{ex2prop}
    The algebra $\operatorname{N}_2$ of Example \ref{ex2} is
    \begin{enumerate}
        \item unimodular, with a two-sided cointegral
        \[
            \Lambda := \left(\sum^3_{a,b, c=0} K^a_1 K^b_2 K^c_3\right) \prod^{t_1}_{k=1}X^+_k \prod^{t_1}_{k=1}X^-_k \prod^{t_2}_{l=1}Z^+_l
            \prod^{t_2}_{l=1}Z^-_l,
        \]
        and a two-sided integral defined on the monomial basis as
        \[
            \lambda(\mb{K}^\mb{v}\mb{X}^\mb{r} ):=
            \begin{cases}
                1 \text{\, if \,} \mb{v} = (0, 0, 0), \mb{r} = (1, 1,\dots,1) \\
                0 \text{\, otherwise},
            \end{cases}
        \]
        \item quasitriangular, with the R-matrix defined for the grouplike generators as 
        \[
            R_{\matfont{z}} := \frac{1}{64}\sum_{\mb{v}, \mb{w} \in \Z^3_4} i^{-\mb{v}\mb{w}^T} \mb{K}^{\mb{w}} \otimes \mb{K}^{\mb{v} \matfont{z}}, 
        \]
        where $\mb{K} = (K_1, K_2, K_3)$ and $\matfont{z} =\begin{pmatrix} 0 && 3 && 2 \\ 1 && 0 && 0 \\ 2 && 0 && 2  \end{pmatrix}$. 
        \\The part containing skew-primitive generators is 
        \[
            R_{\pmb{\alpha}} := \exp\left(\sum^{t_2}_{l = 1} \alpha_l( Z^+_l \otimes L Z^-_l - Z^-_l \otimes L Z^+_l)\right)
        \]
        for some $\pmb{\alpha}=(\alpha_1, \dots, \alpha_{t_2}) \in \C^{2t_2}$, so that 
        \[
            R_{\matfont{z}} R_{\pmb{\alpha}} := \frac{1}{64}\sum_{\mb{v}, \mb{w} \in \Z^3_4} i^{-\mb{v}\mb{w}^T} \mb{K}^{\mb{w}} \otimes \mb{K}^{\mb{v} \matfont{z}} \exp\left(\sum^{t_2}_{l = 1} \alpha_l( Z^+_l \otimes L Z^-_l - Z^-_l \otimes L Z^+_l)\right),
        \]
        \item ribbon, with the ribbon element
        \[
            v:= \exp \left( - 2\sum^{t_2}_{l = 0} \alpha_l Z^+_l Z^-_l  \right).
        \]
        corresponding to the pivotal element $g = L =K^2_3$, which in this example is the only such element available. 
    \end{enumerate}
\end{proposition}
\noindent We can introduce more nilpotent generators while retaining the same quasitriangular structure.

\begin{proposition}
\label{ex3prop}
    The algebra $\operatorname{N}_3$ of Example \ref{ex3} is
    \begin{enumerate}
        \item unimodular, with a two-sided cointegral
        \[
            \Lambda := \left(\sum^3_{a,b, c=0} K^a_1 K^b_2 K^c_3 \right) \prod^{t_1}_{k=1}X^+_k \prod^{t_1}_{k=1}X^-_k 
            \prod^{t_1}_{k=1}Y^+_k \prod^{t_1}_{k=1}Y^-_k 
            \prod^{t_2}_{l=1}Z^+_l
            \prod^{t_2}_{l=1}Z^-_l,
        \]
        and a two-sided integral defined on the monomial basis by
        \[
            \lambda(\mb{K}^\mb{v}\mb{X}^\mb{r} ):=
            \begin{cases}
                1 \text{\, if \,} \mb{v} = (0, 0,  0), \mb{r} = (1, 1,\dots,1) \\
                0 \text{\, otherwise}.
            \end{cases}
        \]

        \item quasitriangular, with the R-matrix defined for the grouplike generators as 
        \[
            R_{\matfont{z}} := \frac{1}{64}\sum_{\mb{v}, \mb{w} \in \Z^3_4} i^{-\mb{v}\mb{w}^T} \mb{K}^{\mb{w}} \otimes \mb{K}^{\mb{v} \matfont{z}}, 
        \]
        where $\mb{K} = (K_1, K_2, K_3)$ and $\matfont{z} =\begin{pmatrix} 0 && 3 && 2 \\ 1 && 0 && 0 \\ 2 && 0 && 2  \end{pmatrix}$. 
        \\The skew-primitive part is 
        \[
            R_{\pmb{\alpha}} := \exp\left(\sum^{t_2}_{l = 1} \alpha_l( Z^+_l \otimes L Z^-_l - Z^-_l \otimes L Z^+_l)\right)
        \]
        for some $\pmb{\alpha}=(\alpha_1, \dots, \alpha_{t_2}) \in \C^{2t_2}$, so that 
        \[
            R_{\matfont{z}} R_{\pmb{\alpha}} := \frac{1}{64}\sum_{\mb{v}, \mb{w} \in \Z^3_4} i^{-\mb{v}\mb{w}^T} \mb{K}^{\mb{w}} \otimes \mb{K}^{\mb{v} \matfont{z}} \exp\left(\sum^{t_2}_{l = 1} \alpha_l( Z^+_l \otimes L Z^-_l - Z^-_l \otimes L Z^+_l)\right),
        \]
        \item ribbon, with ribbon element
        \[
            v:= \exp \left( -2 \sum^{t_2}_{k = 1} \alpha_l Z^+_l Z^-_l  \right),
        \]
        corresponding to the unique pivotal element $g = L = K^2_3$. 
    \end{enumerate}
\end{proposition}

\begin{proof}[Proof of Propositions \ref{ex1prop}, \ref{ex2prop}, \ref{ex3prop}]
The propositions can be verified by a direct check of unomodularity, quasitriangularity and ribbon axioms, or equivalently using Theorem \ref{unimodTheorem} for the cointegral, \ref{QTthm} for the R-matrix and \ref{ribbonThm} for the ribbon element.
\end{proof}
\begin{proposition}
\label{ExAreSNF}
    The Hopf algebras of examples \ref{ex1}, \ref{ex2}, \ref{ex3} admit only strongly non-factorizable quasitriangular structures.
\end{proposition}
\begin{proof}
    The theorem is verified by essentially checking all choices of bialgebra maps $f:H^* \rightarrow H^{cop}$. Indeed, since the relations are diagonal, the grouplike generators $\mathcal{K}_b$ have to be mapped to grouplike generators $K_a$ and nilpotent generators $\mathcal{X}_j, \mathcal{X}^\pm_j, \mathcal{Y}^\pm_j, \mathcal{Z}^\pm_l$ to nilpotent generators $X_k, X^\pm_k, Y^\pm_k, Z^\pm_m$. We highlight the key points.
    

    First, we explain why the examples do not admit any quasitriangular structures such that generators $X^\pm_k$, as well as  $Y^\pm_k$ in Example \ref{ex3} appear in $R_{\pmb{\alpha}}$. One can simply check all possible assignment and find that the only way to satisfy (C1) is to map 
    \[
        f(\mathcal{X}_k^\pm) = (K_1K_2)^{\mp1} X_k^\pm, \text{ resp. } f(\mathcal{Y}_k^\pm) = (K^2_1K_2)^{\mp1} Y_k^\pm,
    \]
    but it can be verified it fails (C4) for all $K_a$, $a=1, \dots, s$.
    The mapping 
    \[
        f(\mathcal{X}_k^\pm) = (K_1K_2)^{\pm1} X_k^\mp, \text{ resp. } f(\mathcal{Y}_k^\pm) = (K^2_1K_2)^{\pm1} Y_k^\mp,
    \]
    analogous to the terms involving $Z_l^\pm$, satisfies (B1) only if $(K_1K_2) = (K_1K_2)^{-1}$, which is not the case for our choice of relations. If the image of $\mathcal{X}_k^\pm$ or $\mathcal{Y}_k^\pm$ involves $Z^\pm_l$, (B1) is not satisfied, which can be checked directly. 
 
    Finally, from the definitions of the integral and the R-matrix, we observe that all nilpotent generators except $Z^\pm_l$ (or all in Example \ref{ex1}) are absent from the respective R-matrices. There is no relation available that introduces them, as the commutations are necessarily diagonal, thus they are absent from the monodromy matrix as well. Moreover, the antipode does not change the content of nilpotent generators in any monomial. Since the integral is non-zero only on the linear combinations of monomials containing the top element $T = \mb{X}^{(1, \dots, 1)}$, and since $T$ contains all the nilpotent generators, we retrieve the result. 
\end{proof}
We can also deduce the following about the ribbon structure.
\begin{corollary}
\label{CorNenciuAF}
    The Hopf algebras of examples \ref{ex1}, \ref{ex2}, \ref{ex3} admit only anomalous ribbon structures.
\end{corollary}
\begin{proof}
    In the proof of Proposition \ref{ExAreSNF}, we showed that there is a subset of nilpotent generators that cannot appear in the R-matrix, $R$. It is a straightforward consequence of the diagonality of all relations that these generators are also necessarily absent from the corresponding Drinfeld element $u = S(R'')R'$ and so is the ribbon element $v$, as in all cases the pivotal element $g\in \C[\mb{K}]$ and $v:=g^{-1}u$. Hence by the same reasoning $\lambda(v) = 0$ and so from Defintion \ref{DefAF} only anomalous ribbon structures appear. 
\end{proof}

We will later need formulas for the monodromy matrices corresponding to the class of R-matrices $R_{\matfont{z}}R_{\pmb{\alpha}}$ that appeared extensively in the examples.
\begin{proposition}
    \label{monodromyProp}
    Let $H(\mb{m}, t, \matfont{d}, \matfont{u})$ be a Hopf algebra containing nilpotent generators $Z_l^\pm$, for $l=1, \dots, t_2$. Moreover let $Z^+_l$ and $Z^-_l$ be of opposite type and $\{Z_{l_1}^\pm, Z_{l_2}^\pm\}=\{Z_{l_1}^\pm, Z_{l_2}^\mp\}=0$. Denote $L_l:=\mb{K}^{\matfont{u}_l}$  Let $R:=R_{\matfont{z}}R_{\pmb{\alpha}}$ be an R-matrix of the form as in Theorem \ref{R-matrixConstr}.
    Then the monodromy matrix is given by
    \[
        M:=R_{21}R =\exp\left(2\sum^{t_2}_{l=1} \alpha_l (Z_l^+ \otimes L_l Z_l^{-} - Z_l^- \otimes L_l Z_l^{+})\right).
    \] 
\end{proposition}
    \begin{proof}
        We first consider the behaviour of the terms $L_l Z_l^\mp \otimes Z_l^\pm$, that appear in $R_{\pmb{\alpha}, 21}$, when they commute past $R_{\matfont{z}}$.
        We have 

         \begin{align*}
            (L Z_l^\mp \otimes Z_l^\pm)R_{\matfont{z}}&=(L_l Z_l^\mp \otimes Z_l^\pm)\frac{1}{\prod^s_{a=1} m_a}\sum_{\mb{v}, \mb{w} \in \Z_\mb{m}} \pmb{\xi}^{-\mb{v}\cdot \mb{w}} \mb{K}^{\mb{w}} \otimes \mb{K}^{\mb{v}\matfont{z}} \\
            &= \frac{1}{\prod^s_{a=1} m_a}\sum_{\mb{v}, \mb{w} \in \Z_\mb{m}} \pmb{\xi}^{-\mb{v}\cdot \mb{w}+\mb{w}\cdot\matfont{d}_l - (\mb{v}\matfont{z})\cdot \matfont{d}_l} \mb{K}^{\mb{w}+\matfont{u}_l} \otimes \mb{K}^{\mb{v}\matfont{z}-\matfont{u}_l}(Z_l^\mp\otimes L_l Z_l^\pm),
        \end{align*}
        where we commuted $Z_l^\pm$ past $\mb{K}^\mb{w}$ and $\mb{K}^{\mb{v}\matfont{z}}$, absorbing the $L_l$ in $\mb{K}^\mb{w}$ and factoring it out of $\mb{K}^{\mb{v}\matfont{z}}$ in every element of the sum. Note that this is allowed because $L^2=\1$. Now using (C1), (C2) of Theorem \ref{R-matrixConstr}, we note that $(\mb{v}\matfont{z})\cdot \matfont{d}_l = \mb{v} \cdot \matfont{u}_l$. Now we introduce the shifts in the sum $\mb{w}\mapsto \mb{w}-\matfont{u}_l$ and $\mb{v}\matfont{z} \mapsto \mb{v}\matfont{z}+ \matfont{u}_l$, that is $\mb{v}\mapsto \mb{v}+\matfont{d}_l$. Now the coefficient in the sum becomes
        \[
            \pmb{\xi}^{-(\mb{v}+\matfont{d}_l)\cdot (\mb{w}-\matfont{u}_l) +(\mb{w}-\matfont{u}_l)\cdot\matfont{d}_l  - (\mb{v}\matfont{z} + \matfont{u}_l)\cdot \matfont{d}_l}
        \]
        which again by (C1), (C2) of Theorem \ref{R-matrixConstr} simplifies to 
        \[
            \pmb{\xi}^{-\mb{v}\cdot \mb{w} - \matfont{d}_l\cdot \matfont{u}_l} = \pmb{\xi}^{-\mb{v}\cdot \mb{w}}(-1),
        \]
        since by Definition \ref{genDef}, $\pmb{\xi}^{\matfont{d}_l\cdot \matfont{u}_l} = -1$.
        Thus, the commutation stands 
        \[
            \frac{1}{\prod^s_{a=1} m_a}\sum_{\mb{v}, \mb{w} \in \Z_\mb{m}} \pmb{\xi}^{-\mb{v}\cdot \mb{w}} \mb{K}^{\mb{w}} \otimes \mb{K}^{\mb{v}\matfont{z}}(-Z_l^\mp\otimes L_l Z_l^\pm).
        \]
        Now we can consider the exponential term in full to notice that 
        \[
             \exp\left(\sum^{t_2}_{l=1} \alpha_l ( L_l Z_l^{-}\otimes Z_l^+ - L_l Z_l^{+}\otimes Z_l^-)\right)R_{\matfont{z}} =\\
             R_{\matfont{z}} \exp\left(\sum^{t_2}_{l=1} \alpha_l (Z_l^+ \otimes L_l Z_l^{-} - Z_l^- \otimes L_l Z_l^{+})\right)=R_{\matfont{z}}R_{\pmb{\alpha}}.
        \]
        So we retrieve
        \[
            R_{\matfont{z}, 21}R_{\pmb{\alpha}, 21} R_{\matfont{z}} R_{\alpha} = R_{\matfont{z}, 21} R_{\matfont{z}}R_{\alpha} R_{\alpha} = 
            \\ \exp\left(2\sum^{t_2}_{l=1} \alpha_l (Z_l^+ \otimes L_l Z_l^{-} - Z_l^- \otimes L_l Z_l^{+})\right).
        \]
        \end{proof}

\section{Non-factorizable biproducts with $u_q \mathfrak{sl}_2$}

In this section we construct a family of Hopf algebras, which augment $u_q \mathfrak{sl}_2$, for $q$ a primitive $r$-th root of unity and $r \equiv 0 \mod 4$ by a Hopf algebra of Nenciu type $H = H(\mb{m}, t, \matfont{d}, \matfont{u})$ so that both pieces interact. The main property of the construction is that the result is ribbon, with the new R-matrix and ribbon element constructed from those carried by $U$ and $H$. The idea of extending $u_q \mathfrak{sl}_2$ using nilpotent generators is inspired by the work of Majid, see for instance \cite{majid_2000} Example 10.2.13, but the strategy used to find the non-factorizable examples did not appear in the literature to our best knowledge.

\subsection{The small quantum group $u_q \mathfrak{sl}_2$}
We use this section to set the conventions for the small quantum group. We follow mostly the conventions of \cite[Section 9.1]{chari_pressley_1995} and \cite[Section 5.7]{etingof_gelaki_nikshych_ostrik_2016}, used also in \cite{beliakova_derenzi_2023}. Let $r'=r/2$ and $r'' = r/4$ throughout.

\begin{definition}
\label{smallQG}
    The Hopf algebra $u_q \mathfrak{sl}_2$ where $q$, is a primitive $r$-th root of unity, $r \equiv 0 \mod 4$, is the algebra generated by elements $K, E,$ and $F$, satisfying the relations
    \[
        K^{r'} = \1, \;\; E^{r'} = F^{r'} = 0
    \]
    \[
        KEK^{-1} = q^2 E, \;\; KFK^{-1} = q^{-2} F, \;\; [E, F] = \frac{K-K^{-1}}{q-q^{-1}}. 
    \]
    The Hopf structure is given by
    \[
        \epsilon(K) = 1, \;\; \epsilon(E) = \epsilon(F) = 0,
    \]
    \[
        \Delta(K) = K \otimes K, \;\; \Delta(E) = \1 \otimes E + E \otimes K, \;\; \Delta(F) = F \otimes \1 + K^{-1} \otimes F
    \]
    and 
    \[
        S(K) = K^{-1} = K^{r'-1}, \;\; S(E) = -E K^{-1}, \;\; S(F) = -K F.
    \]
\end{definition}
\noindent Recall also the notation for quantum integers, for each $k \in \Z_{r'}$:
\[
    \{k\}:= q^k - q^{-k},\;\; [k]:=\frac{\{k\}}{\{1\}},\;\; [k]!:= [k][k-1]\dots[1].
\]
\begin{proposition}
\label{uqsl2prop}
    The Hopf algebra  $u_q \mathfrak{sl}_2$ is
    \begin{enumerate}
        \item unimodular with a two-sided cointegral
        \[
            \Lambda:=\frac{\{1\}^{r'-1}}{\sqrt{r''}[r'-1]!}\sum^{r'-1}_{a=0}  E^{r'-1} F^{r'-1} K^a,
        \]
        and a left integral
        \[
            \lambda(E^b F^c K^a) := \frac{\sqrt{r''}[r'-1]!}{\{1\}^{r'-1}} \delta_{a, r'-1}\delta_{b, r'-1}\delta_{c, r'-1},
        \]
        \item quasitriangular with $R:= D \Theta$, where
        \[
            D:= \frac{1}{r'}\sum^{r'-1}_{b, c=0} q^{-2bc} K^b \otimes K^c,
        \]
        and 
        \[
            \Theta:= \sum^{r'-1}_{a=0} \frac{\{1\}^a}{[a]!}q^{\frac{a(a-1)}{2}} E^a \otimes F^a.
        \]
        \item ribbon, with the ribbon element for $r \equiv 4 \mod 8$ 
        \[
            v = \frac{i^{\frac{r''-1}{2}}}{\sqrt{r''}} \sum^{r'-1}_{a=0}\sum^{r''-1}_{b=0} \frac{\{-1\}^a}{[a]!} q^{-\frac{(a+3)a}{2} + \frac{(r''+1)^3}{2}(2b-1)^2}  E^aF^a K^{-a-2b},
        \]
        and for $r \equiv 0 \mod 8$
        \[
            v = \frac{1-i}{\sqrt{r'}} \sum^{r'-1}_{a=0}\sum^{r''-1}_{b=0} \frac{\{-1\}^a}{[a]!} q^{-\frac{(a+3)a}{2} + 2b^2} E^aF^a  K^{-a-2b-1}.
        \]
        corresponding to the pivotal element $g=K$.
    \end{enumerate}
\end{proposition}
\begin{remark}
\label{DasRz}
    To make a connection with the notation of Theorem \ref{QTthm} for $D$, note that $q^{2r'}=1$, but $K^{r'}=\1$. Thus, we can write
    \[
        D =  \frac{1}{r'}\sum_{v, w \in \Z_{r'}} (q^2)^{-wv} K^w \otimes K^{zv}
    \]
    which corresponds to a $R_{\matfont{z}}$ for $\xi = q^2$ and $\matfont{z}=1$.
\end{remark}

\subsection{The semi-direct biproduct structure}
Towards the goal of constructing non-factorizable extensions of $u_q \mathfrak{sl}_2$ we introduce versions of \textit{semi-direct} (or \textit{smash}) product and coproduct. While these constructions are well-known in the study of Hopf algebras, see for instance \cite{molnar_1977}, \cite{andruskiewitsch2006}, we use slightly non-standard conventions to achieve the properies studied in the examples below. 
\begin{notation}
The Sweedler notation will often apear in this section when considering a semi-direct biproduct $U \ltimes H$ of two Hopf algebras $U$ and $H$. Then it refers to the corresponding coproducts $\Delta_U$ and $\Delta_H$ of the respective tensor components, rather than the new coproduct $\Delta$ of $U \ltimes H$.
\end{notation}
\begin{definition}
    \label{smashBiproduct}
    Let $U,H$ be a Hopf algebras. Let
    \[
        \lhd: H \otimes U\rightarrow H
    \]
    be a right algebra action of $U$ on $H$ (see for instance \cite{kassel_2012}, Chapter III.5). We define the \textit{right semi-direct product} of $U$ and $H$ to be the algebra $U\otimes H$ with the unit $\1_U \otimes\1_H$ and the product $\mu$ given by 
    \[
     \mu(g\otimes h, g' \otimes h') = (g\otimes h)(g' \otimes h') := g g'_{(1)} \otimes (h \lhd g'_{(2)}) h'.
    \]
   Let also
    \[
        \delta: H \rightarrow H \otimes U
    \]
    be a right coaction of $U$ on $H$ (see for instance \cite{kassel_2012} Chapter III.6). We adopt the Sweedler notation
    \[
        \delta : h\mapsto h^{[0]} \otimes h^{[1]}.
    \]
    We define the \textit{right semi-direct coproduct} of $U$ and $H$ as the coalgebra $U\otimes H$ with the counit $\epsilon_U \otimes \epsilon_H$, and the coproduct
    \[
        \Delta(g\otimes h) := (g_{(1)}\otimes h^{[0]}_{(1)})\otimes(g_{(2)} h^{[1]}_{(1)} \otimes h_{(2)}).
    \]
     Let $S:U\otimes H \rightarrow U\otimes H$ be the invertible map given by \[
        S(g \otimes h) \mapsto S_U(g h^{[1]}) \otimes S_H(h^{[0]}).
     \] If $(U \otimes H, \mu,  \1_U \otimes 1_H, \Delta, \epsilon_U \otimes \epsilon_H ,S)$ is a Hopf algebra then we will call it a \textit{semi-direct biproduct} of $U$ and $H$ and denote $U\ltimes H$ see for instance \cite[Section 6]{majid_2000}.
     
\end{definition}
\begin{remark}
    We emphasize that all tensor products between $U$ and $H$ in Definition \ref{smashBiproduct} are taken in the category $\text{\textbf{Vect}}_\C$, which carries the trivial braiding. In particular, this is a different setting than the \textit{Radford biproduct} which is taken in the braided category of Yetter-Drinfeld modules over the acting algebra, see \cite{andruskiewitsch2006}.
\end{remark}
\begin{remark}
    We do not attempt to determine the conditions under which the semi-direct product and semi-direct coproduct in Definition \ref{smashBiproduct} give rise to the semi-direct biproduct. We merely present in the sequel a construction where the requirements are indeed satisfied.
\end{remark}

\subsection{Non-factorizable biproducts with $u_q \mathfrak{sl}_2$}
In this section we build biproducts of $u_q \mathfrak{sl}_2$ and $H(\mb{m}, t, \matfont{d}, \matfont{u})$, and we show they can have the properties of unimodularity, quasitriangularity and ribbonness, and can be (strongly) non-factorizable. For clarity of discussion and the applications we keep in mind, we restrict the order of $q$, $r$ to be a multiple of 8.
\begin{definition}
\label{extHopfDef}
    Let $U= u_q \mathfrak{sl}_2$ and $H = H(\mb{m}, t, \matfont{d}, \matfont{u})$, where $q$, is a primitive $r$-th root of unity, $r \equiv 0 \mod 8$ and $r'=r/2$, $r''=r/4$. We define $U \ltimes H$, where
    \begin{itemize}
        \item the right action $\lhd: H(\mb{m}, t, \matfont{d}, \matfont{u}) \otimes  u_q \mathfrak{sl}_2 \rightarrow H(\mb{m}, t, \matfont{d}, \matfont{u})$ for $K_a \in H$, $a=1, \dots, s$ by
         \[
             K_a \lhd K = K_a, \;\; K_a \lhd E = K_a \lhd F = 0 
         \]
        and for $X_k \in H$ by
        \[
            X_k \lhd K = -X_k, \;\; X_k \lhd E = X_k \lhd F = 0, 
        \]
        \item the right coaction $\delta: H(\mb{m}, t, \matfont{d}, \matfont{u}) \rightarrow H(\mb{m}, t, \matfont{d}, \matfont{u}) \otimes  u_q \mathfrak{sl}_2$  for $K_a \in H$ by
        \[
            K_a \mapsto K_a \otimes \1_U,
        \]
        for $X_k \in H$ by
        \[
            X_k \mapsto X_k \otimes K^{r''},
        \]
        and on an arbitrary monomial $\mb{K}^\mb{v} \mb{X}^\mb{r} \in H$ for $\mb{v}\in \Z_\mb{m}$ and $\mb{r}\in \Z^t_2$ by
        \[
            \mb{K}^\mb{v} \mb{X}^\mb{r} \mapsto \mb{K}^\mb{v} \mb{X}^\mb{r} \otimes K^{|\mb{r}|r''}.
        \]
    \end{itemize}
    
\end{definition}
\begin{theorem}
\label{extHopfThm}
    Let $u_q \mathfrak{sl}_2 \ltimes H(\mb{m}, t, \matfont{d}, \matfont{u})$ be as in Definition \ref{extHopfDef}. Then it is a Hopf algebra with the following semi-direct biproduct structure:
    \begin{itemize}
        \item the new unit $\1 = \1_U \otimes \1_H$ and counit $\epsilon := \epsilon_U \otimes \epsilon_H$,
        \item the algebra structure where we suppress the tensor product structure 
    \begin{align*}
       K:= K \otimes \1_H && E:= E \otimes \1_H && F:= F \otimes \1_H && 
    \end{align*}
    \begin{align*}
        K_a := \1_U \otimes K_a  && X_k := \1_U \otimes X_k,
    \end{align*}
    so that for any $a, b, c, \in \Z_{r'}$, $\mb{v}\in \Z_\mb{m}$ and $\mb{r}\in \Z^t_2$
    \[
        E^a F^b K^c \mb{X}^\mb{r}\mb{K}^\mb{v}:= E^a F^b K^c \otimes \mb{X}^\mb{r}\mb{K}^\mb{v}
    \]
    and find the new algebra relations between the the elements of $U$ and $H$ to be 
    \begin{align*}
        [K, K_a] = [E, K_a] = [F, K_a] = 0, && \{K, X_k\} = \{E, X_k\} = [F, X_k]=0. 
    \end{align*} 
    \item new Hopf structure for $X_k$ give by
    \begin{align*}
        \Delta(X_k) = \1 \otimes X_k + X_k \otimes K^{r''}\mb{K}^{\matfont{u}_k} && S(X_k)  =  -X_k K^{r''}\mb{K}^{-\matfont{u}_k},
    \end{align*}
    while for the remaining generators $E, F, K$ and $K_a$ in this notation the coproduct $\Delta$ and antipode $S$ are unchanged.
    \end{itemize}
\end{theorem}
\begin{proof}
    See Section \ref{proof:extHopfThm} in Appendix A.
\end{proof}

\begin{proposition}
\label{extMonomialBasis}
The Hopf algebra $u_q \mathfrak{sl}_2 \ltimes H(\mb{m}, t, \matfont{d}, \matfont{u})$ has a monomial basis
    \[
       \{E^e F^f K^k \mb{K}^\mb{w}\mb{X}^\mb{r}|e, f, k =1, \dots, r', \mb{w} \in \Z_\mb{m}, \mb{r} \in \Z^t_2\}.
    \]
\end{proposition}
\begin{proof}
    This follows directly from the fact that both $U$ and $H$ have monomial bases, by the Poincar\'e-Birkhoff-Witt theorem (see for instance \cite{kassel_2012}) and Proposition \ref{monomialBasis}, respectively. Indeed, all newly introduced relations between the two bases are diagonal by Theorem \ref{extHopfDef}. Hence, by the same argument as in the proof of Proposition \ref{monomialBasis}, we can establish the existence of the monomial basis.
\end{proof}
\begin{remark}
\label{SESproposition}
    Let $U,H$ and $U\ltimes H$ be as in Definition \ref{extHopfDef}, in particular let the pivotal element of $H$ be $g_H$. Then we claim the following maps are Hopf algebra maps, which do \emph{not} give a short exact sequence:
    \[
        U\hookrightarrow U\ltimes H \twoheadrightarrow H.
    \]
    The map $\iota : U\hookrightarrow U\ltimes H$ is given by $u \mapsto u \otimes \1_H$ for $u \in U$, and has a left inverse $j: U \ltimes H\twoheadrightarrow U$ given by
    \begin{align*}
        j(K) = K, && j(E) = E, && j(F) = F, && j(K_a) = \1_U, && j(X_k) = 0. 
    \end{align*}
    The map $p : U\ltimes H \twoheadrightarrow H$, is given by 
    \begin{align*}
        p(K) = g_H, && p(E) = p(F) = 0, && p(K_a) = K_a, && p(X_k) = X_k.  
    \end{align*}
\end{remark}
\begin{remark}
\label{NonSplit}
    The map $p$ in Remark \ref{SESproposition} has a right inverse as an algebra map, given by the map $$q :H \hookrightarrow H\ltimes U, \;\; q(h)=\1_U\otimes h,$$ but it is \textit{not} a Hopf algebra map. For instance $S(X_k) = K^{r''} \mb{K}^{\matfont{u}_k}X_k$, so the image of $H$ under $q$ is not closed under taking antipodes (as it is not under taking coproducts).
    Nevertheless, it makes sense to speak, for instance, of elements $u, h \in U\ltimes H$ as belonging to $u \in U$ and $h\in H$, and we will often slightly abuse the notation in this way.
\end{remark}
\begin{remark}
    Let $U,H$ and $U\ltimes H$ be as in Theorem \ref{extHopfThm}. Then in light of Remark \ref{SESproposition} $U\ltimes H$ is not isomorphic to $U\otimes H$ as Hopf algebra in the usual sense. In particular, $U\otimes H$ has no non-trivial commutation relations between the elements belonging to $U$ and $H$. Similarly, $U\ltimes H$ is not isomorphic to $U\oplus H$, that is the short exact sequence is not fully split, as per Remark \ref{NonSplit}.
\end{remark}

Now we state the main theorem of the section.

\begin{theorem}
\label{extMainThm}
    Let $U= u_q \mathfrak{sl}_2$ where $q$ is a primitive $r$-th root of unity, $r\equiv 0 \mod 8$ as in Definition \ref{smallQG} carrying the integral $\lambda_U$, the cointegral $\Lambda_U$, the R-matrix $D\Theta$ and ribbon element $v_U$ corresponding to a pivotal element $g_U$, as in Proposition \ref{uqsl2prop}.
    Let also $H = H(\mb{m}, t, \matfont{d}, \matfont{u})$ with
    \begin{enumerate}
        \item a two-sided cointegral
        \[
            \Lambda_H := \prod^s_{a=1} \left(\sum^{m_a-1}_{b=0} K^b_a \right) \prod^t_{k=1} X_k,
        \]
        \item a two-sided integral
        \[
            \lambda_H(\mb{K}^\mb{v} \mb{X}^\mb{r}):=
            \begin{cases}
                1 \text{\, if \,}\mb{K}^\mb{v} \mb{X}^\mb{r}  = T \\
                0 \text{\, otherwise},
            \end{cases}
        \]
        \item an R-matrix $R_H := R_{\matfont{z}}R_{\pmb{\alpha}}$ where
        \[
            R_{\matfont{z}}=\frac{1}{\prod^s_{a=1} m_a}\sum_{\mb{v}, \mb{w} \in \Z_\mb{m}} \pmb{\xi}^{-\mb{v}\cdot \mb{w}} \mb{K}^{\mb{w}} \otimes \mb{K}^{\mb{v}\matfont{z}}
        \]
        \[
            R_{\pmb{\alpha}} = \exp\left( \sum^{t_2}_{l=0} \alpha_l(Z_l^+ \otimes L_lZ_l^- - Z_l^- \otimes L_lZ_l^+) \right),
        \]
        where $L_l := \mb{K}^{\matfont{u}_l}$ 
        \item a ribbon element
        \[
            v_H:= \exp \left( -2 \sum^{t_2}_{k = 1} \alpha_l Z^+_l Z^-_l  \right)
        \]
        and the corresponding pivotal element $g$ and the Drinfeld element
        \[
            u_H:= g^{-1}\exp \left( -2 \sum^{t_2}_{k = 1} \alpha_l Z^+_l Z^-_l  \right).
        \]
    \end{enumerate}
    Then, $U\ltimes H$, for the choices of action and coaction as in Definition \ref{extHopfDef}, carries the following structures
    \begin{enumerate}
        \item a two-sided cointegral
        \[
            \Lambda:=\Lambda_U \Lambda_H = \frac{\{1\}^{r'-1}}{\sqrt{r''}[r'-1]!}\sum^{r'-1}_{a=0}  E^{r'-1} F^{r'-1} K^a 
             \prod^s_{a=1} \left(\sum^{m_a-1}_{b=0} K^b_a \right) \prod^t_{k=1} X_k,
        \]
        \item a left integral $\lambda\in (U\ltimes H)^*$, $\lambda:=\lambda_U \otimes \lambda_H$, equivalently defined on the monomial basis by
        \[
             \lambda(E^b F^c K^a\mb{K}^\mb{v}\mb{X}^\mb{r}) :=
            \begin{cases}
                \frac{\sqrt{r''}[r'-1]!}{\{1\}^{r'-1}} \ \text{\, if \,} \mb{v} = (0, \dots, 0), \mb{r} = (1,\dots,1), a=b=c=r'-1 \\
                0 \text{\, otherwise},
            \end{cases}
        \]        \item an R-matrix $R:=DR_{\matfont{z}}\Theta\Bar{R}_{\pmb{\alpha}}$, where 
        \[
            \Bar{R}_{\pmb{\alpha}}:= \exp\left( \sum^{t_2}_{l=0} \alpha_l(Z^+ \otimes \Bar{L}_l Z_l^- - Z_l^- \otimes \Bar{L}_l Z_l^+) \right)
        \]
        for $\Bar{L}_l = K^{r''} L_l$,
        \item a ribbon element 
        \begin{multline*}
            v:=v_U u_H = \\
            \frac{1-i}{\sqrt{r'}} \sum^{r'-1}_{a=0}\sum^{r''-1}_{b=0} \frac{\{-1\}^a}{[a]!} q^{-\frac{(a+3)a}{2} + 2b^2} E^aF^a  K^{-a-2b-1} \prod^{t_2}_{l=1} L_l\exp \left( -2 \sum^{t_2}_{k = 1} \alpha_l Z^+_l Z^-_l  \right)
        \end{multline*}
        with the pivotal element $g= K$, and the Drinfeld element
        \begin{multline*}
             u:=u_U u_H = \\
            \frac{1-i}{\sqrt{r'}} \sum^{r'-1}_{a=0}\sum^{r''-1}_{b=0} \frac{\{-1\}^a}{[a]!} q^{-\frac{(a+3)a}{2} + 2b^2} E^aF^a  K^{-a-2b} \prod^{t_2}_{l=1} L_l\exp \left( -2 \sum^{t_2}_{k = 1} \alpha_l Z^+_l Z^-_l  \right).
        \end{multline*}
    \end{enumerate}
\end{theorem}

\begin{proof}
The proof is rather technical and lenghty, so we postpone it to Section \ref{proof:extMainThm} in Appendix A.
\end{proof}

\begin{proposition}
\label{extNonFact}
    Let $u_q \mathfrak{sl}_2 \ltimes H(\mb{m}, t, \matfont{d}, \matfont{u})$ be as above. If $H(\mb{m}, t, \matfont{d}, \matfont{u})$ carries a strongly non-factorizable ribbon structure, then it induces a strongly non-factorizable ribbon structure on  $u_q \mathfrak{sl}_2 \ltimes H(\mb{m}, t, \matfont{d}, \matfont{u})$.
\end{proposition}
\begin{proof}
    By the reasoning of Proposition \ref{ExAreSNF}, if $H = H(\mb{m}, t, \matfont{d}, \matfont{u})$ is strongly non-factorizable, the top element $\mb{X}^{(1, \dots, 1)}$ is missing from the monodromy matrix $M_H \in H\otimes H$. Since all new relations are also diagonal, and as before the antipode does not change the content of nilpotent generators in any monomial, it will be missing from the new monodromy matrix $M \in (U\ltimes H)\otimes(U\ltimes H)$, so by the same argument the result follows.
\end{proof}

\subsection{Examples}
As in the preceding section, we restrict to $r\equiv 0 \mod 8$, $r\geq 8$ so we retain a non-trivial relation of $[E, F]$ in the $u_q \mathfrak{sl}_2$. We start with an example related to the Hopf algebra $\operatorname{SF}_2$ (the version involving $\operatorname{SF}_{2n}$ is analogous), but one which generically admits non-strongly non-factorizable quaistriangular structures.

\begin{example}
    Let $u_q \mathfrak{sl}_2\ltimes H$, where $H = \operatorname{SF}_2$ be a Hopf algebra generated by $K, E, F$ with relations and morphisms of Definition \ref{smallQG}, as well as $L$, $Z^\pm$ with the following relations
    \[
        L^2= \1, \;\;(Z^\pm)^2 = 0, \;\; 
    \]
    \[
        K Z^\pm K^{-1} =q^{r'}Z^\pm = -Z^\pm, \;\; L Z^{\pm} K^{-1}_2 =- Z^\pm 
    \]
    \[
        [L, K] = [L, E] = [L, F]=0
    \]
    \[
        EZ^\pm = q^{r'} Z^\pm E= - Z^\pm E, \;\; [Z^\pm, F]=0. 
    \]
    Let $\Bar{L}:=K^{r''} L$ for convenience. The Hopf structure is 
    \[
        \epsilon(L) = 1, \;\; \epsilon(Z^\pm) = 0, 
    \]
    \[
        \Delta(L) = L \otimes L, \;\;  \Delta(Z^\pm) = \1 \otimes Z^\pm + Z^\pm \otimes \Bar{L}^{\pm 1}.
    \]
    and
    \[
        S(L) = L^{-1}, \;\; S(Z^\pm) = - Z^\pm \Bar{L}^{\mp 1}.
    \]
\end{example}
\noindent The smallest instance with a non-trivial commutation relation $[E, F]$ occurs at $r=8$. The dimension is then $4^3\times 2^3  = 2^{8}=512$.
\begin{proposition}
    The algebra $u_q \mathfrak{sl}_2 \ltimes \operatorname{SF}_2$ is 
    \begin{enumerate}
        \item unimodular with a two-sided cointegral
        \[
            \Lambda:=\frac{\{1\}^{r'-1}}{8\sqrt{r''}[r'-1]!}\sum^{r'-1}_{a=0}\sum^2_{b}   E^{r'-1} F^{r'-1}K^a L^b Z^+ Z^-,
        \]
        and a left integral expressed on the monomial basis by 
        \[
            \lambda(  E^e F^f K^a L^b (Z^+)^g (Z^-)^h) :=
            \\ \frac{\sqrt{r''}[r'-1]!}{\{1\}^{r'-1}} \delta_{a, r'-1}\delta_{b, 0}\delta_{e, r'-1}\delta_{f, r'-1}\delta_{g, 1}\delta_{h, 1},
        \]
        \item quasitriangular, with the R-matrix $R:=R_{\matfont{z}}D\Theta \Bar{R}_{\pmb{\alpha}}$, where $D, \Theta$ were defined in Proposition \ref{uqsl2prop},  $R_{\matfont{z}}$ in Proposition \ref{SFprop}, and
        \[
            \Bar{R}_{\pmb{\alpha}} := \exp\left( \alpha( Z^+ \otimes \Bar{L} Z^- - Z^- \otimes \Bar{L} Z^+)\right)
        \]
        for $\Bar{L} = K^{r''}L$,
        \item ribbon, with the ribbon element
        \[
            v = \frac{1-i}{\sqrt{r'}} \sum^{r'-1}_{a=0}\sum^{r''-1}_{b=0} \frac{\{-1\}^a}{[a]!} q^{-\frac{(a+3)a}{2} + 2b^2}  E^aF^a K^{-a-2b-1} L\exp \left( -2 \alpha Z^+ Z^-  \right)
        \]
        corresponding to the pivotal element $g=K$.
    \end{enumerate}
\end{proposition}
\noindent The presence of $R_{\pmb{\alpha}}$ can lead to $\lambda(M')\otimes M'' \neq 0$ being non-zero. In particular, the resulting expression is dependent on $\alpha$ and turns out to be non-zero as long as $\alpha\neq 0$. Thus, $u_q \mathfrak{sl}_2\ltimes \operatorname{SF}_2$ is not generically strongly non-factorizable. We now construct an example admitting only strongly non-factorizable quasitriangular structures. To do it, we introduce one more Nenciu-type algebra, supressing the Nenciu data for brevity.
\begin{example}
\label{auxExample}
    Let $\operatorname{N}_4:=H(\mb{m}, t, \matfont{d}, \matfont{u})$ be a Nenciu type algebra generated by grouplike $K_1, K_2$ and nilpotent $X^\pm$ generators, and the following relations
    \[
        K^4_1 = K^4_2 = \1, \;\;(X^\pm)^2 = 0, \;\; 
    \]
    \[
        K_1 X^{\pm} K^{-1}_2 = K_2 X^{\pm} K^{-1}_2 = \pm i X^\pm. 
    \]
    Let $L:= K_1 K_2$ for convenience. The Hopf structure is 
    \[
        \epsilon(K_1) = \epsilon(K_2)=1, \;\; \epsilon(X^\pm) = 0, 
    \]
    \[
        \Delta(K_1) = K_1 \otimes K_1, \;\; \Delta(K_2) = K_2 \otimes K_2, \;\; \Delta(X^\pm) = \1 \otimes X^\pm + X^\pm \otimes L^{\pm 1}
    \]
    and
    \[
        S(K_1) = K_1^{-1}, \;\; S(K_2) = K_2^{-1}, \;\; S(X^\pm) = - X^\pm L^{\mp 1}.
    \]
\end{example}
\begin{proposition}
    The Hopf algebra $\operatorname{N}_4$ of Example \ref{auxExample} is
    \begin{enumerate}
        \item unimodular, with two-sided cointegral
        \[
            \Lambda := \left(\sum^3_{a,b=0} K^a_1 K^b_2 \right) X^+X^- 
        \]
        and a two-sided integral expressed on the monomial basis by
        \[
            \lambda(\mb{K}^\mb{v}\mb{X}^\mb{r} ):=
            \begin{cases}
                1 \text{\, if \,} \mb{v} = (0, 0), \mb{r} = (1, 1) \\
                0 \text{\, otherwise}.
            \end{cases}
        \]
        \item triangular, with the R-matrix defined for the grouplike elements as 
        \[
            R_{\matfont{z}} := \frac{1}{16}\sum_{\mb{v}, \mb{w} \in \Z^2_4} i^{-\mb{v}\mb{w}^T} \mb{K}^{\mb{w}} \otimes \mb{K}^{\mb{v} \matfont{z}}, 
        \]
        where $\mb{K} = (K_1, K_2)$, $\pmb{\xi} = (i, i)$, and $\matfont{z} =\begin{pmatrix} 2 && 3\\ 1 && 0 \end{pmatrix}$. 
        \item and trivially ribbon, with the ribbon element $\1$, corresponding to the unique pivotal element $g = K^2_1$. 
    \end{enumerate}
\end{proposition}
\begin{remark}
\label{H4Rem}
    Note $\operatorname{N}_4\subset \operatorname{N}_1$ of Example \ref{ex1}. In analogy with $\operatorname{N}_1$ it can be checked directly from Theorem \ref{QTthm} that the Hopf algebra $\operatorname{N}_4$ admits only triangular quasitriangular structures.
\end{remark}
\begin{example}
\label{extexample1}
    Let $u_q \mathfrak{sl}_2 \ltimes \operatorname{N}_4$, be the Hopf algebra generated by $K, E, F$ with relations and morphisms of Definition \ref{smallQG}, as well as $K_1, K_2$, $X^\pm$ with the following relations
    \[
        K^4_1 = K^4_2 = \1, \;\;(X^\pm)^2 = 0, \;\; 
    \]
    \[
        K X^\pm K^{-1} =q^{r'}X^\pm = -X^\pm, \;\; K_1 X^{\pm} K^{-1}_2 = K_2 X^{\pm} K^{-1}_3 = \pm i X^\pm 
    \]
    \[
        [K_1, K_2] = [K_1, K] = [K_1, E] = [K_1, F] = [K_2, K] = [K_2, E]=[K_2, F]=0
    \]
    \[
        EX^\pm = q^{r'} X^\pm E= - X^\pm E, \;\; [X^\pm, F]=0. 
    \]
    Let $\Bar{L}:=K^{r''} K_1 K_2$ for convenience. The Hopf structure is 
    \[
        \epsilon(K_1) = \epsilon(K_2)=1, \;\; \epsilon(X^\pm) = 0, 
    \]
    \[
        \Delta(K_1) = K_1 \otimes K_1, \;\; \Delta(K_2) = K_2 \otimes K_2, \;\; \Delta(X^\pm) = \1 \otimes X^\pm + X^\pm \otimes \Bar{L}^{\pm 1}
    \]
    and
    \[
        S(K_1) = K_1^{-1}, \;\; S(K_2) = K_2^{-1}, \;\; S(X^\pm) = - X^\pm \Bar{L}^{\mp 1}.
    \]
\end{example}

\noindent The smallest instance with a non-trivial commutation relation $[E, F]$ occurs at $r=8$. The dimension is then $4^3\times 2^4\times2^2 = 2^{12}=4096$. Also, $X^\pm$ share many features with $X_a$ of Example \ref{ex1}. We have to modify the unimodular and ribbon structure accordingly. 
\begin{proposition}
\label{extProp1}
    The algebra $u_q \mathfrak{sl}_2 \ltimes \operatorname{N}_4$ is 
    \begin{enumerate}
        \item unimodular with a two-sided cointegral
        \[
            \Lambda:=\frac{\{1\}^{r'-1}}{\sqrt{r''}[r'-1]!}\sum^{r'-1}_{a=0}\sum^4_{b, c = 0}   E^{r'-1} F^{r'-1}K^a K_1^b K_2^c X^+ X^-,
        \]
        and a left integral expressed on the monomial basis by 
        \[
            \lambda(  E^e F^f K^a K_1^b K_2^c (X^+)^g (X^-)^h) =
            \\ \frac{\sqrt{r''}[r'-1]!}{\{1\}^{r'-1}} \delta_{a, r'-1}\delta_{b, 0}\delta_{c, 0}\delta_{e, r'-1}\delta_{f, r'-1}\delta_{g, 1}\delta_{h, 1},
       \]
        \item quasitriangular, with the R-matrix $R:=R_{\matfont{z}}D\Theta$, where $D, \Theta$ were defined in Proposition \ref{uqsl2prop}, and 
        \[
            R_{\matfont{z}} := \frac{1}{16}\sum_{\mb{v}, \mb{w} \in \Z^2_4} i^{-\mb{v}\mb{w}^T} (K_1, K_2)^{\mb{w}} \otimes (K_1, K_2)^{\mb{v} \matfont{z}}, 
        \] where $\matfont{z} =\begin{pmatrix} 2 && 3\\ 1 && 0 \end{pmatrix}$,
        \item ribbon, with the ribbon element
        \[
            v = \frac{1-i}{\sqrt{r'}} \sum^{r'-1}_{a=0}\sum^{r''-1}_{b=0} \frac{\{-1\}^a}{[a]!} q^{-\frac{(a+3)a}{2} + 2b^2}  E^aF^a K^{-a-2b-1} K^2_2.
        \]
        corresponding to the pivotal element $g=K$.
    \end{enumerate}
\end{proposition}
\begin{proof}
    See Section \ref{proof:extProp1} in Appendix A.
\end{proof}
\begin{proposition}
\label{SNFextension}
    The Hopf algebra $u_q \mathfrak{sl}_2\ltimes H$ carries a strongly non-factorizable ribbon Hopf algebra, and does not admit any non-strongly non-factorizable quasitraingular structures.
\end{proposition}
\begin{proof}
As pointed out in Remark \ref{H4Rem}, applying the argument of Proposition \ref{ExAreSNF} to $\operatorname{N}_4$, we see can easily see it admits only strongly non-factorizable (in fact triangular) quasitriangular structures. Then by Proposition \ref{extNonFact} we find that the quasitriangular structures induced on $u_q \mathfrak{sl}_2\ltimes H$ are necessarily strongly non-factorizable.
    
Thus, we aim to use Theorem \ref{classicalTheorem} to rule out quasitriangular structures not induced in the above way, that is to show the are no bialgebra maps $f: (u_q \mathfrak{sl}_2\ltimes H) ^* \rightarrow (u_q \mathfrak{sl}_2\ltimes H )^{cop}$ that assign non-zero values to $\mathcal{X}^\pm$, the image of $X^\pm$ under duality as given by Proposition 
\ref{dualProp} (in fact $\operatorname{N}_4$ is found to be self-dual so we drop the calligraphic notation in this proof), as well as ones containing $X^\pm$ in their image. We first understand which map produces the $D\Theta$ part of the R-matrix, and as we already established that the $H$ piece contributes only triangular structures, it suffices rule out any potential cross-terms. It is enough to verify the claims on the relevant generators, as any bialgebra map needs to, in particular, respect the monomial basis and all Hopf morphisms.
\begin{enumerate}
    \item  It is a well-known fact (see for instance \cite{kassel_2012}) that the map in Theorem \ref{classicalTheorem} producing the R-matrix of $u_q \mathfrak{sl}_2$ in Proposition \ref{QTthm} is the Cartan involution $f: (u_q \mathfrak{sl}_2)^* \rightarrow (u_q \mathfrak{sl}_2)^{cop}$  
    \begin{align*}
        f(E)=F, && f(F) = E, && f(K)=K^{-1}.
    \end{align*}
    Crucially, $u_q \mathfrak{sl}_2$ is self-dual, as explained in \cite{beliakova_derenzi_2023} and we do not introduce any new notation for the duals of $E$ and $F$. The co-opposite algebra is given by the image of this map, the $E$ and $F$ generators are exchanged and $K$ inverted. Since by Proposition \ref{dualProp} $\operatorname{N}_4$ is also self-dual, we readily verify $(U\ltimes H)^*\cong U\ltimes H$, as Hopf algebras. 


    \item Firstly, we need to exclude assignments of the form $X^\pm \mapsto b E +c F$, for some constants $b, c \in \C$. Since $f$ needs to be an algebra map we seek $ \Delta(f(X^\pm)^2) = 0$. We find 
    \[
       f(X^\pm)^2 = (bE+cF)^2 = b^2 E^2 + c^2 F^2 +bc(EF+FE) \neq 0,
    \]
    unless $b=c=0$. Hence, necessarily $f(X^\pm)=0$.
   \item Secondly, we need to exclude assignments of the form $E\mapsto F + a L^{\mp 1} X^\pm$ or $F\mapsto E + a' L^{\mp 1} X^\pm$ for $a, a' \in \C$. Consider 
    \[
        \Delta(E)  = \1 \otimes E + E\otimes K 
    \]
    against
    \[
        \Delta^{cop}(F) + \Delta^{cop}(L^{\mp1}X^\pm)=\\
        \1\otimes (F+L^{\mp 1} X^\pm) + F\otimes K^{-1} + a L^{\mp 1} X^\pm \otimes  L^{\mp 1}. 
    \]
    We need $f$ to be a coalgebra map so
    \[
    f(\Delta(E)) = f(\1)\otimes f(E) + f(k) \otimes f(E),
    \]
    which requires,
    \[
    F\otimes K^{-1} + a L^{\mp 1} X^\pm \otimes  L^{\mp 1} = E \otimes f(K)
    \]
    so $a=0$ and $f(K) = K^{-1}$ or $a=1$ and $f(K) = K^{-1} + L^{\pm 1}$. But in the latter case $f$ is not a coalgebra map as the image of the grouplike $K$ is the non-grouplike $K^{-1} + L^{\pm 1}$. Thus, $a=1$, and $f(K) = K^{-1}$ and $f(E) = F$. The assignment of $f(F) = E$ is fixed by the same argument and we retrieve the assignments of the Cartan involution.
\end{enumerate}
Thus, we see the top monomial picked up by the integral of $U\ltimes H$ is absent from the monodromy matrix, as in Proposition \ref{ExAreSNF}. We conclude the proof with the remark that we need not ask about the fate of generators $K_a$ under the map $f$, as this does not affect the strong non-factorizability condition in this case.
\end{proof}
\begin{corollary}
\label{CorBiprodAF}
    The Hopf algebra of example \ref{extexample1} admits only anomalous ribbon structures.
\end{corollary}
\begin{proof}
    Analogous to the proof of Corollary \ref{CorNenciuAF}.
\end{proof}

It is also possible to include the exponential piece of the R-matrix as in Example \ref{ex2}, while retaining strong non-factorizability.

\begin{example}
\label{extexample2}
    Let $u_q \mathfrak{sl}_2\ltimes \operatorname{N}_2$, be the Hopf algebra generated by $K, E, F$ with relations and morphisms of def. \ref{smallQG}, as well as $K_1, K_2, K_3$, $X^\pm$, $Z^\pm$ with the following relations, for $a=1, 2, 3$
    \[
      [K_a, K] = [K_a, E] = [K_a, F]=0,\;\;  K X^\pm K^{-1} = -X^\pm, \;\;K_a^4 = \1, \;\;  \;\; K_a X^{\pm} = \pm i X^{\pm} K_a
    \]
    \[
         K Z^\pm K^{-1} = -Z^\pm\;\;K_1 Z^{\pm} K^{-1}_2= Z^{\pm} , \;\; K_2 Z^{\pm} K^{-1}_3= - Z^{\pm}, \;\;  K_3 Z^{\pm} K^{-1}_4= \pm i Z^{\pm}
    \]
    \[
        \{X^{\pm}, X^{\pm}\} = \{X^{\pm}, X^{\mp}\} = \{Z^{\pm}, X^{\pm}\} = \{Z^{\pm}, X^{\mp}\} = \{Z^{\pm}, Z^{\pm}\}=\{Z^{\pm}, Z^{\mp}\} =  0
    \]
    \[
        EX^\pm = - X^\pm E,\;\; EZ^\pm = - Z^\pm E, \;\; [X^\pm, F]=[Z^\pm, F]=0. 
    \]
    Let also $\Bar{L}:=K^{r''} K_3^2$ as a shorthand. The Hopf structure is defined by
    \[
        \epsilon(K_a) = 1, \;\; \epsilon(X^{\pm})=\epsilon(Z^{\pm}) = 0,
    \]
    \[
        \Delta(K_a) = K_a \otimes K_a, \;\; \Delta(X^{\pm}) = \1 \otimes X^{\pm} + X^{\pm} \otimes (K^{r''}K_1 K_2)^{\pm 1}, \;\; \Delta(Z^{\pm}) = \1 \otimes Z^{\pm} + Z^{\pm} \otimes \Bar{L}
    \]
    \[
       S(K_a)=K^{-1}_a \;\; S(X^{\pm}) = -X^{\pm} (K^{r''} K_1 K_2)^{\mp1} \;\;  S(Z^{\pm}) = -Z^{\pm} \Bar{L}.
    \]
\end{example}
\noindent At $r=8$ the dimension is $4^4\times 2^4\times2^2 \times 2^2= 2^{16}=65536$. 
\begin{proposition}
\label{extprop2}
   The algebra $u_q \mathfrak{sl}_2\ltimes \operatorname{N}_2$ is 
    \begin{enumerate}
        \item unimodular with a two-sided cointegral
        \[
            \Lambda:=\frac{\{1\}^{r'-1}}{\sqrt{r''}[r'-1]!}\sum^{r'-1}_{a=0}\sum^4_{b, c, d= 0}   E^{r'-1} F^{r'-1} K^a K_1^b K_2^c K_3^dX^+ X^-Z^+ Z^-,
        \]
        and a left integral expressed on the monomial basis by
        \begin{align*}
             &\lambda(E^e F^f K^a K_1^b K_2^c K^d_4 (X^+)^g (X^-)^h  (Z^+)^i (Z^-)^j)  \\=&\frac{\sqrt{r''}[r'-1]!}{\{1\}^{r'-1}} \delta_{a, r'-1}\delta_{b, 0}\delta_{c, 0}\delta_{d, 0}\delta_{e, r'-1}\delta_{f, r'-1}\delta_{g, 1}\delta_{h, 1}\delta_{a, 1}\delta_{j, 1}.
        \end{align*}
        \item quasitriangular, with the R-matrix $R:=R_{\matfont{z}}D\Theta R_{\pmb{\alpha}}$, where $D, \Theta$ were defined in Proposition \ref{uqsl2prop}, and 
        \[
            R_{\matfont{z}} := \frac{1}{64}\sum_{\mb{v}, \mb{w} \in \Z^3_4} i^{-\mb{v}\mb{w}^T} (K_1, K_2, K_3)^{\mb{w}} \otimes (K_1, K_2, K_3)^{\mb{v} \matfont{z}}, 
        \] for $\matfont{z} =\begin{pmatrix} 0 && 3 && 2 \\ 1 && 0 && 0 \\ 2 && 0 && 2  \end{pmatrix}$,
        and
        \[
            R_{\alpha} := \exp\left(\alpha( Z^+ \otimes \Bar{L} Z^- - Z^- \otimes \Bar{L} Z^+)\right).
        \]
        \item ribbon, with the ribbon element
        \[
            v = \frac{1-i}{\sqrt{r'}} \sum^{r'-1}_{a=0}\sum^{r''-1}_{b=0} \frac{\{-1\}^a}{[a]!} q^{-\frac{(a+3)a}{2} + 2b^2} K^{-a-2b-1} K^2_3 E^aF^a\exp \left( - 2 \alpha Z^+ Z^-  \right).
        \]
        corresponding to the pivotal element $g=K$.
    \end{enumerate}
\end{proposition}
\begin{proof}
    See \ref{proof:extprop2} in Appendix A.
\end{proof}

\section{Topological applications}
\label{sectopology}
As mentioned in the introduction, the main motivation to study non-factorizable ribbon Hopf algebras is that they produce invariants of 4-dimensional 2-handlebodies as described by Beliakova and De Renzi in \cite{beliakova_derenzi_2023}. More details can be found in \cite{beliakova_bobtcheva_derenzi_piergallini_2023}, which relies partly on \cite{bobtcheva_piergallini_2006}. We first introduce some definitions, and then recall some of the results we mentioned. For a precise introduction to the theory of handles and handlebodies, and more precise definitions, we refer to \cite[Sections 4.1 and 4.2]{gompf_stipsicz_1999}.
\subsection{Handlebodies}

Recall that a \textit{$n$-dimensional $k$-handle} is a a copy of the $4$-ball $D^k \times D^{n-k}$ destined to be attached to another manifold along $(\partial D^k) \times D^{n-k}$. The \emph{core} of the handle is a $k$-cell $D^k \times \lbrace 0 \rbrace$, and the handle is a thickening of it.
\begin{definition}
A \textit{$4$-dimensional $2$-handlebody} is a smooth manifold $W$ together with a filtration
 $$W_{0} \subset W_1 \subset W_2,$$ where $W_0$ is a finite disjoint union of $4$-balls, $W_1$ is obtained from $W_0$ by gluing $1$-handles and $W_2$ is obtained from $W_1$ by gluing $2$-handles.
\end{definition}

\begin{remark}
Up to $1$-deformations (defined below), if $W$ is connected, we can always assume that $W_0$ is connected, i.e. that there is a single $4$-ball. Hence from now on we shall consider only this case.
\end{remark}

Note that such an object can be encoded in a Kirby diagram as explained in \cite[Section 4.4]{gompf_stipsicz_1999}. For us, a regular Kirby diagram is a diagram of a link in $S^3$ that splits as a disjoint union of a trivial unlink diagram whose components are dotted, and a framed link. Such a datum is meant to represent the attaching maps of the handles in $S^3 = \partial D^4$. The framed knots represent embedded tori where the 2-handles are attached, while the dotted unlink represents $1$-handles as explained in \cite[Section 5.4]{gompf_stipsicz_1999}. Roughly speaking, instead of adding $1$-handles, one can equivalently push a disk bounded by an undotted compontent inside $D^4$ and excise a neighborhood of it.

A famous theorem of Cerf implies that two diffeomorphic $2$-handlebodies are related by a sequence of handle slides, creation/removal of cancelling pairs of handles, and isotopies respecting the filtration. This motivates the following definition.

\begin{definition}
\label{def2eq}
Let $1 \leq l \leq 4$. An \textit{$l$-deformation} of a $4$-dimensional $2$-handlebody is a finite sequence of
\begin{itemize}
\item isotopies of attaching maps of $i$-handles in $\partial W_{i-1}$, for $1 \leq i \leq l$,
\item handle slides of i-handles over other i-handles for $1 \leq i \leq l$
\item creations/removals of cancelling pairs of $(i/i-1)$-handles for $1 \leq i \leq l$.
\end{itemize} Two handlebodies related by an $l$-deformation are called \textit{$l$-equivalent.}
\end{definition}

\noindent It is known (see \cite[Section 1.2]{bobtcheva_piergallini_2006}) that two diffeomorphic handlebodies are $3$-equivalent, but it is not known if they are always $2$-equivalent:

\begin{conjecture}[Generalized Gompf conjecture]
\label{conj}
Two diffeomorphic $4$-dimensional $2$-handlebodies are $2$-equivalent. 
\end{conjecture}

There exist potential counterexamples to this conjecture, that is pairs of $3$-equivalent handlebodies that have not been proven not to be $2$-equivalent. Hence it is interesting to be provided with invariants of handlebodies up to $2$-equivalence.

\subsection{Invariants of handlebodies up to $2$-equivalence}

In \cite{beliakova_bobtcheva_derenzi_piergallini_2023}, a category $\operatorname{4HB}$ whose objects are $3$-dimensional $1$-handlebodies and whose morphisms are $4$-dimen\-sional $2$-handlebodies up to $2$-equivalence is defined, and proven to be equivalent as a braided category to the strict braided monoidal category freely generated by a single so-called \textit{$\BP$-Hopf algebra} $H$. Here $\BP$ is short for Bobtcheva-Piergallini, and the definition of a $\BP$-Hopf algebra object can be found at \cite[Definition 2.4.1]{beliakova_bobtcheva_derenzi_piergallini_2023}. We shall not need this result in its full extent, but we observe that any BP-Hopf algebra in a braided category $\mathcal{C}$ yields a functor from $\operatorname{4HB}$ to $\mathcal{C}$. In particular, the objects of $\operatorname{4HB}$ are in bijection with the natural numbers, and the morphisms from $0$ to $0$ are precisely the set of $4$-dimensional $2$-handlebodies up to $2$-equivalence. Thus any BP-Hopf algebra gives an invariant of $4$-dimensional $2$-handlebodies up to $2$-equivalence.

A very important feature of a unimodular ribbon Hopf algebra $H$ is that it can be ``transmuted'' into a $\operatorname{BP}$-Hopf algebra $\underline{H}$ in the category $H$-$\operatorname{mod}$ (as a special case of \cite[Prop. 7.3]{beliakova_derenzi_2023}), where $\underline{H}$ is $H$ as an $H$-module for the left adjoint action from Proposition \ref{adjointAct}. Summarizing, a unimodular ribbon Hopf algebra yields an invariant $J_H$ of $2$-handlebodies through a rather complicated process. Fortunately, we do not have to go through the details since a direct algorithm computing the invariant is given in \cite[Section 8]{beliakova_derenzi_2023}. Let us sketch the algorithm, since it will have an importance in the proofs below.

\begin{algorithm}[Bead Algorithm, \cite{beliakova_derenzi_2023}]
\label{beadalgo}
Let $W$ be a $4$-dimensional $2$-handlebody and $L$ be an associated regular Kirby diagram, then the invariant $J_H(W)$ can be computed in the following way. From the initial diagram, we get another diagram with labeled beads by \begin{enumerate}
\item deleting crossings and placing beads labeled with the components of the R-matrix, \begin{center}
\raisebox{-\height/2}{\includegraphics[scale=0.9]{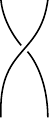}}  \hspace{2mm} $\mapsto$ \hspace{2mm}   \raisebox{-\height/2}{\includegraphics[scale=0.9]{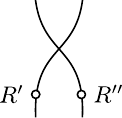}} and \hspace{2mm}  \textbf{\raisebox{-\height/2}{\includegraphics[scale=0.9]{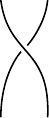}}  \hspace{2mm} $\mapsto$ \hspace{2mm}   \raisebox{-\height/2}{\includegraphics[scale=0.9]{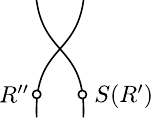}}}
\end{center}
\item deleting 1-handles and placing beads corresponding to an iterated coproduct of the cointegral, \begin{center}
\raisebox{-\height/2}{\includegraphics[scale=0.9]{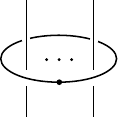}}  \hspace{2mm} $\mapsto$ \hspace{2mm}   \raisebox{-\height/2}{\includegraphics[scale=0.9]{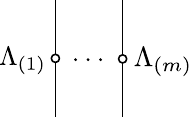}} \end{center}
\item deleting framings and placing beads labeled with copies of the pivotal element $g$ \begin{center}
\raisebox{-\height/2}{\includegraphics[scale=0.9]{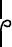}}.  \hspace{2mm} $\mapsto$ \hspace{2mm}   \raisebox{-\height/2}{\includegraphics[scale=0.9]{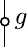}} \hspace{2mm} and \hspace{2mm}  \textbf{\raisebox{-\height/2}{\includegraphics[scale=0.9]{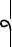}}  \hspace{2mm} $\mapsto$ \hspace{2mm}   \raisebox{-\height/2}{\includegraphics[scale=0.9]{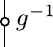}}}
\end{center}
\end{enumerate} We collect the beads together using the rules \begin{center}
\raisebox{-\height/2}{\includegraphics[scale=0.9]{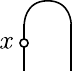}}  \hspace{2mm} $=$ \hspace{2mm}   \raisebox{-\height/2}{\includegraphics[scale=0.9]{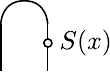}}  \hspace{2mm}  , \hspace{2mm}  \textbf{\raisebox{-\height/2}{\includegraphics[scale=0.9]{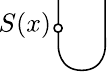}}  \hspace{2mm} $=$ \hspace{2mm}   \raisebox{-\height/2}{\includegraphics[scale=0.9]{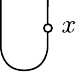}}}
\end{center}\begin{center}
and \hspace{4mm }\raisebox{-\height/2}{\includegraphics[scale=0.9]{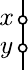}} \hspace{2mm}$=$ \hspace{2mm}   \raisebox{-\height/2}{\includegraphics[scale=0.9]{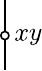}},
\end{center}  and isotope the components untill we get a trivial diagram of an unlink with a single labeled bead on every component:
\begin{center}
\includegraphics[scale=0.9]{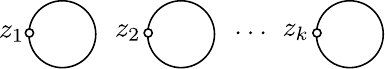}.
\end{center}

Then $$J_H(W) = \prod_{j = 1}^k \lambda(z_j),$$ where the value is computed using the Sweedler notation. That is, here, $J_H(W)$ is computed as a sum of products, that is a linear combination of terms in $\underline{H}^{\otimes k}$, and not a product.
\end{algorithm}

\subsection{Cases of degeneration of the invariant}

When the unimodular ribbon Hopf algebra is factorizable Beliakova and De Renzi proved that the invariant we obtained is actually an invariant of the boundary of the handlebody together with its signature. We call this situation \emph{degeneration to dimension $3$}.

\begin{proposition}[\cite{beliakova_derenzi_2023}, Theorem 1.2 and Section 8.2]
If the unimodular ribbon category $H$ is factorizable, then the invariant $J_H$ degenerates to dimension $3$. 
\end{proposition}

Furthermore, it is known that two diffeomorphic $4$-dimensional $2$-handlebodies become $2$-equivalent after connected summing a certain number of times with $S^2 \times D^2$. Since the invariant $J_4$ is additive with respect to the connected summing, this implies that in order to tackle Conjecture \ref{conj}, we would like the invariant $J_4(S^2 \times D^2)$ to be trivial. This is tantamount to $H$ not being cosemisimple.

\begin{proposition}[\cite{beliakova_derenzi_2023}, Appendix C]
If the unimodular ribbon category $H$ is cosemisimple, and $W_1$ and $W_2$ are diffeomorphic handlebodies, then $J_H(W_1) = J_H(W_2)$.
\end{proposition}

Finally, rougly speaking, the invariant can degenerate to dimension $2$. Indeed any $4$-dimensional $2$-handlebody $W$ deformation-retracts onto its \emph{spine}: the $\operatorname{2CW}$-complex of dimension $2$ determined by the cores of its handles. This complex, in turn, yields a presentation of the fundamental group of the $4$-dimensional $2$-handlebody. A $2$-deformation of the handlebody is reflected as an Andrews-Curtis (or 2-equivalence) move for the group presentation. We refer to \cite{bobtcheva_2023} for more details. If an invariant of $4$-dimensional $2$-handlebody cannot distinguish between handlebodies with $2$-equivalent group presentations of their fundamental groups, then we say it \emph{degenerates to their spines}. The following is a consequence of \cite{bobtcheva_2023}.
\begin{proposition}
\label{2dDegenProp}
If the unimodular ribbon Hopf algebra $H$ is triangular (i.e. has trivial ribbon element), then the invariant $J_H$ degenerates to its spine.
\end{proposition}

The generalized Andrews-Curtis conjecture states that any balanced presentation of the trivial group can be reduced to the empty presentation through a sequence of Andrews-Curtis moves. Hence an invariant that degenerates to the spines is still a useful invariant for tackling the Andrews-Curtis conjecture.

\noindent We have justified the interest in searching for non-factorizable non-cosemisimple unimodular ribbon Hopf algebras. Nevertheless we show that, for most of the non-factorizable examples $H$ we gave in the previous sections, the invariant $J_H$ is still not rich enough for the Gompf conjecture.

\subsection{Degeneration of the invariant associated to the Nenciu algebras}

In this section, we show that for any of the strongly non-factorizable Hopf algebra $$H = H(\pmb{\alpha}):=(H(\mb{m}, \mb{t}, \matfont{d}, \matfont{u}), R(\pmb{\alpha}))$$ constructed using the approach of the previous sections, the invariant $J_H$ degenerates to the spines. For the purpose of exposition, as in Example \ref{ex2}, we will call $Z^\pm_l$, $l=1, \dots, t_2$ nilpotent generators participating in the R-matrix, and $X_k$, $k=1, \dots, t_1$ the ones appearing only in the cointegral. In the notation, we emphasized the tuple $\pmb{\alpha}\in \C^l$ that parametrizes the R-matrix as in in Example \ref{ex2}.
\begin{theorem}
    Consider a Kirby diagram for $W$, with $k_1$ 1-handles and $k_2$ 2-handles. Then the invariant $J_H$ obtained using an strongly non-factorizable Hopf algebra $H(\pmb{\alpha})$ constructed as in previous sections, such as the one in Example \ref{ex2} behaves as follows.
    \begin{itemize}
        \item If $k_1 = k_2$, then $J_H(W) = J_{H(0)}(W)$, where the latter invariant degenerates to the spine.
        \item If $k_1 \neq k_2$, the invariant vanishes: $J_H(W)=0$.
    \end{itemize}
\end{theorem}
\begin{figure}[h!]
    \centering
    \includegraphics[width=0.7\linewidth]{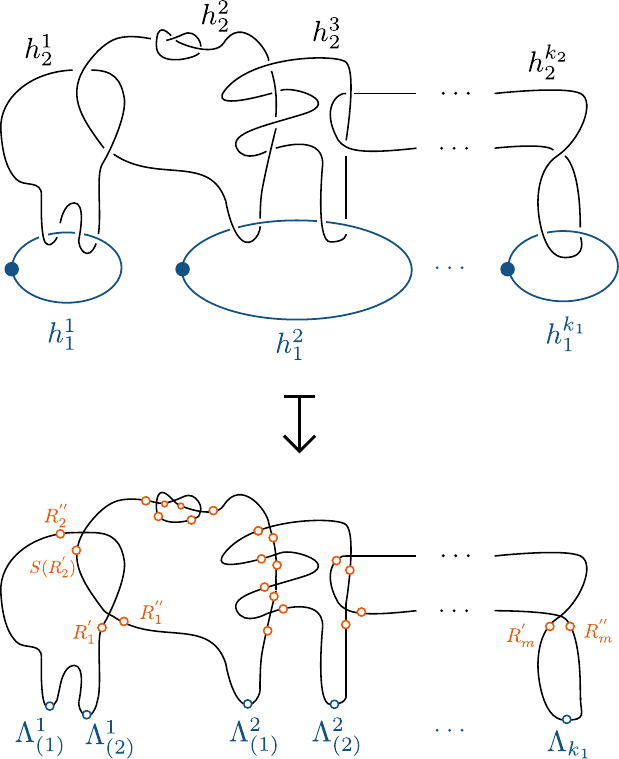}
    \caption{A Kirby diagram with $k_1$ 1-handles isotoped to the bottom, $k_2$ (possibly knotted) 2-handles and $m$  crossings, and the image under the bead algorithm, where each $2$-handle is linked to at least one $1$-handle. The beads have been colored arbitrarily and not all labeled for a better readability.}
\label{fig:2dDegen}
\end{figure}
\begin{proof}
To compute the invariant, we isotope the 1-handles to the bottom of the diagram, as in Figure \ref{fig:2dDegen}, and we use Algorithm \ref{beadalgo}. Recall that in $H$, the generators $X_k$ do not appear in the R-matrix nor the ribbon element, and that the integral $\lambda$ will vanish on a product of generators unless it contains precisely once each skew-primitive generator.
First of all, we exclude some simple cases, where the invariant vanishes regardless of $k_1$ and $k_2$. \begin{itemize}
    \item If the diagram contains an isolated 1-handle (not attached to 
    any 2-handles), then there is a contribution of $\epsilon(\Lambda)=0$ by non-semisimplicity of the algebra $H(\pmb{\alpha})$.
    \item If there is an isolated (unlinked from 1-handles and 2-handles) 2-handle in the diagram, carrying an bead with a label $P$, this gives a contribution of $\lambda(P)=0$, as $P$ is a sum of products where none of the $X_k$ generators appear. 
    \item Similarly, if there is a 2-handle linked only to another 2-handle and no 1-handles, carrying an bead with a label $P$, there is a contribution $\lambda(P)=0$, as $P$ is obtained from copies of the R-matrix and the ribbon element, corresponding to crossings and framings.
\end{itemize}  With these cases excluded, we can ensure that if $k_1 > 0$ and every 2-handle links to at least one 1-handle, as in the top of Fig. \ref{fig:2dDegen}. We apply Algorithm \ref{beadalgo} and find the bottom picture of Fig. \ref{fig:2dDegen}. Using the commutation relations, we arrive at the following value of the invariant
\[
    J_H(W) = \prod^{k_2}_{j=1} C \lambda_{j}(\Lambda^{1_j}_{(1_{1_j})}\dots\Lambda^{1_j}_{(i_{1_j})}\dots\Lambda^{m_j}_{(1_{m_j})}\Lambda^{m_j}_{(i_{m_j})} R'_a\dots S(R''_b)\dots),
\]
where the Sweedler notation hides that it is a sum over all indices. Here $C\in \C$ is a constant resulting from the commutations between the generators, and depends both on $j$ and the choice of the term in the sum. The $j$-th 2-handle links to $m_j$ 1-handles, indexed with $1_j, \dots, m_j$, and it links precisely  $i_{l_j}$ times with the $l$-th $1$-handle. Every linking carries a piece of $\Delta(\Lambda)$ in the Sweedler notation. This is followed by a product of copies of $R', R''$ and their antipodes.  
Once again, recall that $R$ is generated by grouplike generators $\mb{K}$ and nilpotent generators $Z_l^\pm$, $l=1, \dots, t_2$, while $\Lambda$ contains also $X_k$, $k = 1, \dots, t_1$, and $\lambda$ will vanish on any product of generators that does not contain precisely once each of the variable $X_k$. When taking coproducts of the integral element, the $t_1$ variables $X_k$ are distributed among the factors of the tensor product. Then, for fixed indices the term
$$\prod^{k_2}_{j=1} C \lambda_{j}(\Lambda^{1_j}_{(1_{1_j})}\dots\Lambda^{1_j}_{(i_{1_j})}\dots\Lambda^{m_j}_{(1_{m_j})}\Lambda^{m_j}_{(i_{m_j})} R'_a\dots S(R''_b)\dots)$$ contains exacatly $k_1t_1$ of the variables $X_k$. But if we want this term to be non-trivial, then there should be also $t_1$ variables $X_k$ for each $j$, in order to ``feed" the integral $\lambda$. Thus we should have $k_1t_1= k_2t_1$. Then we should have $k_1 = k_2$. As a heuristic summary, we have the following three cases.\\\\
\textbf{Case 1: ``underfed"} If $k_1 < k_2$, then the invariant vanishes, as in every term, there will be at least one $j$ such $\lambda_j$ does not have enough $X_k$ generators to ``eat'' in order to be non-zero.\\\\
\textbf{Case 2: ``overfed"} If $k_1 > k_2$, then in every term, there is a $j$ with too many of the $X_k$ variables. Thus, the invariant vanishes again.\\\\
\textbf{Case 3: ``balanced"} If $k_1=k_2$, then, at least in some terms, there is exactly enough copies of each $X^\pm_k$ to feed all integrals. Since $\Delta(\Lambda)$ produces all shuffles of $T$ (with grouplike prefactors), for at least one configuration the obstacles of cases 1 and 2 can be avoided, since we assumed $m_j>0$ for all $j$ if $k_1>0$. But now the cointegrals introduced also the correct number of copies of $Z^\pm_l$. But $$R = 1 \otimes 1 +\alpha_1 (Z_1^+\otimes L_1 Z_1^- - Z_1^-\otimes L_1 Z_1^+)+\dots,$$ so that every term involving the variable $\alpha_l$ will have too many copies of variables $Z$. Hence, setting the parameter $\alpha_l$ to $0$ does not change the final value of the invariant.

In conclusion, the invariant can be non-vanishing only in the balanced case, and then is tantamount to the invariant obtained from the triangular Hopf algebra $H(\pmb{\alpha} = \pmb{0})$.
\end{proof}


\begin{remark}
\label{feedingRem}
Here is a series of facts clarifying where the algebras $H(\pmb{\alpha})$ lie regarding the Gompf conjecture.
    \begin{enumerate}
        \item Notice that the Euler characteristic of the handlebody is given by $1- k_1 + k_2$. Hence if the Euler characteristic of $W$ is different from $1$, then $J_H(W) = 0$.
        \item Since two diffeomorphic handlebodies have the same Euler characteristic, tackling the Gompf conjecture with the invariant $J_H$ would be possible only for handlebodies of Euler characteristic $1$, corresponding to the balanced case, and then we could alternatively work with the invariant $J_{H}(\pmb{\alpha} = \pmb{0})$ associated to a triangular algebra. Hence, we would be looking at the Andrews-Curtis conjecture.
        \item The algebras are applicable indeed to the Andrews-Curtis case, but by the classification of triangular Hopf algebras of \cite{andruskiewitsch_etingof_gelaki_2002} and their correspondence to Hopf supergroups, the cases where $\pmb{\alpha} = \pmb{0}$ are incarnations of the latter.
        \item The above reasoning for the balanced case breaks down if we introduce non-diagonal relations involving some of the $Z^\pm_l$ which do not preserve the number of them in the expression, but produces a non-zero expression nevertheless. An example is $[E, F]$ in $u_q \mathfrak{sl}_2$. Then in the balanced case the generators coming from $R$ (possibly in various powers) can in principle survive up to the integrals. This motivated the construction of the semi-direct biproducts $U\ltimes H$.
    \end{enumerate}
\end{remark}

Let $U = u_q\mathfrak{sl}_2$ and $H = H(\mb{m}, t, \matfont{d}, \matfont{u})$ and let $U \ltimes H$ be their semidirect biproduct as in Definition \ref{SNFextension}. The authors observed through direct computation on many examples where $W$ is a $0$ to $0$ morphism in the category \textbf{4HB}, that the invariant $J_{U\ltimes H}$ factorizes: it satisfies $J_{U\ltimes H}(W)=J_U(W)J_H(W)$. In particular, it does not seem to carry any truly 4-dimensional information, thus we leave this as a conjecture
\begin{conjecture}
\label{failedConjecture}
    Let $U = u_q\mathfrak{sl}_2$ and $H = H(\mb{m}, t, \matfont{d}, \matfont{u})$ and let $U \ltimes H$ be their semidirect biproduct as in Definition \ref{SNFextension}. Let also $W$ be a $0$ to $0$ morphism in the category \textbf{4HB}. Then the invariant given by Algorithm \ref{beadalgo} satisfies
    \[
        J_{U\ltimes H}(W)=J_U(W)J_H(W).
    \]
\end{conjecture}

It has to be noted, however, that for higher-rank, say 1 to 0 or 1 to 1, morphisms of \textbf{4HB}, the interaction of $U$ and $H$ pieces becomes manifest, so there is a possibility the invariant can detect 4-dimensional information there.


\newpage
\appendix
\section{Proofs}
In this appendix, in a series of propositions, we prove Theorem \ref{extHopfThm}, Theorem \ref{extMainThm} as well as Propositions \ref{extProp1} and \ref{extprop2}.
\subsection{Proof of Theorem \ref{extHopfThm}}
\label{proof:extHopfThm}
\begin{proof}[Proof of Theorem \ref{extHopfThm}]
    Let $X_k \in H(\mb{m}, t, \matfont{d}, \matfont{u})$, for $k=1, \dots, t_t$. Let also $L_k:=\mb{K}^{\matfont{u}_k}$ for convenience.
    Since there is no non-trivial action on the elements $K_a \in H$, we concentrate on the generators $X_k \in H$. From the formulas of Definition \ref{smashBiproduct} we first fix for any $u \in U$, $h \in H$
    \[
        u\otimes h:=(u \otimes \1_H)(\1_U \otimes h) 
    \]
     and obtain the following algebra relations on generators
     \[
        (\1_U \otimes K_a)(K \otimes \1_H) = (K \otimes \1_H)(\1_U \otimes K_a) = (K \otimes K_a)
    \]
    \[
        (\1_U \otimes K_a)(E \otimes \1_H) = (E \otimes \1_H)(\1_U \otimes K_a) = (E \otimes K_a)
    \]
    \[
        (\1_U \otimes K_a)(F \otimes \1_H) = (F \otimes \1_H)(\1_U \otimes K_a) = (F \otimes K_a).
    \]
    \[
        (\1_U \otimes X_k)(K \otimes \1_H) = -(K \otimes \1_H)(\1_U \otimes X_k) = -(K \otimes X_k)
    \]
    \[
        (\1_U \otimes X_k)(E \otimes \1_H) = -(E \otimes \1_H)(\1_U \otimes X_k) = -(E \otimes X_k)
    \]
    \[
        (\1_U \otimes X_k)(F \otimes \1_H) = (F \otimes \1_H)(\1_U \otimes X_k) = (F \otimes X_k).
    \]
    We can easily extend these relations to arbitrary monomials $E^a F^b K^c \in U$ and $\mb{K}^\mb{v}\mb{X}^\mb{r} \in H$ for $a, b, c \in \Z_{r'}$, $\mb{v} \in \Z_\mb{m}$ and $\mb{r} \in \Z^t_2$, resulting in
    \[
        (E^a F^b K^c \otimes \1_H)(\1_U \otimes \mb{K}^\mb{v}\mb{X}^\mb{r}) = E^a F^b K^c \otimes \mb{K}^\mb{v}\mb{X}^\mb{r},
    \]
    and
    \begin{equation}
        \label{relBetweenPieces}
    \end{equation}
    \[
        (\1_U \otimes \mb{K}^\mb{v}\mb{X}^\mb{r})(E^a F^b K^c \otimes \1_H) = (-1)^{|\mb{r}|(a+c)} E^a F^b K^c \otimes \mb{K}^\mb{v}\mb{X}^\mb{r}.
    \]
    Here we used the fact the the map $\iota$ of Remark \ref{SESproposition} is an algebra map which can be checked directly. The list of the relations between the generators is an exhaustive one, as the choices of action and coaction force them to be diagonal. Moreover, with the operations defined for arbitrary monomials as above, we see that the new relations respect the monomial bases of both pieces. Since both $U$ and $H$ are fully defined by the relations on generators and their extension over the PBW- and monomial bases respectively, and this fact is respected by the new operations, we can conclude there are no relations between the pieces $U$ and $H$ not resulting from the relation of Equation \ref{relBetweenPieces}. That is, we obtain the exhaustive list of relations of $U\ltimes H$ by deriving them from the relations between the generators and applying the fact they respect the new operations of the Hopf algebra.\\
    Using the algebra structure that is now fixed, it is proved in Proposition \ref{extMonomialBasis} that $U\ltimes H$ has a monomial basis consisting of terms of the form $E^a F^b K^c \otimes \mb{K}^\mb{v} \mb{X}^\mb{r}$, so we can suppress the tensor product by defining
    \[
        E^a F^b K^c \mb{K}^\mb{v} \mb{X}^\mb{r}:=E^a F^b K^c \otimes \mb{K}^\mb{v} \mb{X}^\mb{r}
    \]
    for the entire monomial basis. 
    \\
    We now compute the coproduct of any such monomial in $U \ltimes H$ and show that $\Delta$ is an algebra map. Define the "parity" of monomials $|E^a F^b K^c|:=a+c$ and of $|\mb{K}^\mb{v}\mb{X}^\mb{r}| := |\mb{r}|$. Firstly, it is clear from Definition \ref{smashBiproduct} that 
    \[
        \Delta(E^a F^b K^c\otimes \1_H) = ((E^a F^b K^c)_{(1)}\otimes \1_H)\otimes ((E^a F^b K^c)_{(2)}\otimes \1_H), 
    \]
    so it is an algebra map for any product of monomials purely in $U$.
    Then, we also have 
    \[
        \Delta(\1_U \otimes \mb{K}^\mb{v}\mb{X}^\mb{r}) = \left(\1_U \otimes (\mb{K}^\mb{v} \mb{X}^\mb{r})_{(1)}\right)\otimes \left(K^{|(\mb{K}^\mb{v} \mb{X}^\mb{r})_{(1)}|r''} \otimes (\mb{K}^\mb{v} \mb{X}^\mb{r})_{(2)}\right)
    \]
    now let $\1_U \otimes \mb{K}^{\mb{v}_1}\mb{X}^{\mb{r}_1}$ and $\1_U \otimes \mb{K}^{\mb{v}_2}\mb{X}^{\mb{r}_2}$ be two such monomials.  Then 
    \begin{align*}
        &\Delta(\1_U \otimes \mb{K}^{\mb{v}_1}\mb{X}^{\mb{r}_1})\Delta(\1_U \otimes \mb{K}^{\mb{v}_2}\mb{X}^{\mb{r}_2})\\=&
            \left(\1_U \otimes (\mb{K}^{\mb{v}_1} \mb{X}^{\mb{r}_1})_{(1)}\right)\otimes \left(K^{|(\mb{K}^{\mb{v}_1} \mb{X}^{\mb{r}_1})_{(1)}|r''} \otimes (\mb{K}^{\mb{v}_1} \mb{X}^{\mb{r}_1})_{(2)}\right)\\& \hspace{5mm}
            \left(\1_U \otimes (\mb{K}^{\mb{v}_2} \mb{X}^{\mb{r}_2})_{(1)}\right)\otimes \left(K^{|(\mb{K}^{\mb{v}_2} \mb{X}^{\mb{r}_2})_{(1)}|r''} \otimes (\mb{K}^{\mb{v}_2} \mb{X}^{\mb{r}_2})_{(2)}\right).
    \end{align*}
    Now $K^{r''}$ commutes with all $K_a$ and $X_k$, so using the fact $\Delta_H$ is an algebra map on $H$, the expression rearranges to 
    \begin{align*}
        &\left(\1_U \otimes (\mb{K}^{\mb{v}_1} \mb{X}^{\mb{r}_1}\mb{K}^{\mb{v}_2} \mb{X}^{\mb{r}_2})_{(1)}\right)\otimes \left(K^{|(\mb{K}^{\mb{v}_1} \mb{X}^{\mb{r}_1}\mb{K}^{\mb{v}_2} \mb{X}^{\mb{r}_2})_{(1)}|r''} \otimes (\mb{K}^{\mb{v}_1} \mb{X}^{\mb{r}_1}\mb{K}^{\mb{v}_2} \mb{X}^{\mb{r}_2})_{(2)}\right)\\
        =& \Delta(\1_U \otimes \mb{K}^{\mb{v}_1} \mb{X}^{\mb{r}_1}\mb{K}^{\mb{v}_2} \mb{X}^{\mb{r}_2}).
    \end{align*}
    Thus we checked that $\Delta$ behaves as an algebra map on products of monomials coming from $H$. It is left to check the cross-terms. First consider

\begin{align*}
     &\Delta(E^a F^b K^c\otimes \1_H) \Delta(\1_U \otimes \mb{K}^\mb{v}\mb{X}^\mb{r})  \\ =&
        \left(((E^a F^b K^c)_{(1)}\otimes \1_H)\otimes ((E^a F^b K^c)_{(2)}\otimes \1_H)\right) \\& \hspace{3mm}\left((\1_U \otimes (\mb{K}^\mb{v} \mb{X}^\mb{r})_{(1)}) \otimes K^{|(\mb{K}^\mb{v} \mb{X}^\mb{r})_{(1)}|r''} \otimes (\mb{K}^\mb{v} \mb{X}^\mb{r})_{(2)}\right)\\=&
         \left((E^a F^b K^c)_{(1)}) \otimes (\mb{K}^\mb{v} \mb{X}^\mb{r})_{(1)} \right)  \otimes \left( (E^a F^b K^c)_{(2)} K^{|(\mb{K}^\mb{v} \mb{X}^\mb{r})_{(1)}|r''} \otimes (\mb{K}^\mb{v} \mb{X}^\mb{r})_{(2)}\right) \\=&
        \Delta(E^a F^b K^c\otimes \mb{K}^\mb{v}\mb{X}^\mb{r}).
\end{align*}
Then for $d, e, f \in \Z_{r'}$ using the above
\begin{align*}
    &\Delta(E^a F^b K^c\otimes \1_H) \Delta(E^e F^f K^d  \otimes \mb{K}^\mb{v}\mb{X}^\mb{r}) \\
    = &\Delta(E^a F^b K^c\otimes \1_H) \Delta(E^e F^f K^d  \otimes \1_H)\Delta(\1_U \otimes \mb{K}^\mb{v}\mb{X}^\mb{r})\\
    =&\Delta(E^a F^b K^c E^e F^f K^d\otimes \1_H)\Delta(\1_U \otimes \mb{K}^\mb{v}\mb{X}^\mb{r})\\
    =&\Delta(E^a F^b K^c E^e F^f K^d\otimes \mb{K}^\mb{v}\mb{X}^\mb{r}).
\end{align*}
Similarly, 
\begin{align*}
    &\Delta(E^a F^b K^c\otimes \mb{K}^{\mb{v}_1}\mb{X}^{\mb{r}_1}) \Delta(\1_U  \otimes \mb{K}^{\mb{v}_2}\mb{X}^{\mb{r}_2}) \\
    = &\Delta(E^a F^b K^c\otimes \1_H) \Delta(\1_U  \otimes \mb{K}^{\mb{v}_1}\mb{X}^{\mb{r}_1})\Delta(\1_U  \otimes \mb{K}^{\mb{v}_2}\mb{X}^{\mb{r}_2})\\
    = &\Delta(E^a F^b K^c\otimes \1_H) \Delta(\1_U  \otimes \mb{K}^{\mb{v}_1}\mb{X}^{\mb{r}_1}\mb{K}^{\mb{v}_2}\mb{X}^{\mb{r}_2})\\
    = &\Delta(E^a F^b K^c \otimes \mb{K}^{\mb{v}_1}\mb{X}^{\mb{r}_1}\mb{K}^{\mb{v}_2}\mb{X}^{\mb{r}_2}).
\end{align*}
Thus, it is now enough to establish
    \begin{align*}
        &\Delta(\1_U \otimes \mb{K}^\mb{v}\mb{X}^\mb{r}) \Delta(E^a F^b K^c\otimes \1_H)  \\
        =&\left((\1_U \otimes (\mb{K}^\mb{v} \mb{X}^\mb{r})_{(1)}) \otimes (K^{|(\mb{K}^\mb{v} \mb{X}^\mb{r})_{(1)}|r''} \otimes (\mb{K}^\mb{v} \mb{X}^\mb{r})_{(2)})\right) \\
        & \hspace{3mm} \left(((E^a F^b K^c)_{(1)}\otimes \1_H)\otimes ((E^a F^b K^c)_{(2)}\otimes \1_H \right) \\
        =&(-1)^{|(E^a F^b K^c)_{(1)}| |(\mb{K}^\mb{v}\mb{X}^\mb{r})_{(1)}| + |(E^a F^b K^c)_{(2)}| |(\mb{K}^\mb{v}\mb{X}^\mb{r})_{(2)}|}\\
        &\hspace{3mm} \left((E^a F^b K^c)_{(1)}) \otimes (\mb{K}^\mb{v} \mb{X}^\mb{r})_{(1)} \right)  \otimes \left( (E^a F^b K^c)_{(2)} K^{|(\mb{K}^\mb{v} \mb{X}^\mb{r})_{(1)}|r''} \otimes (\mb{K}^\mb{v} \mb{X}^\mb{r})_{(2)}\right).
    \end{align*}
    But after considering the sums hidden in the Sweedler notation we have the sign factor
    \[
        (-1)^{|(E^a F^b K^c)_{(1)}| |(\mb{K}^\mb{v}\mb{X}^\mb{r})_{(1)}| + |(E^a F^b K^c)_{(2)}| |(\mb{K}^\mb{v}\mb{X}^\mb{r})_{(2)}|} =(-1)^{|E^a F^b K^c| |\mb{K}^\mb{v}\mb{X}^\mb{r}|}
    \]
    so the result is 
    \[
        (-1)^{|E^a F^b K^c| |\mb{K}^\mb{v}\mb{X}^\mb{r}|} \Delta(E^a F^b K^c\otimes \mb{K}^\mb{v}\mb{X}^\mb{r}).
    \]
    The sign factor is exactly $(-1)^{|\mb{r}|(a+c)}$ as above. 
    Then we have the general formula
    \begin{align*}
        &\Delta(E^a F^b K^c \otimes \mb{K}^{\mb{v}_1}\mb{X}^{\mb{r}_1}) \Delta(E^e F^f K^d\otimes \mb{K}^{\mb{v}_2}\mb{X}^{\mb{r}_2})\\
        =&\Delta(E^a F^b K^c \otimes \1_H)\Delta(\1_U \otimes \mb{K}^{\mb{v}_1}\mb{X}^{\mb{r}_1})\Delta(E^e F^f K^d \otimes \1_H)\Delta(\1_U\otimes \mb{K}^{\mb{v}_2}\mb{X}^{\mb{r}_2})\\
        =&(-1)^{|E^e F^f K^g| |\mb{K}^{\mb{v}_1}\mb{X}^{\mb{r}_1}|}\Delta(E^a F^b K^c E^e F^f K^d\otimes \mb{K}^{\mb{v}_1}\mb{X}^{\mb{r}_1}\mb{K}^{\mb{v}_2}\mb{X}^{\mb{r}_2}).
    \end{align*}
    Thus, $\Delta$ is an algebra map, and $U\ltimes H$ is a bialgebra. Now, the antipode is unique if it exists, and a similar check shows it is given by the formula in Definition \ref{smashBiproduct}. Thus, we can consistently suppress the internal tensor product, and we proved $U\ltimes H$ is indeed a Hopf algebra as required.
\end{proof}
\subsection{Proof of Theorem \ref{extMainThm}}
\label{proof:extMainThm}
We prove Theorem \ref{extMainThm} in a series of propositions.
\begin{proposition}
    Let $U$ carry a left integral $\lambda_U$ and a two-sided cointegral $\Lambda_U$, and $H$ carry a two-sided integral $\lambda_H$ and a two-sided cointegral $\Lambda_H$ as in the statement of Theorem \ref{extMainThm}. Then $U\ltimes H$ admits the following structures 
    \begin{enumerate}
        \item a two-sided cointegral $\Lambda:=\Lambda_U\Lambda_H$
        \[
            \Lambda:=\frac{\{1\}^{r'-1}}{\sqrt{r''}[r'-1]!}\sum^{r'-1}_{a=0}  E^{r'-1} F^{r'-1} K^a 
             \prod^s_{a=1} \left(\sum^{m_a-1}_{b=0} K^b_a \right) \prod^t_{k=1} X_k,
        \]
        \item a left integral $\lambda\in (U\ltimes H)^*$ defined on the monomial basis of Proposition 
        \begin{multline*}
             \lambda(E^b F^c K^a\mb{K}^\mb{v}\mb{X}^\mb{r}) := \\
            \begin{cases}
                \frac{\sqrt{r''}[r'-1]!}{\{1\}^{r'-1}} \ \text{\, if \,} \mb{v} = (0, \dots, 0), \mb{r} = (1,\dots,1), a=b=c=r'-1 \\
                0 \text{\, otherwise}.
            \end{cases}
        \end{multline*}
    \end{enumerate}
\end{proposition}
\begin{proof}
We need to verify the properties of the integrals according to the Definition \ref{defunimod}, except for the two-sidedness of the integral which is already not satisfied by the $u_q \mathfrak{sl}_2$ piece as remarked in \cite{beliakova_derenzi_2023}, since the distinguished grouplike element is $K^2\neq \1$. By Theorem \ref{extHopfThm}, it suffices to do it on the monomial basis.\\
\begin{enumerate}
    \item \underline{$h\Lambda = \Lambda h=\epsilon(h)\Lambda$}:\\
    Consider $\Lambda$ multiplying any other element of the monomial basis $h \in H$. If $h = E^b F^c \mb{X}^\mb{r}$ with $b, c \in \{1, \dots, r'-1\}$ and $\mb{r}\in \Z^t_2$, that is it contains nilpotent generators, we clearly get $h \Lambda = \Lambda h = \epsilon(h)\Lambda=0$, because for the generators $E, F$ the cointegral of Proposition \ref{uqsl2prop} is two-sided, and $T$ is central (recall $T:=\mb{X}^{(1, \dots, 1)}$ is the top element of $H$). Moreover, $E^{r'-1} F^{r'-1} T$ commutes with all grouplike generators, so since $K \sum^{r'-1}_{i=0} K^i= \sum^{r'-1}_{i=0} K^i$ and $K_a \sum^{m_a-1}_{i=0} K_a^i= \sum^{m_a-1}_{i=0} K_a^i$, we find $K \Lambda = \Lambda K = \epsilon(K)\Lambda=\Lambda$ and $K_a \Lambda = \Lambda K_a = \epsilon(K_a)\Lambda=\Lambda$ for $a=1, 2, 3$. Thus, $\Lambda$ is a two-sided cointegral for $H$.\\
    \item \underline{$(\Id \otimes \lambda)\Delta(h) = \1\lambda(h)$}: \\
    We know that only the monomial $E^{r'-1} F^{r'-1} K^{r'-1}  T$ is picked up by the integral, and it is vanishing for any other element of the monomial basis. We need to verify only those monomial basis elements where all nilpotent generators appear, so of the form $ E^{r'-1} F^{r'-1} K^a \mb{K}^\mb{w} T$ for some $a\in \{0, \dots, r'\}$ and $\mb{w}\in \Z_\mb{m}$. Then the only term in the coproduct carrying the top piece in the right factor is
    \[
        \Delta( E^{r'-1} F^{r'-1}\mb{K}^\mb{w} T) =  K^{1-r'} \mb{K}^\mb{w}\otimes  E^{r'-1} F^{r'-1} K^{a}\mb{K}^\mb{w} T +\dots 
    \]
    If we now act with $\Id \otimes \lambda$, we see that $K_{a}\mb{K}^\mb{w} = K^{r'-1}$, so $a=r'-1, \mb{w}=(0, \dots, 0)$ is required for the top term to survive, all other vanish. But then we have
    \[
        \Delta( E^{r'-1} F^{r'-1} K^{r'-1}_1 T) = \1  \otimes  E^{r'-1} F^{r'-1} K^{r'-1} T +\dots
    \]
    so the left factor is proportional to $\1$ as required. Thus, after acting by $\Id \otimes \lambda$, we retrieve the axiom.\\
    \item \underline{$\lambda(\Lambda)$=1}:\\  
    Immediately satisfied by construction.
\end{enumerate}
\end{proof}
\begin{proposition}
\label{extQTprop}
Let $U = u_q \mathfrak{sl}_2$ carry the R-matrix $R_U = D\Theta$ and $H = H(\mb{m}, t, \matfont{d}, \matfont{u})$ carry the R-matrix $R_H = R_{\matfont{z}}R_{\pmb{\alpha}}$ respectively, such that
    \[
        R_{\pmb{\alpha}} = \exp\left( \sum^{t_2}_{l=0} \alpha_l(Z^+ \otimes LZ^- - Z^- \otimes LZ^+) \right).
    \]
     Then $u_q \mathfrak{sl}_2 \ltimes H(\mb{m}, t, \matfont{d}, \matfont{u})$ admits the R-matrix $R:=D R_{\matfont{z}}\Theta \Bar{R}_{\pmb{\alpha}}$, where 
    \[
        \Bar{R}_{\pmb{\alpha}}:= \exp\left( \sum^{t_2}_{l=0} \alpha_l(Z^+ \otimes \Bar{L}Z^- - Z^- \otimes \Bar{L}Z^+) \right)
    \]
    for $\Bar{L} = K^{r''} L$.
\end{proposition}
\begin{proof}
Firstly, we check whether $D R_{\matfont{z}}\Bar{R}_{\pmb{\alpha}}$ fulfills (QT1)-(QT4) of Definition \ref{QTthm}. But the only change is the appearance of the $K^{r''}$ factors which commute with the entire $H(\mb{m}, t, \matfont{d}, \matfont{u})$, and $(K^{r''})^2 = \1$. Furthermore $D$ and $R_\mb{z}$ can be collected to one R-matrix containing all grouplike generators, by including $K$ and  $\mb{K}$ in the tuple $\Bar{\mb{K}}:=(K, K_1, \dots, K_s)$, as well as $q^2$ and $\pmb{\xi}$ in $\Bar{\pmb{\xi}} := (q^2, \xi_1, \dots, \xi_s)$, and extending $\matfont{z}$ into the block matrix 
\[
    \Bar{\matfont{z}}=\begin{pmatrix}
    \matfont{z} && 0\\ 0 &&1
\end{pmatrix},
\]
so that 
\[
     R_{\Bar{\matfont{z}}} = \sum_{\mb{v}, \mb{w} \in \Z_{r'} \oplus \Z_\mb{m}} \Bar{\pmb{\xi}}^{-\mb{v}\cdot \mb{w}} \Bar{\mb{K}}^{\mb{w}} \otimes \Bar{\mb{K}}^{\mb{v}\matfont{z}}= D R_{\matfont{z}}.
\]
This is easily seen to satisfy (QT1)-(QT4) on its own, directly or by the relevant parts of Theorem \ref{QTthm}. Similarly, it satisfies (QT5) for all $X_k \in H$ as well.
Secondly, it is an easy check that 
\[
    [E\otimes F, Z^\pm \otimes \Bar{L} Z^\mp] = [E\otimes F, Z^\mp \otimes \Bar{L} Z^\pm] = 0 
\]
so that $\Theta$ commutes with $\Bar{R}_{\pmb{\alpha}}$, and the whole product $D R_{\matfont{z}}\Theta \Bar{R}_{\pmb{\alpha}}$ can be checked to fulfill (QT1)-(QT4), which can be checked directly using that $U$, $H$ are both quasitriangular.\\
Thirdly, we need to verify (QT5) for all generators. But we see that 
\[
    [\Delta(K), R_{\matfont{z}}\Bar{R}_{\pmb{\alpha}}] = [\Delta(E), R_{\matfont{z}}\Bar{R}_{\pmb{\alpha}}] = [\Delta(F), R_{\matfont{z}}\Bar{R}_{\pmb{\alpha}}] = 0
\]
so for $K, E, F$ (QT5) follows from its fulfillment in $u_q\mathfrak{sl}_2$. Also for nilpotent generators $X_k$,  we have
\[
    [\Delta(K_a), \Theta] = [\Delta(X_k), \Theta] = 0,
\]
so for $K_a$, $X_k$ follows from $DR_{\matfont{z}}\Bar{R}_{\pmb{\alpha}}$ being a valid Nenciu type R-matrix. That is $DR_{\matfont{z}}\Bar{R}_{\pmb{\alpha}}$ satisfying (QT5) for $K_a$ and $X_k$ follows from Theorem \ref{QTthm}.\\
\end{proof}
\begin{lemma}
\label{extMonodrLemma}
     For the setting of Proposition \ref{extQTprop}, the corresponding monodromy matrix is 
    \[
        M = D_{21}\Theta_{21}D\Theta (\Bar{R}_{\pmb{\alpha}})^2:= M_U (\Bar{R}_{\pmb{\alpha}})^2.
    \]
    where $M_U:= D_{21}\Theta_{21}D\Theta$.
\end{lemma}
\begin{proof}
    Since $\Theta$ and $\Bar{R}_{\pmb{\alpha}}$ commute we can rewrite it as
    \[
        M = R_{\matfont{z}, 21} D_{21} \Theta_{21} \Bar{R}_{\pmb{\alpha}, 21}R_{\matfont{z}} D \Bar{R}_{\pmb{\alpha}} \Theta 
    \]
    Now, we can check that $\Bar{R}_{\pmb{\alpha}, 21} R_{\matfont{z}}D R_{\pmb{\alpha}} = R_{\matfont{z}}D (\Bar{R}_{\pmb{\alpha}})^2$. Moreover, $R_{\matfont{z}}$ commutes with $D, \Theta$, so can be brought to $R_{\matfont{z}, 21}R_{\matfont{z}}=\1\otimes \1$ and we find 
    \[
        =D_{21} \Theta_{21} D (\Bar{R}_{\pmb{\alpha}})^2 \Theta.
    \]
    We can commute $(\Bar{R}_{\pmb{\alpha}})^2$ past $\Theta$ and retrieve $M = M_U (\Bar{R}_{\pmb{\alpha}})^2$, as expected.
\end{proof}
\begin{proposition}
\label{extRibProp}
Let, using the setting of Proposition \ref{extQTprop}, $U = u_q \mathfrak{sl}_2$ carry the ribbon structure $(D\Theta, K u_U)$ and $H = H(\mb{m}, t, \matfont{d}, \matfont{u})$ carry the ribbon structure $(R_{\matfont{z}}R_{\pmb{\alpha}}, g_H u_H)$, where $u_U$, $u_H$ are Drinfeld elements of $D\Theta$ and $R_{\matfont{z}} R_{\pmb{\alpha}}$ respectively, such that
    \[
        R_{\pmb{\alpha}} = \exp\left( \sum^{t_2}_{l=0} \alpha_l(Z_l^+ \otimes L_l Z_l^- - Z_l^- \otimes L_l Z_l^+) \right),
    \]
and 
    \[
        v_H:= \exp \left( -2 \sum^{t_2}_{k = 1} \alpha_l Z^+_l Z^-_l  \right)
    \]
    with the corresponding pivotal element satisfying  $g^2_H=\1$ and the Drinfeld element 
    \[
        u_H:= g_H \exp \left( -2 \sum^{t_2}_{k = 1} \alpha_l Z^+_l Z^-_l  \right).
    \]
Then $u_q \mathfrak{sl}_2 \ltimes H(\mb{m}, t, \matfont{d}, \matfont{u})$ carries the ribbon structure $(D R_{\matfont{z}}\Theta \Bar{R}_{\pmb{\alpha}}, v_U u_H )$, with the ribbon elment
        \begin{multline*}
            v:=v_U u_H = \\
            \frac{1-i}{\sqrt{r'}} \sum^{r'-1}_{a=0}\sum^{r''-1}_{b=0} \frac{\{-1\}^a}{[a]!} q^{-\frac{(a+3)a}{2} + 2b^2} E^aF^a  K^{-a-2b-1} g_H \exp \left( -2 \sum^{t_2}_{k = 1} \alpha_l Z^+_l Z^-_l  \right),
        \end{multline*}
        for $\Bar{L} = K^{r''} L$, with the pivotal element $g= K^{-1}$, and the Drinfeld element 
        \[
            u:=u_U u_H = \\
            \frac{1-i}{\sqrt{r'}} \sum^{r'-1}_{a=0}\sum^{r''-1}_{b=0} \frac{\{-1\}^a}{[a]!} q^{-\frac{(a+3)a}{2} + 2b^2} E^aF^a  K^{-a-2b} g_H \exp \left( -2 \sum^{t_2}_{k = 1} \alpha_l Z^+_l Z^-_l  \right).
        \] 
\end{proposition}
\begin{proof}
    It is straightforward to find the Drinfeld element of $R = D R_{\matfont{z}}\Theta \Bar{R}_{\pmb{\alpha}}$ to be $u:=u_U u_H$ since the two pieces coming from the respective algebras commute and the piece coming from $\Bar{R}_{\pmb{\alpha}}$ is in fact just $u_H$ as for $R_{\pmb{\alpha}}$, because $(K^{r''})^2=\1$. We need to find a pivotal element $g \in G(U \ltimes H)$ that would make the ribbon element $v:=g^{-1}u$ central. This is achieved by choosing $g=K$, as it verifies the criterion for $K, E, F$ by assumption, and $\{K, X_k\}=0$, which is the criterion necessary for $S^2(X_k) = K X_k K^{-1} = -X_k$. The ribbon axioms can be checked directly, using the relations in $U\ltimes H$ and the fact that both $U, H$ carry ribbon elements that we know explicitly. Let us denote them $v_H$ and $v_U$. Since $S(v_H) = v_H$ and $S(g_H) = g_H$, we have $S(u_H) = u_H$, because $[g_H, u_H]=0$. Then we have
    \[
        S(v) = S(v_U u_H) = S(u_H) S(v_U) =  u_H v_U = v_U  u_H = v.
    \]
    In a similar fashion (R2) is immediate. Finally, for (R3) we found the monodromy matrix in Lemma \ref{extMonodrLemma} above
\[
    M = D_{21}\Theta_{21}D\Theta (\Bar{R}_{\pmb{\alpha}})^2:= M_U (\Bar{R}_{\pmb{\alpha}})^2,
\]
so we have, using the commutations established so far
\[
    M\Delta(v) =   M_U (\Bar{R}_{\pmb{\alpha}})^2 \Delta(v_U) \Delta(u_H). 
\]
Now since $\Delta(g_H)^2 \Delta(K) \Delta(K^{-1})= \1\otimes \1$ we have
\[
    (\Bar{R}_{\pmb{\alpha}})^2  \Delta(u_H) =(\Bar{R}_{\pmb{\alpha}})^2 \Delta(K) \Delta(K^{-1}) \Delta(g_H)^2  \Delta(u_H).
\]
Before we proceed, let us consider the subalgebra $\C\langle K, K_a, X_k \rangle \subset U\ltimes H$, for $a=1, \dots, s$ and $k=1, .., t$. It is a Nenciu algebra. Using Theorems \ref{QTthm} and \ref{ribbonThm} it can be shown that $D R_{\matfont{z}}\Bar{R}_{\pmb{\alpha}}$ is an R-matrix with the monodromy $(\Bar{R}_{\pmb{\alpha}})^2$ for this subalgebra and $\Bar{v}_H:= K^{-1}g_H u_H$ is a compatible ribbon element with the corresponding pivotal element $g= K$. Hence we have, 
\[
    (\Bar{R}_{\pmb{\alpha}})^2  \Delta(\Bar{v}_H) = \Bar{v}_H \otimes \Bar{v}_H. 
\]
Therefore, since $[\Bar{R}_{\pmb{\alpha}}, K]=0$,
\[
     (\Bar{R}_{\pmb{\alpha}})^2  \Delta(K) \Delta(u_H) =\Delta(K)\Delta(g_H) (\Bar{R}_{\pmb{\alpha}})^2  \Delta(\Bar{v}_H)  = \Delta(K)\Delta(g_H)(\Bar{v}_H \otimes \Bar{v}_H) = (u_H\otimes u_H).
\]
Plugging it back to the main equation we find
\[
    M_U (\Bar{R}_{\pmb{\alpha}})^2 \Delta(v_U) \Delta(u_H) =  M_U  \Delta(v_U) (\Bar{R}_{\pmb{\alpha}})^2 \Delta(u_H)=\\
     =(v_U \otimes v_U)(u_H \otimes u_H) = v\otimes v. 
\]
which is what we sought.
\end{proof}
\subsection{Proofs of Propositions \ref{extProp1} and \ref{extprop2}}

\begin{proof}[Proof of Proposition \ref{extProp1}]
\label{proof:extProp1}
What follows is the proof of Theorem \ref{extMainThm} applied to the Hopf algebra of Example \ref{extexample1}. We think it informative to spell it out in some details. We split the proof into several pieces.
\begin{enumerate}
    \item We need to verify the properties of the integrals, see Preliminaries, except for the two-sidedness of the integral which is already not satisfied by the $u_q \mathfrak{sl}_2$ piece \cite{beliakova_derenzi_2023}. \\
    \underline{$(\Id \otimes \lambda)\Delta(h) = \1\lambda(h)$}: \\
    We know that only the element $E^{r'-1} F^{r'-1} K^{r'-1}  X^+ X^-$ is picked up by the integral, and it is vanishing for any other element of the monomial basis. We need to verify only those monomial basis elements where all nilpotent generators appear, so of the form $ E^{r'-1} F^{r'-1} \mb{K}^\mb{w} X^+ X^-$ for some $\mb{w}\in \Z^3_4$. Then the only term in the coproduct carrying the top piece in the right factor is
    \[
        \Delta( E^{r'-1} F^{r'-1}\mb{K}^\mb{w} X^+ X^-) =  K^{1-r'} \mb{K}^\mb{w}\otimes  E^{r'-1} F^{r'-1} \mb{K}^\mb{w} X^+ X^- +\dots 
    \]
    If we now act with $\Id \otimes \lambda$, we see that $\mb{K}^\mb{w} = K^{r'-1}$ is required for the top term to survive, all other vanish. But then we have
    \[
        \Delta( E^{r'-1} F^{r'-1} K^{r'-1}_1 X^+ X^-) = \1  \otimes  E^{r'-1} F^{r'-1} K^{r'-1} X^+ X^- +\dots
    \]
    so the left factor is proportional to $\1$ as required. Thus, after acting by $\Id \otimes \lambda$, we retrieve the axiom.\\
    \underline{$h\Lambda = \Lambda h=\epsilon(h)\Lambda$}:\\
    Now consider $\Lambda$ multiplying any other element of the monomial basis $h \in H$. If $h$ contains nilpotents, we clearly get $h \Lambda = \Lambda h = \epsilon(h)\Lambda=0$, because the cointegral of Proposition \ref{uqsl2prop} is two-sided, and $X^+ X^-$ is central. Moreover, $E^{r'-1} F^{r'-1} X^+ X^-$ commutes with all grouplike generators, so by the usual reasoning $K \Lambda = \Lambda K = \epsilon(K)\Lambda=\Lambda$ and $K_a \Lambda = \Lambda K_a = \epsilon(K_a)\Lambda=\Lambda$ for $a=1, 2$. Thus, $\Lambda$ is a two-sided cointegral for $H$.\\
    \underline{$\lambda(\Lambda)$=1}:\\  
    Satisfied by construction.

    \item It is easy to see that $R_{\matfont{z}}$ commutes with $D$ and $\Theta$. From Remark \ref{H4Rem} and Proposition \ref{ex1prop} we know that it fulfills (QT1)-(QT4), as does the product $D\Theta$. Hence the entire R-matrix does. It is left to check (QT5). Firstly, for $E$ and $F$, we see that $[E, R_{\matfont{z}}]=[F, R_{\matfont{z}}]=0$, so for these two generators (QT5) follows from it being true for $u_q \mathfrak{sl}_2$. Then, for $X^\pm$, notice that $[\Delta(X^\pm), E\otimes F]=0$ in $H\otimes H$. Thus, we determine (QT5) only for $R_{\matfont{z}}D$. The two can be written as $R_{\Bar{\matfont{z}}}$ for 
    \[
        \Bar{\matfont{z}} = \begin{pmatrix}1&&0&&0 \\ 0 && 2 && 3\\ 0 && 1 &&0 \end{pmatrix}.
    \]
    This is clearly a block matrix, so (QT5) can be checked for $K$ and $K_1, K_2$ separately. For the latter this is the same as in Example \ref{ex1}, for the former, using Remark \ref{DasRz}, we have
    \[
        \Delta^{cop}(X^\pm) D = (X^\pm \otimes \1 + (K^{r''} K_1 K_2)^{\pm 1}\otimes X^\pm) \sum_{v, w \in \Z_{r'}} (q^2)^{-wv} K^w \otimes K^{v}=
    \]
    \[
        \sum_{v, w \in \Z_{r'}} (q^2)^{-wv  \mp w r''} K^{w}  \otimes K^{v \mp r''} (X^\pm \otimes (K^{r''})^{\pm 1}) +
    \]
    \[
        + \sum_{v, w \in \Z_{r'}} (q^2)^{-wv  \mp vr''} K^{w \pm r''}  \otimes K^{v} ((K_1 K_2)^{\pm 1} \otimes X^\pm).
    \]
    The two sums can be restored to the original form by the substitutions $v \mapsto v \pm r''$ and $w \mapsto w \pm r''$, respectively. The remaining $(K_1 K_2)^{\pm1}$ factors are dealt with in the same way by $R_{\matfont{z}}$.

    \item Note that the modification to the ribbon element is given by $K^2_2$, which recall from Remark \ref{H4Rem} and Proposition \ref{ex1prop} is the Drinfeld element $u_H$ of $R_{\matfont{z}}$. Moreover, $u_H$ commutes with $v_U$, so the only thing left to check is whether $K$ remains the pivotal element. Indeed, $\{K, X^\pm\}=0$, which is consistent with the order of the antipode $S^4(X^\pm)=X^\pm$. \\
    We can also check the ribbon axioms, as follows. Let $v_U$ be the ribbon element of Proposition \ref{uqsl2prop}. For (R1), notice that $v = v_U K^2_2$. Immediately we have $\epsilon(v) = \epsilon(v_U)\epsilon(K^2_2) = 1$, since $v_U$ is the ribbon element for $u_q \mathfrak{sl}_2$. Similarly, $S(v) = S(K^2_2)S(v_U) = K_1^2 v_U = v_U K^2_2$, since $K_1^2$ is of order 2 and commutes with $v_U$. The axiom (R2) is immediate.\\
    Towards (R3), to verify $R_{21}R \Delta(v) = (v\otimes v)$, note first that the monodromy matrix is
    \[
        R_{21}R = R_{\matfont{z}, 21} D_{21} \Theta_{21} R_{\matfont{z}} D \Theta =  R_{\matfont{z}, 21} R_{\matfont{z}} D_{21} \Theta_{21}  D \Theta = M_U, 
    \]
    where $M_U = D_{21} \Theta_{21}  D \Theta$ is the monoromy matrix of $u_q \mathfrak{sl}_2$. This follows as $R_{\matfont{z}}$ is on its own triangular and commutes with the other components. Now, we have 
    \[
    M\Delta(v) = M_U\Delta(v_U)\Delta(K^2_2) = (v_U\otimes v_U)(K^2_2 \otimes K^2_2) = v_U K^2_2 \otimes v_U K^2_2,
    \]
    as required.
\end{enumerate}
\end{proof}

\begin{proof}[Proof of Proposition \ref{extprop2}]
\label{proof:extprop2}
    Most of the proof is the same as in the Proposition \ref{extProp1}. We have to establish that the new R-matrix and ribbon element fulfil the axioms of a ribbon Hopf algebra. \\
    Firstly it is easy to check that 
    \[
    [Z^\pm\otimes L Z^\mp, E\otimes F] = [Z^\pm\otimes L Z^\mp, \Delta(E)]=[Z^\pm\otimes L Z^\mp, \Delta(F)]=0.
    \]
    Thus, the $R_{\pmb{\alpha}}$ commutes with $\Theta$, while $R_{\matfont{z}}D$ is its triangular piece of the R-matrix. So we can directly again execute the (QT1)-(QT4) axioms, as well as (QT5), since the appearance $R_{\pmb{\alpha}}$ does not interrupt their holding for $E$ and $F$. In light of Theorem \ref{classicalTheorem}, we can note that the map is block on generators - are no cross-terms between $E, F$ and $Z^\pm$.\\
    To verify ribbon axioms call again $v_U$ the ribbon element of Proposition \ref{uqsl2prop}, and let $u_H = K^2_4\exp \left( - 2 \alpha Z^+ Z^-  \right)$, so that $v=v_U u_H$. Then we obviously have $\epsilon(v) = \epsilon(v_U)\epsilon(u_H) = 1$, and we check $S(v_U u_H) = S(u_Hv_U) = S(v_U)S(u_H) = v_U S(u_H)$ since both pieces commute. This is because, while $u_H$ is not the ribbon element for the subalgebra $\langle K, K_1, K_2, K_3, X^\pm, Z^\pm\rangle$, but $\exp \left( - 2 \alpha Z^+ Z^-  \right)$ is, so we compute
    \[
        S(K^2_4\exp \left( - 2 \alpha Z^+ Z^-  \right)) =  \exp \left( - 2 \alpha Z^+ Z^-  \right)K^2_4 = u_H,
    \]
    since $\exp \left( - 2 \alpha Z^+ Z^-  \right)$ is central. \\
    For the final axiom we need the monodromy matrix. Call $M_U$ the monodromy matrix of $u_q \mathfrak{sl}_2$ and $M_H$ of the subalgebra $\langle K, K_1, K_2, K_3, X^\pm, Z^\pm\rangle$, which can be easily shown to be $M_H = R^2_\alpha$. We can now write
    \[
        M = R_{\matfont{z}, 21} D_{21} \Theta_{21} R_{\alpha, 21}R_{\matfont{z}} D \Theta R_{\alpha}.
    \]
    Since $\Theta$ and $R_{\pmb{\alpha}}$ commute we can rewrite this as
    \[
        = R_{\matfont{z}, 21} D_{21} \Theta_{21} R_{\alpha, 21}R_{\matfont{z}} D R_{\alpha} \Theta 
    \]
    Now, we can check that $R_{\alpha, 21} R_{\matfont{z}}D R_{\pmb{\alpha}} = R_{\matfont{z}}D R^2_\alpha = R_{\matfont{z}}D M_H$. Moreover, $R_{\matfont{z}}$ commutes with $D, \Theta$, so can be brought to $R_{\matfont{z}, 21}R_{\matfont{z}}=\1\otimes \1$ and we find 
    \[
        =D_{21} \Theta_{21} D M_H \Theta.
    \]
    Finally, since $M_H = R^2_\alpha$ we can commute it past $\Theta$ and retrieve $M = M_U M_H$. With this we now verify
    \[
        M_U M_H\Delta(v_U u_H) = M_U M_H \Delta(u_H)\Delta(v_U) = M_U M_H \Delta(K^2_4)\Delta (\exp \left( - 2 \alpha Z^+ Z^-  \right)) \Delta(v_U).
    \]
    Since $\exp \left( - 2 \alpha Z^+ Z^-  \right)$ is central, we easily rearrange 
    \[
        M_H \Delta(K^2_4)\Delta( \exp \left( - 2 \alpha Z^+ Z^-  \right) )=\Delta(K^2_4)  M_H\Delta( \exp \left( - 2 \alpha Z^+ Z^-  \right) ) = u_H\otimes u_H,
    \]
    since the exponential part is ribbon in $H$ with respect to $M_H$. We found so far 
    \[
        M_UM_H\Delta(v_U u_H) = M_U\Delta(v_U)M_H\Delta(u_H)=M_U\Delta(v_U)(u_H \otimes u_H),
    \]
    but this rearranges to 
    \[
        M_U\Delta(v_U)(u_H \otimes u_H) = (v_U\otimes v_U) (u_H \otimes u_H) = (v_Uu_H \otimes v_Uu_H),
    \]
    as required.
\end{proof}

\section{Summary of examples}
In this appendix we collect all the examples of both constructions, with their corresponding structures. Throughout we denote 
\begin{itemize}
    \item by $\Lambda$ the two-sided cointegral 
    \item by $\lambda$ (respectively $\lambda_L$) the two-sided (respectively left) integral 
    \item by $R$ the R-matrix
    \item by $u$ the Drinfeld element
    \item by $g$ the pivotal element
    \item by $v$ the ribbon element.
\end{itemize}
\subsection{Nenciu type examples}

\begin{example}
\label{appex1}
    Let $N_1$ be the Hopf algebra generated by $K_a$, $X_k$, $Z^{\pm}_l$ for $a\in \{1, 2\}$, $j, k\in \{1, \dots, t_1\}$; $ t_1\in 4\N$ and $l, m\in \{1, \dots, t_2\}$; $t_2\in \N$, subject to the following relations
    \[
         K_a^4 = \1, \;\; K_a X_j = i X_j K_a, \;\; K_1 Z^{\pm}_l = \pm iZ^\pm_l K_1, \;\; K_2 Z^{\pm}_l =   Z^\pm_l K_2,
    \]
    \[
        Z^\pm_l X_k = \mp i X_k Z^\pm_l, \;\;\;  \{X_j, X_k\}= \{Z^{\pm}_l, Z^{\pm}_m\}=\{Z^{\pm}_l, Z^{\mp}_m\} =  0,
    \]
    where $i^2=-1$. The Hopf structure is defined by
    \[
        \epsilon(K_a) = 1, \;\; \epsilon(X_k) =\epsilon(Z^{\pm}_l)= 0,
    \]
    \[
        \Delta(K_a) = K_a \otimes K_a, \;\; \Delta(X_k) = \1 \otimes X_k + X_k \otimes K_1 K_2, \;\; \Delta(Z^{\pm}_l) = \1 \otimes Z^{\pm}_l + Z^{\pm}_l \otimes K_1^{\pm 2} K_2^{\pm 1} ,
    \]
    \[
        S(K_a)=K^{-1}_a, \;\; S(X_k) = -X_k (K_1 K_2)^{-1}, \;\; S(Z^{\pm}_l) = -Z^{\pm}_l K_1^{\mp 2} K_2^{\mp 1}.
    \]
    The Nenciu data for this Hopf algebra is
    \begin{itemize}
        \item  $\mb{m} = (4, 4)$, so $s=2$
        \item  $t = t_1+2t_2$
        \item  $\matfont{d}$ and $\matfont{u}$ are $t\times 2$ matrices
        \begin{align*}
            \matfont{d}=\begin{pmatrix}
                1&&1\\
                \multicolumn{3}{c}{$\vdots$} \\
                1&&1\\
                1&&0\\
                -1&&0\\
                \multicolumn{3}{c}{$\vdots$} \\
                 1&&0\\
                -1&&0\\
            \end{pmatrix},&&
            \matfont{u}=\begin{pmatrix}
                1&&1\\
                \multicolumn{3}{c}{$\vdots$} \\
                1&&1\\
                2&&1\\
                -2&&-1\\
                \multicolumn{3}{c}{$\vdots$} \\
                 2&&1\\
                -2&&-1\\
            \end{pmatrix}
        \end{align*}
    \end{itemize}
\end{example}
\noindent The case of minimal dimension of Example \ref{ex1} occurs for $t_1 = 4$, $t_2 = 1$, and its dimension is $4\times 4 \times 2^6 = 2^{10} = 1024$. It will be shown later this examples admits only triangular ribbon structures.
\begin{proposition}
\label{appex1prop}
    The algebra $N_1$ of Example \ref{ex1} carries
    \begin{align*}
        \Lambda &:= \left(\sum^3_{a,b=0} K^a_1 K^b_2 \right) \prod^{t_1}_{k=1}X_k \prod^{t_2}_{l=1}Z^+_l
            \prod^{t_2}_{l=1}Z^-_l\\
        \lambda(\mb{K}^\mb{v}\mb{X}^\mb{r} )&:=
            \begin{cases}
                1 \text{\, if \,} \mb{v} = (0, 0), \mb{r} = (1, 1,\dots,1) \\
                0 \text{\, otherwise}.
            \end{cases}\\
        R_{\matfont{z}} &:= \frac{1}{16}\sum_{\mb{v}, \mb{w} \in \Z^2_4} i^{-\mb{v}\mb{w}^T} \mb{K}^{\mb{w}} \otimes \mb{K}^{\mb{v} \matfont{z}},\\
        v&:=\1,
    \end{align*}
where $\matfont{z} =\begin{pmatrix} 2 && 3\\ 1 &&0 \end{pmatrix}$ and $g=K^2_1$. 
\end{proposition}

\begin{example}
    \label{appex2}
    Let $N_2$ be the Hopf algebra generated by $K_a$, $X_k^{\pm}$ and $Z^{\pm}_l$, for $a=1, 2, 3$, $j, k=1, \dots, t_1$; $ t_1\in \N$ and $l, m =1, \dots, t_2$, $t_2\in \N$, subject to the following relations
    \[
        K_a^4 = \1, \;\;  \;\; K_a X^{\pm}_k = \pm i X^{\pm}_k K_a
    \]
    \[
       K_1 Z^{\pm}_l = Z^{\pm}_l K_1, \;\; K_2 Z^{\pm}_l = - Z^{\pm}_l K_2, \;\;  K_3 Z^{\pm}_l = \pm i Z^{\pm}_l K_3
    \]
    \[
        \{X^{\pm}_j, X^{\pm}_k\} = \{X^{\pm}_j, X^{\mp}_k\} = \{Z^{\pm}_l, X^{\pm}_k\} = \{Z^{\pm}_l, X^{\mp}_k\} = \{Z^{\pm}_l, Z^{\pm}_m\}=\{Z^{\pm}_l, Z^{\mp}_m\} =  0.
    \]
    Let also $L:=K_3^2$ as a shorthand, note that $L^2 = \1$. The Hopf structure is defined by
    \[
        \epsilon(K_a) = 1, \;\; \epsilon(X^{\pm}_k)=\epsilon(Z^{\pm}_l) = 0,
    \]
    \[
        \Delta(K_a) = K_a \otimes K_a, \;\; \Delta(X^{\pm}_k) = \1 \otimes X^{\pm}_k + X^{\pm}_k \otimes (K_1 K_2)^{\pm 1}, \;\; \Delta(Z^{\pm}_l) = \1 \otimes Z^{\pm}_l + Z^{\pm}_l \otimes L
    \]
    \[
       S(K_a)=K^{-1}_a \;\; S(X^{\pm}_k) = -X^{\pm}_k (K_1 K_2)^{\mp1} \;\;  S(Z^{\pm}_l) = -Z^{\pm}_l L.
    \].
    The Nenciu data for this Hopf algebra is 
    \begin{itemize}
        \item  $\mb{m} = (4, 4, 4)$, so $s=3$
        \item  $t = 2t_1+2t_2$
        \item Then $\matfont{d}$ and $\matfont{u}$ are $t\times 3$ matrices
         \begin{align*}
            \matfont{d}=\begin{pmatrix}
                1&&1&&1\\
                -1&&-1&&-1\\
                \multicolumn{5}{c}{$\vdots$} \\
               1&&1&&1\\
                -1&&-1&&-1\\
                0&&2&&1\\
                0&&-2&&-1\\
                \multicolumn{5}{c}{$\vdots$} \\
                 0&&2&&1\\
                0&&-2&&-1\\
            \end{pmatrix},&&
            \matfont{u}=\begin{pmatrix}
               1&&1&&0\\
                -1&&-1&&0\\
                \multicolumn{5}{c}{$\vdots$} \\
               1&&1&&0\\
                -1&&-1&&0\\
                0&&0&&2\\
                0&&0&&-2\\
                \multicolumn{5}{c}{$\vdots$} \\
                 0&&0&&2\\
                0&&0&&-2\\
            \end{pmatrix}
        \end{align*}
    \end{itemize}
    
\end{example}
\noindent The case of minimal dimension of Example \ref{ex2} occurs for $t_1 = 1$, $t_2 = 1$, and its dimension is $4\times 4\times 4 \times 2^4 = 2^{10} = 1024$. 
\begin{proposition}
\label{appex2prop}
    The algebra $N_2$ of Example \ref{ex2} carries
    \begin{align*}
 \Lambda &:= \left(\sum^3_{a,b, c=0} K^a_1 K^b_2 K^c_3\right) \prod^{t_1}_{k=1}X^+_k \prod^{t_1}_{k=1}X^-_k \prod^{t_2}_{l=1}Z^+_l
            \prod^{t_2}_{l=1}Z^-_l,\\
\lambda(\mb{K}^\mb{v}\mb{X}^\mb{r} )&:=
            \begin{cases}
                1 \text{\, if \,} \mb{v} = (0, 0, 0), \mb{r} = (1, 1,\dots,1) \\
                0 \text{\, otherwise}.
            \end{cases},\\
            R_{\matfont{z}} R_{\pmb{\alpha}} &:= \frac{1}{64}\sum_{\mb{v}, \mb{w} \in \Z^3_4} i^{-\mb{v}\mb{w}^T} \mb{K}^{\mb{w}} \otimes \mb{K}^{\mb{v} \matfont{z}} \exp\left(\sum^{t_2}_{l = 1} \alpha_l( Z^+_l \otimes L Z^-_l - Z^-_l \otimes L Z^+_l)\right),\\
            v&:= \exp \left( - 2\sum^{t_2}_{l = 0} \alpha_l Z^+_l Z^-_l  \right), 
    \end{align*}
 where $\matfont{z} =\begin{pmatrix} 0 && 3 && 2 \\ 1 && 0 && 0 \\ 2 && 0 && 2  \end{pmatrix}$, $\pmb{\alpha}=(\alpha_1, \dots, \alpha_{t_2}) \in \C^{2t_2}$, and $g = L =K^2_3$.
\end{proposition}

\begin{example}
    \label{appex3}
    Let $N_3$ be the Hopf algebra generated by $K_a$, $X_k^{\pm}$, $Y^{\pm}_k$ and $Z^{\pm}_l$, for $a=1, 2, 3$, $j, k=1, \dots, t_1$; $ t_1\in \N$ and $l, m =1, \dots, t_2$; $t_2\in \N$, subject to the following relations
    \[
        K_a^4 = \1, \;\;  \;\; K_a X^{\pm}_k = \pm i X^{\pm}_k K_a, 
    \]
    \[
       K_1 Y^{\pm}_k = \pm i Y^{\pm}_k K_1, \;\; K_2 Y^{\pm}_k = Y^{\pm}_k K_2, \;\;  K_3 Y^{\pm}_k = \pm i Y^{\pm}_l K_3
    \]
    \[
       K_1 Z^{\pm}_l = Z^{\pm}_l K_1, \;\; K_2 Z^{\pm}_l = -Z^{\pm}_l K_2, \;\;  K_3 Z^{\pm}_l = \pm i Z^{\pm}_l K_3
    \]
    \[
         \{X^{\pm}_j, X^{\pm}_k\} = \{X^{\pm}_j, X^{\mp}_k\} =  \{Y^{\pm}_j, Y^{\pm}_k\} = \{Y^{\pm}_j, Y^{\mp}_k\}, \;\; X_j Y_k = i Y_k X_j
    \]
    \[
       \{Z^{\pm}_l, X^{\pm}_k\} = \{Z^{\pm}_l, X^{\mp}_k\} = \{Z^{\pm}_l, Y^{\pm}_k\} = \{Z^{\pm}_l, Y^{\mp}_k\} = \{Z^{\pm}_l, Z^{\pm}_m\}=\{Z^{\pm}_l, Z^{\mp}_m\} =  0.
    \]
    Let also $L:=K_3^2$ as a shorthand. The Hopf structure is defined by
    \[
        \epsilon(K_a) = 1, \;\; \epsilon(X^{\pm}_k)= \epsilon(Y^{\pm}_k)=\epsilon(Z^{\pm}_l) = 0, 
    \]
    \[
        \Delta(K_a) = K_a \otimes K_a, \;\; \Delta(X^{\pm}_k) = \1 \otimes X^{\pm}_k + X^{\pm}_k \otimes (K_1 K_2)^{\pm 1},
    \]
    \[
         \Delta(Y^{\pm}_k) = \1 \otimes Y^{\pm}_k + Y^{\pm}_k \otimes (K_1 K^2_2)^{\pm1} \;\; \Delta(Z^{\pm}_l) = \1 \otimes Z^{\pm}_l + Z^{\pm}_l \otimes L
    \]
    \[
       S(K_a)=K^{-1}_a \;\; S(X^{\pm}_k) = -X^{\pm}_k (K_1 K_2)^{\mp1}, \;\; S(Y^{\pm}_k) = -Y^{\pm}_k (K^2_1 K_2)^{\mp1} \;\;  S(Z^{\pm}_l) = -Z^{\pm}_l L.
    \]

    \begin{itemize}
        \item  $\mb{m} = (4, 4, 4)$, so $s=3$
        \item  $t = 4t_1+2t_2$
        \item Then $\matfont{d}$ and $\matfont{u}$ are $t\times 3$ matrices
       \begin{align*}
            \matfont{d}=\begin{pmatrix}
                1&&1&&1\\
                -1&&-1&&-1\\
                \multicolumn{5}{c}{$\vdots$} \\
               1&&1&&1\\
                -1&&-1&&-1\\
                1&&0&&1\\
                -1&&0&&-1\\
                \multicolumn{5}{c}{$\vdots$} \\
               1&&0&&1\\
                -1&&0&&-1\\
                0&&2&&1\\
                0&&-2&&-1\\
                \multicolumn{5}{c}{$\vdots$} \\
                 0&&2&&1\\
                0&&-2&&-1\\
            \end{pmatrix},&&
            \matfont{u}=\begin{pmatrix}
               1&&1&&0\\
                -1&&-1&&0\\
                \multicolumn{5}{c}{$\vdots$} \\
               1&&1&&0\\
                -1&&-1&&0\\
                2&&1&&0\\
                -2&&-1&&0\\
                \multicolumn{5}{c}{$\vdots$} \\
               2&&1&&0\\
                -2&&-1&&0\\
                0&&0&&2\\
                0&&0&&-2\\
                \multicolumn{5}{c}{$\vdots$} \\
                 0&&0&&2\\
                0&&0&&-2\\
            \end{pmatrix}.
        \end{align*}
    \end{itemize}
\end{example}
\noindent The case of minimal dimension of Example \ref{ex3} occurs for $t_1 = 1$, $t_2 = 1$, and its dimension is $4\times 4\times 4 \times 2^6 = 2^{12} = 4096$.
\begin{proposition}
\label{appex3prop}
    The algebra $N_3$ of Example \ref{ex3} carries
    \begin{align*}
            \Lambda &:= \left(\sum^3_{a,b, c=0} K^a_1 K^b_2 K^c_3 \right) \prod^{t_1}_{k=1}X^+_k \prod^{t_1}_{k=1}X^-_k 
            \prod^{t_1}_{k=1}Y^+_k \prod^{t_1}_{k=1}Y^-_k 
            \prod^{t_2}_{l=1}Z^+_l
            \prod^{t_2}_{l=1}Z^-_l,\\
    \lambda(\mb{K}^\mb{v}\mb{X}^\mb{r} )&:=
            \begin{cases}
                1 \text{\, if \,} \mb{v} = (0, 0,  0), \mb{r} = (1, 1,\dots,1) \\
                0 \text{\, otherwise}.
            \end{cases}\\ 
            R_{\matfont{z}} R_{\pmb{\alpha}} &:= \frac{1}{64}\sum_{\mb{v}, \mb{w} \in \Z^3_4} i^{-\mb{v}\mb{w}^T} \mb{K}^{\mb{w}} \otimes \mb{K}^{\mb{v} \matfont{z}} \exp\left(\sum^{t_2}_{l = 1} \alpha_l( Z^+_l \otimes L Z^-_l - Z^-_l \otimes L Z^+_l)\right),\\
            v&:= \exp \left( -2 \sum^{t_2}_{k = 1} \alpha_l Z^+_l Z^-_l  \right),
    \end{align*}
where $\matfont{z} =\begin{pmatrix} 0 && 3 && 2 \\ 1 && 0 && 0 \\ 2 && 0 && 2  \end{pmatrix}$, $\pmb{\alpha}=(\alpha_1, \dots, \alpha_{t_2}) \in \C^{2t_2}$ and $g = L = K^2_3$.  
\end{proposition}
\begin{example}
\label{appauxExample}
    Let $N_4:=H(\mb{m}, t, \matfont{d}, \matfont{u})$ be a Nenciu type algebra generated by grouplike $K_1, K_2$ and nilpotent $X^\pm$ generators, and the following relations
    \[
        K^4_1 = K^4_2 = \1, \;\;(X^\pm)^2 = 0, \;\; 
    \]
    \[
        K X^\pm K^{-1} =q^{r'}X^\pm = -X^\pm, \;\; K_1 X^{\pm} K^{-1}_2 = K_2 X^{\pm} K^{-1}_2 = \pm i X^\pm. 
    \]
    Let $L:=K^{r'/2} K_1 K_2$ for convenience. The Hopf structure is 
    \[
        \epsilon(K_1) = \epsilon(K_2)=1, \;\; \epsilon(X^\pm) = 0, 
    \]
    \[
        \Delta(K_1) = K_1 \otimes K_1, \;\; \Delta(K_2) = K_2 \otimes K_2, \;\; \Delta(X^\pm) = \1 \otimes X^\pm + X^\pm \otimes L^{\pm 1}.
    \]
    and
    \[
        S(K_1) = K_1^{-1}, \;\; S(K_2) = K_2^{-1}, \;\; S(X^\pm) = - X^\pm L^{\mp 1}.
    \]
\end{example}
\begin{proposition}
    The Hopf algebra $N_4$ of Example \ref{auxExample} carries
    \begin{align*}
            \Lambda &:= \left(\sum^3_{a,b=0} K^a_1 K^b_2 \right) X^+X^- \\
\lambda(\mb{K}^\mb{v}\mb{X}^\mb{r} )&:=
            \begin{cases}
                1 \text{\, if \,} \mb{v} = (0, 0), \mb{r} = (1, 1) \\
                0 \text{\, otherwise}.
            \end{cases}\\
            R_{\matfont{z}} &:= \frac{1}{16}\sum_{\mb{v}, \mb{w} \in \Z^2_4} i^{-\mb{v}\mb{w}^T} \mb{K}^{\mb{w}} \otimes \mb{K}^{\mb{v} \matfont{z}},\\
        v&:=\1, 
    \end{align*}
where $\matfont{z} =\begin{pmatrix} 2 && 3\\ 1 &&0 \end{pmatrix}$ and  $g = K^2_1$.
\end{proposition}

\subsection{Non-factorizable extensions of $u_q \mathfrak{sl}_2$ examples}

\begin{example}
    Let $u_q \mathfrak{sl}_2\ltimes H$, where $H = \operatorname{SF}_2$ be a Hopf algebra generated by $K, E, F$ with relations and morphisms of Definition \ref{smallQG}, as well as $K_1$, $Z^\pm$ with the following relations
    \[
        K^2_2= \1, \;\;(Z^\pm)^2 = 0, \;\; 
    \]
    \[
        K Z^\pm K^{-1} =q^{r'}Z^\pm = -Z^\pm, \;\; K_1 Z^{\pm} K^{-1}_2 =- Z^\pm 
    \]
    \[
        [K_1, K] = [K_1, E] = [K_1, F]=0
    \]
    \[
        EZ^\pm = q^{r'} Z^\pm E= - Z^\pm E, \;\; [Z^\pm, F]=0. 
    \]
    Let $L:=K^{r'/2} K_1$ for convenience. The Hopf structure is 
    \[
        \epsilon(K_1) = 1, \;\; \epsilon(Z^\pm) = 0, 
    \]
    \[
        \Delta(K_1) = K_1 \otimes K_1, \;\;  \Delta(Z^\pm) = \1 \otimes Z^\pm + Z^\pm \otimes L^{\pm 1}
    \]
    and
    \[
        S(K_1) = K_1^{-1}, \;\; S(Z^\pm) = - Z^\pm L^{\mp 1}.
    \]
\end{example}
\noindent The instance of smallest dimension with a non-trivial commutation relation $[E, F]$ occurs at $r=8$. The dimension is then $4^3\times 2^3  = 2^{8}=512$.
\begin{proposition}
    The algebra $u_q \mathfrak{sl}_2 \ltimes \operatorname{SF}_2$ carries 
    \begin{align*}
    \Lambda&:=\frac{\{1\}^{r'-1}}{\sqrt{r''}[r'-1]!}\sum^{r'-1}_{a=0}\sum^2_{b}   E^{r'-1} F^{r'-1}K^a L^b Z^+ Z^-,\\
    \lambda_L&(  E^e F^f K^a L^b (Z^+)^g (Z^-)^h):= \\
    &= \frac{\sqrt{r''}[r'-1]!}{\{1\}^{r'-1}} \delta_{a, r'-1}\delta_{b, 0}\delta_{e, r'-1}\delta_{f, r'-1}\delta_{g, 1}\delta_{h, 1},\\
R&:=R_{\matfont{z}}D\Theta\Bar{R}_{\pmb{\alpha}} \\
            \Bar{R}_{\pmb{\alpha}} &:= \exp\left( \alpha( Z^+ \otimes \Bar{L} Z^- - Z^- \otimes \Bar{L} Z^+)\right) \text{ for } \Bar{L} = K^{r'/2}L,\\
            v := \frac{1-i}{\sqrt{r'}}& \sum^{r'-1}_{a=0}\sum^{r''-1}_{b=0} \frac{\{-1\}^a}{[a]!} q^{-\frac{(a+3)a}{2} + 2b^2}  E^aF^a K^{-a-2b-1} L\exp \left( -2 \alpha Z^+ Z^-  \right),\\
            .
    \end{align*}
    where $D, \Theta$ were defined in Proposition \ref{uqsl2prop},  $R_{\matfont{z}}$ in Proposition \ref{SFprop}, and  $g=K$.
\end{proposition}
\begin{example}
\label{appextexample1}
    Let $u_q \mathfrak{sl}_2 \ltimes N_4$, be the Hopf algebra generated by $K, E, F$ with relations and morphisms of Definition \ref{smallQG}, as well as $K_1, K_2$, $X^\pm$ with the following relations
    \[
        K^4_1 = K^4_2 = \1, \;\;(X^\pm)^2 = 0, \;\; 
    \]
    \[
        K X^\pm K^{-1} =q^{r'}X^\pm = -X^\pm, \;\; K_1 X^{\pm} K^{-1}_2 = K_2 X^{\pm} K^{-1}_3 = \pm i X^\pm 
    \]
    \[
        [K_1, K_2] = [K_1, K] = [K_1, E] = [K_1, F] = [K_2, K] = [K_2, E]=[K_2, F]=0
    \]
    \[
        EX^\pm = q^{r'} X^\pm E= - X^\pm E, \;\; [X^\pm, F]=0. 
    \]
    Let $L:=K^{r'/2} K_1 K_2$ for convenience. The Hopf structure is 
    \[
        \epsilon(K_1) = \epsilon(K_2)=1, \;\; \epsilon(X^\pm) = 0, 
    \]
    \[
        \Delta(K_1) = K_1 \otimes K_1, \;\; \Delta(K_2) = K_2 \otimes K_2, \;\; \Delta(X^\pm) = \1 \otimes X^\pm + X^\pm \otimes L^{\pm 1}
    \]
    and
    \[
        S(K_1) = K_1^{-1}, \;\; S(K_2) = K_2^{-1}, \;\; S(X^\pm) = - X^\pm L^{\mp 1}.
    \]
\end{example}

\noindent The instance of smallest dimension with a non-trivial commutation relation $[E, F]$ occurs at $r=8$. The dimension is then $4^3\times 2^4\times2^2 = 2^{12}=4096$. 
\begin{proposition}
\label{appextProp1}
    The algebra $u_q \mathfrak{sl}_2 \ltimes N_4$ carries 
    \begin{align*}
  \Lambda&:=\frac{\{1\}^{r'-1}}{\sqrt{r''}[r'-1]!}\sum^{r'-1}_{a=0}\sum^4_{b, c = 0}   E^{r'-1} F^{r'-1}K^a K_1^b K_2^c X^+ X^-,\\
            \lambda_L&(  E^e F^f K^a K_1^b K_2^c (X^+)^g (X^-)^h):= \\
            &=
            \frac{\sqrt{r''}[r'-1]!}{\{1\}^{r'-1}} \delta_{a, r'-1}\delta_{b, 0}\delta_{c, 0}\delta_{e, r'-1}\delta_{f, r'-1}\delta_{g, 1}\delta_{h, 1},\\
        R&:=R_{\matfont{z}}D\Theta\\ 
            R_{\matfont{z}} &:= \frac{1}{16}\sum_{\mb{v}, \mb{w} \in \Z^2_4} i^{-\mb{v}\mb{w}^T} (K_1, K_2)^{\mb{w}} \otimes (K_1, K_2)^{\mb{v} \matfont{z}}, \\
            v &:= \frac{1-i}{\sqrt{r'}} \sum^{r'-1}_{a=0}\sum^{r''-1}_{b=0} \frac{\{-1\}^a}{[a]!} q^{-\frac{(a+3)a}{2} + 2b^2}  E^aF^a K^{-a-2b-1} K^2_2,\\
    \end{align*}
    where $D, \Theta$ were defined in Proposition \ref{uqsl2prop}, 
    $\matfont{z} =\begin{pmatrix} 2 && 3\\ 1 &&0 \end{pmatrix}$, and $g=K$.
\end{proposition}
\begin{example}
\label{appextexample2}
    Let $u_q \mathfrak{sl}_2\ltimes N_2$, be the Hopf algebra generated by $K, E, F$ with relations and morphisms of Definition \ref{smallQG}, as well as $K_1, K_2, K_3$, $X^\pm$, $Z^\pm$ with the following relations, for $a=1, 2, 3$:
    \[
      [K_a, K] = [K_a, E] = [K_a, F]=0,\;\;  K X^\pm K^{-1} = -X^\pm\;\;K_a^4 = \1, \;\;  \;\; K_a X^{\pm} = \pm i X^{\pm} K_a
    \]
    \[
         K Z^\pm K^{-1} = -Z^\pm\;\;K_1 Z^{\pm} K^{-1}_2= Z^{\pm} , \;\; K_2 Z^{\pm} K^{-1}_3= - Z^{\pm}, \;\;  K_3 Z^{\pm} K^{-1}_4= \pm i Z^{\pm}
    \]
    \[
        \{X^{\pm}, X^{\pm}\} = \{X^{\pm}, X^{\mp}\} = \{Z^{\pm}, X^{\pm}\} = \{Z^{\pm}, X^{\mp}\} = \{Z^{\pm}, Z^{\pm}\}=\{Z^{\pm}, Z^{\mp}\} =  0
    \]
    \[
        EX^\pm = - X^\pm E,\;\; EZ^\pm = - Z^\pm E, \;\; [X^\pm, F]=[Z^\pm, F]=0. 
    \]
    Let also $L:=K^{r'/2} K_3^2$ as a shorthand, note that this time $L^2 = \1$. The Hopf structure is defined by
    \[
        \epsilon(K_a) = 1, \;\; \epsilon(X^{\pm})=\epsilon(Z^{\pm}) = 0,
    \]
    \[
        \Delta(K_a) = K_a \otimes K_a, \;\; \Delta(X^{\pm}) = \1 \otimes X^{\pm} + X^{\pm} \otimes (K^{r'/2}K_1 K_2)^{\pm 1}, \;\; \Delta(Z^{\pm}) = \1 \otimes Z^{\pm} + Z^{\pm} \otimes L
    \]
    \[
       S(K_a)=K^{-1}_a \;\; S(X^{\pm}) = -X^{\pm} (K^{r'/2} K_1 K_2)^{\mp1} \;\;  S(Z^{\pm}) = -Z^{\pm} L.
    \]
\end{example}
\noindent The instance of smallest dimension occurs when $q$ is a root of unity of order $8$ and the dimension is $4^4\times 2^4\times2^2 \times 2^2= 2^{16}=65536$. 
\begin{proposition}
\label{appextprop2}
   The algebra $u_q \mathfrak{sl}_2\ltimes N_2$ carries 
    \begin{align*}
            \Lambda&:=\frac{\{1\}^{r'-1}}{\sqrt{r''}[r'-1]!}\sum^{r'-1}_{a=0}\sum^4_{b, c, d= 0}   E^{r'-1} F^{r'-1} K^a K_1^b K_2^c K_3^dX^+ X^-Z^+ Z^-,\\
             \lambda_L&(E^e F^f K^a K_1^b K_2^c K^d_4 (X^+)^g (X^-)^h  (Z^+)^i (Z^-)^j):=  \\
             &=\frac{\sqrt{r''}[r'-1]!}{\{1\}^{r'-1}} \delta_{a, r'-1}\delta_{b, 0}\delta_{c, 0}\delta_{d, 0}\delta_{e, r'-1}\delta_{f, r'-1}\delta_{g, 1}\delta_{h, 1}\delta_{a, 1}\delta_{j, 1},\\
        R&:=R_{\matfont{z}}D\Theta R_{\pmb{\alpha}},\\
        R_{\matfont{z}} &:= \frac{1}{64}\sum_{\mb{v}, \mb{w} \in \Z^3_4} i^{-\mb{v}\mb{w}^T} (K_1, K_2, K_3)^{\mb{w}} \otimes (K_1, K_2, K_3)^{\mb{v} \matfont{z}},\\ 
            R_{\mb{\alpha}} &:= \exp\left(\alpha( Z^+ \otimes L Z^- - Z^- \otimes L Z^+)\right),\\
            v := &\frac{1-i}{\sqrt{r'}} \sum^{r'-1}_{a=0}\sum^{r''-1}_{b=0} \frac{\{-1\}^a}{[a]!} q^{-\frac{(a+3)a}{2} + 2b^2} K^{-a-2b-1} K^2_3 E^aF^a\exp \left( - 2 \alpha Z^+ Z^-  \right),
    \end{align*}
    where $D, \Theta$ were defined in Proposition \ref{uqsl2prop}, $\matfont{z} =\begin{pmatrix} 0 && 3 && 2 \\ 1 && 0 && 0 \\ 2 && 0 && 2  \end{pmatrix}$, and $g=K$.
\end{proposition}

\newpage

\bibliographystyle{alpha}
\bibliography{copy}

\Addresses

\end{document}